\documentclass[aop,MSNbibl,seceqn,dvips]{arximspdf}
\usepackage{mathrsfs}
\usepackage{graphicx}

\doi{10.1214/13-AOP842} 
\volume{41}
\issue{6}
\pubyear{2013}
\firstpage{3786}
\lastpage{3878}

\makeatletter

\newcommand{\eqref}[1]{(\ref{#1})}
\newtheorem{theorem}{Theorem}
\newtheorem{lemma}{Lemma}
\newtheorem{proposition}{Proposition}
\newtheorem{corollary}{Corollary}
\newproclaim{remark}{Remark}
\newproclaim{definition}{Definition}

\newcommand{\cT}{\mathcal{T}}
\newcommand{\cL}{\mathcal{L}}

\def\r{\mathbb{R}}
\def\e{\mathbf{E}}
\def\Q{\mathbf{Q}}
\def\ppp{\mathbb{P}}
\def\eee{\mathbb{E}}
\def\P{\mathbb{P}}
\def\F{\mathscr{F}}
\def\G{\mathscr{G}}
\def\I{\mho}
\def\H{\mathscr{H}}
\def\L{\mathscr{L}}
\def\B{\mathscr{B}}
\def\bb{\mathbb{B}}
\def\gg{\mathbb{G}}
\def\D{\mathscr{D}}
\def\zz{\mathscr{Z}}
\def\cc{\mathscr{C}}
\def\n{\mathbb{N}}

\def\ee{\mathrm{e}}
\def\pounds{\textsterling}
\makeatother

\begin{document}
\begin{frontmatter}

\title{The precise tail behavior of the total progeny of a~killed
branching random walk\thanksref{T1}}
\thankstext{T1}{Supported by the French ANR 2010 Blanc 0125 (MEMEMO2)
and the Netherlands Organisation for Scientific Research (NWO).}
\runtitle{Total progeny of a killed branching random walk}

\begin{aug}
\author[A]{\fnms{Elie} \snm{A\"\i d\'ekon}\corref{}\ead[label=e1]{elie.aidekon@tue.nl}},
\author[B]{\fnms{Yueyun} \snm{Hu}\ead[label=e2]{yueyun@math.univ-paris13.fr}}
\and
\author[C]{\fnms{Olivier} \snm{Zindy}\ead[label=e3]{olivier.zindy@upmc.fr}}
\runauthor{E. A\"\i d\'ekon, Y. Hu and O. Zindy}
\affiliation{Eindhoven University of Technology, Universit\'e Paris VI
and\break Universit\'e Paris~XIII}

%
\address[A]{E. A\"\i d\'ekon\\
Department of Mathematics and\\
\quad Computer Science\\
Technische Universiteit Eindhoven\\
P.O. Box 513, 5600 MB Eindhoven\\
The Netherlands\\
\printead{e1}}
\address[B]{Y. Hu\\
D\'epartement de Math\'ematiques, UMR 7539\\
Universit\'e Paris XIII\\
93430 Villetaneuse\\
France\\
\printead{e2}}
\address[C]{O. Zindy\\
Laboratoire de Probabilit\'es UMR 7599\\
Universit\'e Paris VI\\
4 place Jussieu\\
F-75252 Paris Cedex 05\\
France\\
\printead{e3}}
\end{aug}

\received{\smonth{2} \syear{2011}}
\revised{\smonth{1} \syear{2013}}

%
\begin{abstract}
Consider a branching random walk on the real line with a killing
barrier at zero: starting from a nonnegative point, particles reproduce and
move independently, but are killed when they touch the negative
half-line. The population of the killed branching random walk
dies out almost surely in both critical and subcritical cases, where by
subcritical case we mean that the rightmost particle of the branching
random walk without killing has a negative speed, and by critical case,
when this speed is zero. We investigate the
total progeny of the killed branching random walk and give their
precise tail distribution both in the critical and subcritical cases,
which solves an
open problem of Aldous [Power laws and killed branching random walks,
\url{http://www.stat.berkeley.edu/\textasciitilde aldous/Research/OP/brw.html}].
\end{abstract}

%

\end{frontmatter}

\section{Introduction}

We consider a one-dimensional
discrete-time branching random walk $V$ on the real line $\r$. At the
beginning, there is a single particle located at the origin $0$.
Its children, who form the first generation, are positioned
according to a certain point process $\L$ on $\r$. Each of the
particles in the
first generation independently gives birth to new particles that are positioned
(with respect to their birth places) according to a point
process with the same law as $\L$; they form the second generation. And
so on. For any $n\ge1$, each particle at generation $n$ produces new
particles independently of one another and of everything up to the
$n$th generation.

Clearly, the particles of the branching random walk $V$ form a
Galton--Watson tree, which we denote by $\cT$. Call $\varnothing$ the
root. For every vertex $u\in\cT$, we denote by $|u|$ its generation
(then $|\varnothing|=0$) and by $(V(u), |u|=n)$ the positions of the
particles in the $n$th generation. Then $\L= \sum_{|u|=1} \delta_{\{
V(u)\}}$. The tree $\cT$ will encode the genealogy of our branching
random walk.

It will be more convenient to consider a branching random walk $V$
starting from an arbitrary $x\in\r$ [namely, $V(\varnothing)= x$],
whose law is denoted by $\mathbf{P}_x$ and the corresponding expectation
by $\e
_x$. For simplification, we write $\mathbf{P}\equiv\mathbf{P}_0$ and
$\e
\equiv\e
_0$. Let $\nu:= \sum_{|u|=1} 1$ be the number of particles in the first
generation, and denote by $\nu(u)$ the number of children of $u \in
\cT$.



Assume that $ \e[\nu]>1$, namely the Galton--Watson tree $\cT$ is
supercritical. Then the system survives with positive probability
$\mathbf{P}
( \cT=\infty)>0 $. Let us define the logarithmic generating
function for the branching walk
\[
\psi(t):= \log\e\biggl[ \sum_{|u|=1}
e^{t V(u)} \biggr] \in(-\infty, + \infty],\qquad t \in\r.
\]

We shall assume that $\psi$ is finite on an open interval
containing $0$ and that $\operatorname{supp} \L\cap(0, \infty) \neq\varnothing
$ [the later condition is to ensure that $V$ can visit $(0, \infty)$
with positive probability, otherwise the problem that we shall consider
becomes of a different nature]. Assume that there exists $\varrho_* >0$
such that
%
\begin{equation} \psi(\varrho_*) = \varrho_* \psi'(\varrho_*).
\label
{t*}
\end{equation}

We also assume that $\psi$ is finite on an open set
containing $[0, \varrho_*]$. The condition~(\ref{t*}) is not
restrictive: For instance, if we denote by $m^*=\operatorname{esssup}
\operatorname
{supp} \L$, then (\ref{t*}) is satisfied if either $m^*=\infty$ or
$m^*<\infty$ and $\e\sum_{|u|=1} 1_{\{ V(u)=m^*\}}< 1$; see
Jaffuel~\cite{jaffuel} for detailed discussions.

Recall that (Kingman~\cite{kingman}, Hammersley~\cite{hammersley},
Biggins~\cite{biggins}) conditioned on $\{\cT=\infty\}$,
%
\begin{equation} \lim_{n\to\infty} {1\over n} \max_{|u|=n} V(u) =
\psi'(\varrho_*)\qquad \mbox{a.s.,} \label{LGN}
\end{equation}
where $\varrho_*$ is given in (\ref{t*}). According to
$\psi
'(\varrho_*) =0$ or $\psi'(\varrho_*) <0$, we call the case critical or
subcritical. Conditioned on $\{\cT=\infty\}$, the rightmost particle in
the branching random walk without killing has a negative speed in the
subcritical case, while in the critical case it converges almost surely
to $-\infty$ in the logarithmical scale; see \cite{yzpolymer} and
\cite
{addario09} for the precise statement of the rate of almost sure convergence.

We now place a killing barrier at zero: any particle which enters
$(-\infty, 0)$ is removed and does not produce any offspring. Hence at
every generation $n\ge0$, only the particles that always stayed nonnegative
up to time $n$ survive. Denote by $\zz$ the set of all surviving
particles of the killed branching walk,
\[
\zz:= \bigl\{ u\in\cT\dvtx V(v)\ge0, \forall v\in
[\varnothing, u]\bigr\},
\]
where $[\varnothing, u] $ denotes the shortest path in the
tree $\cT$ from $u$ to the root $\varnothing$. We are interested in the
total progeny
\[
Z:= \#\zz.
\]
Then $Z < \infty, a.s., $ in both critical and subcritical cases.
David Aldous made the following conjecture:

\begin{conjecture*}[(Aldous \cite{aldous})]
\begin{longlist}[(ii)]
\item[(i)] (Critical case). If $\psi'(\varrho_*)=0$, then $\e
[Z]<\infty$
and $\e[Z \log Z]=\infty$.
\item[(ii)] (Subcritical case). If $\psi'(\varrho_*)<0$, then there
exists some constant $b>1$ such that $\mathbf{P}(Z>n) = n^{-b + o(1)}$ as
$n \to
\infty$.
\end{longlist}
\end{conjecture*}

Let us define an i.i.d. case when $\L$ is of the form $\L= \sum
_{i=1}^\nu\delta_{\{X_i\}}$ with $(X_i)_{i\ge1}$ a sequence of i.i.d.
real-valued variables, independent of $\nu$. There are several previous
works on the critical and i.i.d. cases: when $(X_i)$ are Bernoulli
random variables, Pemantle \cite{P99} obtained the precise asymptotic
of $\mathbf{P}(Z=n)$ as $n \to\infty$, where the key ingredient of
his proof
is the recursive structure of the system inherited from the Bernoulli
variables $(X_i)$. For general random variables $(X_i)$,
Addario-Berry and Broutin \cite{ab09+} recently confirmed Aldous's
conjecture (i). This was improved later by A\"id\'ekon \cite{aidekon}
who proved that for a regular tree $\cT$ (namely when $\nu$ equals some
integer), for any fixed $x\ge0$,
\begin{eqnarray*}
c_{1}R(x)e^{\varrho_* x}& \le&\liminf_{n \to\infty} n (\log n)^2
\mathbf{P}_x( Z
> n ) \le\limsup_{n \to\infty} n (\log n)^2 \mathbf{P}_x( Z > n )
\\
&\le& c_2
R(x)e^{\varrho_* x},
\end{eqnarray*}
where $c_2>c_{1}>0$ are two constants, and $R(x) $ is some renewal
function which will be defined later. For the continuous setting, the
branching Brownian motion, Maillard \cite{maillard} solved the question
by analytic tools, using link with the F-KPP equation. Berestycki et
al. \cite{berestycki} looked at the genealogy of the branching Brownian
motion with absorption in the near-critical case.

In this paper, we aim at the exact tail behavior of $Z$ both in
critical and subcritical cases and for a general point process $\L$.

Before the statement of our result, we remark that in the subcritical
case ($\psi'(\varrho_*)<0$), there are two real numbers $\varrho_-$ and
$\varrho_+$ such that $0<\varrho_-<\varrho_*< \varrho_+$ and
\[
\psi(\varrho_-)=\psi(\varrho_+)=0,
\]
[the existence of $\varrho_+$ follows from the assumption that
$\operatorname
{supp} \L\cap(0, \infty) \neq\varnothing$].

In the critical case, we suppose that
%
\begin{equation}
\e\bigl[ \nu^{1+\delta^*} \bigr] <\infty, \qquad\sup_{\theta
\in
[-\delta^*, \varrho_*+ \delta^*]}\psi(\theta) < \infty\qquad
\mbox{for some } \delta^*>0. \label{hypcrit}
\end{equation}


In the subcritical case, we suppose that
%
\begin{equation} \qquad\e\biggl[ \biggl( \sum_{|u| =1 } \bigl(1+ e^{\varrho_- V(u)}
\bigr) \biggr)^{{\varrho_+/\varrho_-}+ \delta^*} \biggr] <\infty,\qquad
\sup_{\theta\in[ -\delta^*,\varrho_+ + \delta^*]} \psi(\theta
)<\infty,
\label{hypsub}
\end{equation}
for some $\delta^*>0$. In both cases, we always assume that
there is no lattice that supports $\L$ almost surely.

Our result on the total progeny reads as follows.

\begin{theorem}[(Tail of the total progeny)]
\label{tmain}
Assume $(\ref{t*})$ and that
%
\begin{equation}\label{delta}\quad
\e\bigl[\nu^{\alpha}\bigr]<\infty\qquad \mbox{for some }\cases{
\alpha>2, &\quad $\mbox{in the critical case;}$ \vspace*{2pt}\cr
\displaystyle\alpha> 2 {\varrho_+ \over\varrho_-}, & \quad$\mbox{in the subcritical
case.}$}
\end{equation}
\begin{longlist}[(ii)]
\item[(i)] (Critical case). If $\psi'(\varrho_*)=0$ and (\ref{hypcrit}) holds,
then there exists a constant $c_{\mathrm{crit}}>0$ such that for any $x\ge0$,
\[
\mathbf{P}_x ( Z > n) \sim c_{\mathrm{crit}} R(x) e^{\varrho
_*x}
{1\over n (\log n)^2},\qquad n \to\infty,
\]
where $R(x)$ is a renewal function defined in (\ref{Rx}).

\item[(ii)] (Subcritical case). If $\psi'(\varrho_*) <0$ and (\ref{hypsub})
holds, then there exists a constant $c_{\mathrm{sub}}>0$ such that for any $x\ge
0$,
\[
\mathbf{P}_x ( Z > n) \sim c_{\mathrm{sub}} R(x) e^{\varrho_+ x} n^{-
{\varrho_+ /\varrho_-}},\qquad n \to\infty,
\]
where $R(x)$ is a renewal function defined in (\ref{Rx}).
\end{longlist}
\end{theorem}

The values of $c_{\mathrm{crit}}$ and $c_{\mathrm{sub}}$ are given in Lemma \ref{lcom}.
Let us make some remarks on the assumptions (\ref{hypcrit}) and (\ref
{hypsub}).

\begin{remark}[(I.i.d. case)] If $\L= \sum_{i=1}^\nu\delta_{\{X_i\}}$
with $(X_i)_{i\ge1}$ a sequence of i.i.d. real-valued variables,
independent of $\nu$, then (\ref{hypcrit}) holds if and only if for
some $\delta>0$, $\e[\nu^{1+\delta}]<\infty$ and $\sup_{\theta
\in
[-\delta, \varrho_*+ \delta]}\e[ e^{ \theta X_1} ] <
\infty$
while (\ref{hypsub}) holds if and only if $\e[\nu^{{\varrho_+ /
\varrho_-}+ \delta}]<\infty$ and $\sup_{\theta\in[ -\delta,
\varrho
_++ \delta]} \e[ e^{\theta X_1} ]<\infty$ for some $\delta>0$.
\end{remark}

\begin{remark} By H\"{o}lder's inequality, elementary computations show
that (\ref{hypcrit}) is equivalent to $ \e[ ( \sum_{|u|=1}
(1+e^{ \varrho_* V(u)} ) )^{1+ \delta} ] < \infty$ and
$\sup
_{\theta\in[-\delta, \varrho_*+ \delta]}\psi(\theta) < \infty$, for
some $ \delta>0$.
\end{remark}

To explain the strategy of the proof of Theorem \ref{tmain}, we
introduce at first some notation: for any vertex $u\in\cT$ and $a \in
\r$, we define
%
\begin{eqnarray} \tau_a^+(u)&:=& \inf\bigl\{0\le k\le|u| \dvtx V(u_k)> a\bigr\},
\label{deftau+}\\
\tau_a^-(u)&:= &\inf\bigl\{0\le k\le
|u| \dvtx V(u_k)< a\bigr\}, \label{deftau-}
\end{eqnarray}
with convention $\inf\varnothing:=\infty$ and for $n \ge1$ and for any
$|u |=n$, we write $\{u_0=
\varnothing, u_1,\ldots, u_n\}=[\varnothing, u] $ the shortest path
from the root $\varnothing$ to $u$ ($u_k$
is the ancestor of $k$th generation of $u$).

By using these notation, the total progeny set $\zz$ of the killed
branching random walk can be represented as follows:
\[
\zz= \bigl\{u \in\cT\dvtx \tau^-_0(u) > |u|\bigr\}.
\]

\begin{figure}

\includegraphics{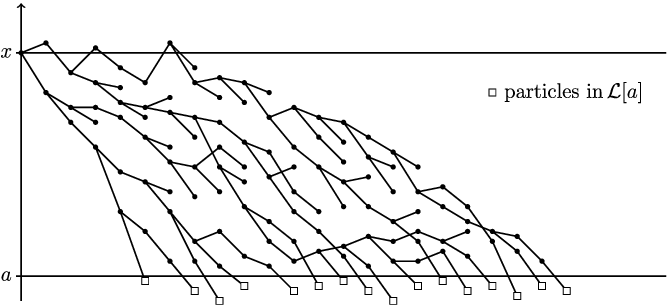}

\caption{The set $\cL[a]$.}\label{fig4}\vspace*{-3pt}
\end{figure}

For $a \le x$, we define $\cL[a]$ as the set of individuals of the
(nonkilled) branching random walk which lie below $a$ for its first
time (see Figure~\ref{fig4}):
%
\begin{equation} \label{l[a]}
\cL[a]:= \bigl\{ u\in\cT\dvtx|u|=\tau_a^-(u) \bigr\},\qquad a \le x.
\end{equation}
Since the whole system goes to $-\infty$, $\cL[a]$ is well defined. In
particular, $\cL[0]$ is the set of leaves of the progeny of the killed
branching walk.
As an application of a general fact for a wide class of graphs, we can
compare the set of leaves $\cL[0]$ with $\zz$. Then it is enough to
investigate the tail asymptotics of $ \#\cL[0]$.

To state the result for $ \#\cL[0]$, we shall need an auxiliary random
walk $S$, under a probability $\Q$, which are defined, respectively, in
(\ref{S}) and in (\ref{Q1}) with the parameter there $\varrho=
\varrho
_*$ in the critical case, and $\varrho= \varrho_+$ in the subcritical
case. We mention that under $\Q$, the random walk $S$ is recurrent in
the critical case and transient in the subcritical case. Let us also
consider the renewal function $R(x)$ associated to $S$ [see (\ref
{Rx})] and $\tau_0^-$ the first time when $S$ becomes negative; see
(\ref{tauS}). For notational simplification, let us write $\Q[\xi]$ for
the expectation of $\xi$ under~$\Q$. Then we have the following theorem.

\begin{theorem}[(Tail of the number of leaves)] \label{tleaves}
Assume (\ref{t*}).
\begin{longlist}[(ii)]
\item[(i)] (Critical case). If $\psi'(\varrho_*)=0$ and (\ref{hypcrit}) holds,
then for any $x\ge0$, we have when $n\to\infty$
\[
\mathbf{P}_x\bigl(\#\cL[0]>n\bigr) \sim c_{\mathrm{crit}}' R(x) e^{ \varrho_* x}
{1\over n
(\log n)^2},
\]
where $c_{\mathrm{crit}}':= \Q[e^{- \varrho_* S_{\tau_0^-}}] - 1
$.\eject

\item[(ii)] (Subcritical case). if $\psi'(\varrho_*)<0$ and (\ref{hypsub})
holds, then we have for any $x\ge0$ when $n\to\infty$,
\[
\mathbf{P}_x\bigl(\#\cL[0]>n\bigr) \sim c_{\mathrm{sub}}' R(x)\ee^{\varrho_+ x}
n^{-{\varrho
_+/\varrho_-}}
\]
for some constant $c'_{\mathrm{sub}} >0$.
\end{longlist}
\end{theorem}

We stress that $\Q$, $S$ and $R(\cdot) $ depend on the parameter
$\varrho=\varrho_*$ (critical case) or $\varrho=\varrho_+$ (subcritical
case). If $\sum_{|u|=1} (1+ e^{\varrho_- V(u)})$ has some larger
moments, then we can give, as in the critical case (i), a probabilistic
interpretation of the constant $c'_{\mathrm{sub}}$ in the subcritical case.

\begin{lemma}\label{Lconst} Under (\ref{t*}) with $\psi'(\varrho_*)<0$
and (\ref{hypsub}). Let us assume furthermore that
%
\begin{equation}\label{hyp-extra} \e\biggl[ \biggl( \sum_{|u|=1} \bigl(1+
e^{\varrho_- V(u)}\bigr) \biggr)^{ { \varrho_+/\varrho_-} +1+\delta
}\biggr] <\infty\qquad \mbox{for some } \delta>0,
\end{equation}
then
\[
c'_{\mathrm{sub}}= c_{\varrho_-} \bigl(c^*_{\mathrm{sub}}\bigr)^{\varrho_+/\varrho_-} \Q
\bigl(\tau
_0^- =\infty\bigr),
\]
where $c_{\varrho_-}$ and $c^*_{\mathrm{sub}}$ are given,
respectively,
by (\ref{liuvarrho-}) and Lemma \ref{Lconjecture2c} [$\Q(\tau_0^-
=\infty)>0$ since the random walk $S$ under $\Q$ drifts to $\infty$].
\end{lemma}

The next lemma establishes the relation between $\#\cL[0]$ and the
total progeny $Z=\#\zz$. Recall that $ \e[\nu]>1$.

\begin{lemma}\label{lcom}
Assume (\ref{delta}). Then Theorem \ref{tleaves} implies Theorem
\ref
{tmain} with:
\begin{longlist}[(ii)]
\item[(i)] in the critical case: $c_{\mathrm{crit}}=( \e[\nu]-1)^{-1} c_{\mathrm{crit}}'$;
\item[(ii)] in the subcritical case: $c_{\mathrm{sub}} =
( \e[\nu]-1)^{- {\varrho_+ /\varrho_-}} c_{\mathrm{sub}}'$.
\end{longlist}
\end{lemma}

\begin{figure}

\includegraphics{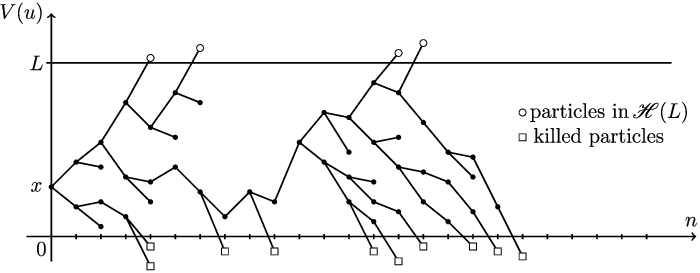}

\caption{The set $\H(L)$.} \label{fig3}
\end{figure}

The above lemma will be proven in Section~\ref{scomparison}, and the
rest of this paper is devoted to the proof of Theorem \ref{tleaves}.
To this end, we shall investigate the maximum of the killed branching
random walk and its progeny. Define for any $L >0$,
%
\begin{equation}
H(L):= \sum_u 1_{\{ \tau_0^-(u) >
\tau_L^+(u)= |u|\}} = \#\H(L),\qquad L>0, \label{hl}
\end{equation}
where
%
\begin{equation} \label{defhl2}
\H(L):= \bigl\{u \in\cT\dvtx \tau_0^-(u) > \tau_L^+(u)= |u|\bigr\}
\end{equation}
denotes the set of particles of the branching random walk on $[0,L]$
with two killing barriers which were absorbed at level $L$ [then $\H(L)
\subset\zz$]. Finally, we define
%
\begin{equation} \label{zl}Z[0, L]:= \sum_u 1_{\{ \tau_0^-(u) = |u| <
\tau_L^+(u)\}},\qquad L>0,
\end{equation}
the number of particles (leaves) which touch $0$ before $L$; see
Figure~\ref{fig3}.\eject

The following result may have independent interest: The first two
parts give a precise estimate on the probability that a level $t$ is
reached by the killed branching random walk. In the third part,
conditioning on the event that the level $t$ is reached, we establish
the convergence in distribution of the overshoots at level $t$ seen as
a random point process.

\begin{theorem}\label{thmyaglom-critiq} Assume $(\ref{t*})$.
\begin{longlist}[(iii)]
\item[(i)]
Assuming $\psi'(\varrho_*)=0$ (critical case) and (\ref{hypcrit}),
we have
\[
\mathbf{P}_x\bigl(H(t)>0\bigr) \sim{\Q[\Re^{-1}] \over C_R } R(x) e^{\varrho
_* x}
\frac{e^{-\varrho_* t}}{t},\qquad t \to\infty,
\]
where $\Q$ is defined in (\ref{Q1}), the random variable $\Re$ is
given in (\ref{eqref-R}) with $\varrho=\varrho_*$ and $C_R>0$ is a
constant given in (\ref{defcr}).

\item[(ii)]
Assuming $\psi'(\varrho_*)<0$ (subcritical case) and (\ref{hypsub}),
we have
\[
\mathbf{P}_x\bigl(H(t)>0\bigr) \sim{\Q[\Re^{-1}] \over C_R } R(x) e^{\varrho_+
x}
e^{-\varrho_+ t},\qquad t \to\infty,
\]
where $\Q$ is defined in (\ref{Q1}), the random variable $\Re$ is
given in (\ref{eqref-R}) with $\varrho=\varrho_+$ and $C_R>0$ is a
constant given in (\ref{defcr}).

\item[(iii)] In both cases and under $ \mathbf{P}_x (\cdot \vert H(t)
>0)$, the point process $\mu_t:= \sum_{ u \in\H(t)} \delta_{\{
V(u)-t\}}$
converges in distribution toward a point process $\widehat\mu_\infty$
on $(0, \infty)$, where $\widehat\mu_\infty$ is distributed as $\mu
_\infty$ under the probability measure ${\Re^{-1}\over\Q[\Re
^{-1} ]} \cdot\Q, $ with $\mu_\infty$ defined in (\ref
{eqref-muinfty}).
\end{longlist}
\end{theorem}


%
%
%
%
%
%
%
%
%
%
%
%
%
%
%
%
%
%
%
%
%
%
%
%
%
%
%
%
%
%
%
%
%
%
%
%
%
%
%
%
%
%
%
%
%
%
%
%
%
%
%
%
%
%
%
%
%
%
%
%
%
%
%
%
%
%
%
%
%
%
%
%
%
%
%
%
%
%
%
%
%
%
%
%
%
%
%
%
%
%
%
%
%
%
%
%
%
%
%
%
%
%
%
%
%
%
%
%
%
%
%
%
%
%

The Yaglom-type result of Theorem \ref{thmyaglom-critiq} plays a
crucial role in the proof of Theorem \ref{tleaves}. Loosely speaking,
to make the total progeny $Z$ (or the set of leaves $\cL[0]$) as large
as possible, the branching walk will reach some level $L$ as high as
possible, and the main contribution to $\#\cL[0]$ comes from the
descendants of those particles which have hit $L$. We control the
contribution from the other particles by computing the moments of $Z[0,
L]$ which are the most technical parts in the proof of Theorem \ref{tleaves}.

In the computations of moments of $Z[0, L]$, we have to distinguish the
contributions of \textit{good} particles from \textit{bad} particles.
By good
particle, we mean that its children do not make extraordinary jumps
[and the number of its children is not too big; see \eqref{defbad} and
\eqref{defbad2} for the precise definitions]. 
To describe separately the numbers of good and bad particles in
$Z[0,L]$, we shall modify the Yaglom-type result Theorem \ref
{thmyaglom-critiq}(iii) as follows.

Denote by $\Omega_f$ the set of $\sigma$-finite measures on $\r$. For
any individual $u\neq \varnothing$, let ${\buildrel\leftarrow\over
u}$ be the parent of $u$ and define
\[
\Delta V(u):= V(u) - V({\buildrel\leftarrow\over u}).
\]

Let us fix a measurable function $\B\dvtx\Omega_f\to\r_+$ and
write by a slight abuse of notation
\[
\B(u)\equiv\B\biggl( \sum_{ {\buildrel\leftarrow\over v}
={\buildrel
\leftarrow\over u }, v \not=u } \delta_{\{ \Delta V(v)\}}\biggr)\qquad
\forall u \in\cT\setminus\{\varnothing\},
\]
and $\B(u)=0$ if $u$ does not have any brothers. We assume
some integrability: there exists some $\delta_1>0$ such that
%
\begin{equation} \label{hypb} \e\biggl[ \sum_{|u|=1} \bigl( 1+ 1_{\{
\varrho
=\varrho_*\}} \bigl\vert V(u)\bigr\vert\bigr) e^{ \varrho V(u)} \B(u)^{\delta_1}
\biggr]< \infty,
\end{equation}
where $\varrho=\varrho_*$ if $\psi'(\varrho_*)=0$ and
$\varrho
=\varrho_+$ if $\psi'(\varrho_*)<0$. For the function $\B$
appearing in
this paper, for instance, $\B(\theta):= ( {1\over\lambda}\int
(1+e^{\varrho x}) \theta({d}x))^2$ in the critical case and
$\B(\theta
):= ( {1\over\lambda} \int\theta({d}x) (1+ e^{\varrho_-
x}))^{1/\varrho_-}$ in the subcritical case (see Sections~\ref{prcritical} and \ref{prsubcritical} where the constant $\lambda$ is
introduced) for $\theta\in\Omega_f$, (\ref{hypb}) will always be a
consequence of (\ref{hypcrit}) or (\ref{hypsub}) by taking a
sufficiently small~$\delta_1$.

Define for $u \in\cT$,
%
\begin{equation}\label{blu} \beta_L(u):= \inf\bigl\{1\le j \le|u|\dvtx
\B(u_j)
> e^{ L- V(u_{j-1})} \bigr\},\qquad L>0,
\end{equation}
with the convention that $\inf\varnothing=\infty$. We consider
\[
\H_\B(L):= \bigl\{ u \in\cT\dvtx \tau_0^-(u) >
\tau
_L^+(u)= |u|, \beta_L(u)=\infty\bigr\}.
\]

In other words, $\H_\B(L)$ only contains those particles $u$
in $\H(L)$ such that $ \B(u_j), j \le|u|, $ are not very large.
Obviously, $\H_\B\equiv\H$ if $\B=0$. We get an extension of Theorem
\ref{thmyaglom-critiq}(iii) as follows:

\begin{proposition}\label{pyaglom} Assume (\ref{hypb}) and the
hypothesis of Theorem \ref{thmyaglom-critiq}. Under $ \mathbf{P}_x
(\cdot
\vert H(t) >0)$, the point process $\mu_{\B,t}:=\sum_{ u \in\H
_\B
(t)} \delta_{\{V(u)-t\}}$ converges in distribution toward a point
process $\widehat\mu_{\B,\infty}$ on $(0, \infty)$, where
$\widehat
\mu_{\B,\infty}$ is distributed as $\mu_{\B,\infty}$ under the
probability measure ${\Re^{-1}\over\Q[\Re^{-1} ]}
\cdot
\Q, $ with $\mu_{\B,\infty}$ defined in (\ref{eqref-muinfty2}).
\end{proposition}

To prove Theorems \ref{tleaves}, \ref{thmyaglom-critiq} and
Proposition \ref{pyaglom}, we shall develop a spinal decomposition for
the killed branching random walk up to some stopping lines. Viewed from
the stopping lines, the branching walk on the spine behaves as a
two-dimensional Markov chain: The first coordinate is a real-valued
random walk (sometimes conditioned to stay positive) until some first
passage times, and the second coordinate takes values in the space of
point measures, whose laws we shall describe through a family of Palm
measures. As the parameter of the stopping lines goes to infinity, we
shall also need some accurate estimates on the real-valued random walk
and establish a convergence in law for the time-reversal random walk,
in both transient and recurrent cases.

The rest of this paper is organized as follows:

\begin{itemize}
\item Section~\ref{ssketch}: we explain the main ideas in the proofs
of Theorems \ref{tleaves} and \ref{thmyaglom-critiq}.
\item Section~\ref{scomparison}: we prove Lemma \ref{lcom}. Then the
rest of this paper is devoted to the proofs of Theorems \ref{tleaves},
\ref{thmyaglom-critiq}, Lemma \ref{Lconst} and Proposition \ref{pyaglom}.
\item Section~\ref{1dRW}: we collect several preliminary results on the
one-dimensional real-valued random walk, both in recurrent and
transient cases; in particular, we establish a result of convergence in
law for a time reversal random walk. The proofs of these results are
postponed in Section~\ref{proofsRW}.
\item Section~\ref{secspinal-decompo}: we develop the spinal
decompositions for the killed and nonkilled branching random walks,
which are the main theoretical tools in the proofs.
\item Section~\ref{secyaglom}: by admitting three technical lemmas
(whose proofs are postponed in Section~\ref{proofsRW}), we prove
Theorem \ref{thmyaglom-critiq} and Proposition \ref{pyaglom}.
\item Sections~\ref{prcritical} and \ref{prsubcritical}: based on
Theorem \ref{thmyaglom-critiq} and Proposition \ref{pyaglom}, we
prove Theorem \ref{tleaves} in the critical and subcritical cases,
respectively. We also prove Lemma \ref{Lconst} in this section.
\item Section~\ref{proofsRW} contains the proofs of the technical
lemmas stated in Sections~\ref{1dRW} and \ref{secyaglom}.
\end{itemize}

Throughout this paper, we adopt the following notation: For a point
process $\Theta= \sum_{i=1}^m \delta_{\{x_i\}}$, we
write $\langle f, \Theta\rangle= \sum_{i=1}^m f(x_i)$. Unless stated
otherwise, we denote by $c$ or $ c'$ (possibly with some subscript)
some unimportant positive constants whose values may change from one
paragraph to another, and by $f (t) \sim g(t) $ as $ t \to t_0 \in[0,
\infty]$ if $\lim_{t \to t_0} {f(t) \over g(t)}=1$; We also write $\e[
X, A]\equiv\e[X 1_A]$ when $A$ is an event and $\e[X]^k = \e[ X^k]
\neq ( \e[X])^k$ when $X$ does not have a short expression.

\section{\texorpdfstring{Heuristics in the proofs of Theorems \protect\ref{tleaves} and \protect\ref{thmyaglom-critiq}}
{Heuristics in the proofs of Theorems 2 and 3}} \label{ssketch}

For brevity, we consider $x=0$ in both Theorems \ref{tleaves} and \ref{thmyaglom-critiq}.



\subsection{\texorpdfstring{Sketch of the proof of Theorem \protect\ref{tleaves}}
{Sketch of the proof of Theorem 2}}

To make $\#\cL[0] \ge n$ very large, the killed branching random walk
needs to hit a high level, say $L$. Recalling \eqref{hl}, \eqref
{defhl2} and \eqref{zl}, we have
\[
\#\cL[0]= Z[0, L]+ \sum_u 1_{\{ \tau_0^-(u) = |u| > \tau_L^+(u)\}}.
\]
Observe that in the above sum over $u$ (if such $u$ exists), the
particle $u$ must be a descendant of some $v \in\H(L)$.\vspace*{1pt} Let us order
the set of particles in $\H(L)$ (possibly empty) in an arbitrary way:
$\H(L)= \{v^{(i)}, 1\le i \le H(L)\}$. Denote by
$\#\cL^{(i)}[0]$ the number of descendants of $v^{(i)} $ which
are absorbed at~$0$ (namely the number of the leaves of the subtree
rooted at $v^{(i)}$). Then we have
\[
\#\cL[0]= Z[0,L]+ \sum_{i=1}^{H(L)} \#\cL^{(i)}[0].
\]

The proof of Theorem \ref{tleaves} is divided into three
main steps:

\begin{longlist}[(1)]
\item[(1)] With a suitable choice of $L=L(n)$, we show that $Z[0,L]$ is
negligible in the event $\{\#\cL[0]>n\}$, which will be a consequence
of Lemmas \ref{lbad} and \ref{L2moment} in the critical case and of
Lemmas \ref{lmomentsubZ}, \ref{lsub-bad} in the subcritical case. The
proof of this fact relies on the computations of the moments of
$Z[0,L]$ by distinguishing the \textit{good} and the \textit{bad} particles.
A particle is either \textit{good} or \textit{bad}; see \eqref
{defbad} and
\eqref{defbad2} for the precise definitions in both critical and
subcritical cases. Roughly saying, a particle is called \textit{bad} if
one of its ancestors makes an extraordinary large jump. The \textit{bad}
particles are few and it is enough to compute the first moment to
control their contributions to $\#\cL[0]$, whereas for the \textit{good}
particles we need to control their higher moments. The computations of
moments are technical and follow from the change of probabilities
(spinal decomposition) and the estimates for random walks presented in
Section~\ref{1dRW}.

Let us denote by $ Y_1 \approx Y_2$ when $\ppp(Y_1 > n) \sim\ppp(Y_2
>n)$ as $n \to\infty$, where the probability $\ppp$ may be $\mathbf{P}$
or $\Q
$ whose choice will be fixed in the proof according to the random
variable $Y_1$ or $Y_2$. It follows that
\[
\#\cL[0]= Z[0,L] + \sum_{i=1}^{H(L)} \#\cL^{(i)}[0] \approx\sum
_{i=1}^{H(L)} \#\cL^{(i)}[0].
\]

Let $H_g(L) $ be the number of some subset $\H_g(L)$ of \textit{good}
particles in $\H(L)$; see~\eqref{defhgl} and \eqref{defhgl2}. Denote by
$\{u^{(j)}, 1\le j \le H_g(L)\}$ the set $\H_g(L)$. For notational
brevity we continue to use the notation $\#\cL^{(j)}[0]$ for the number
of leaves of the subtree rooted at $u^{(j)}$.
Since \textit{bad} particles in $\H(L)$ are negligible as those in $Z[0,
L]$, we have
\[
\#\cL[0] \approx\sum_{i=1}^{H(L)} \#\cL^{(i)}[0] \approx\sum
_{j=1}^{H_g(L)} \#\cL^{(j)}[0].
\]

\item[(2)] Let us consider now the critical case. By a linear
transform we may assume that $\rho^*=1$. By Nerman \cite{nerman81}, on
$\{V(u^{(j)})=y\}$, $\#\cL^{(j)}[0]$ is of order $e^y\over y$ as $y
\to
\infty$. More precisely, if we denote by $B^{(j)}:= e^{- V(u^{(j)})}
V(u^{(j)}) \#\cL^{(j)}[0]$, then under $\mathbf{P}$, conditioning on
$\{
V(u^{(j)}), 1\le j \le H(L)\}$ and letting $n \to\infty$ [hence
$L=L(n) \to\infty$], $B^{(j)}$ converges in law to $ c^* \partial
W^{(j)}_\infty$ where $c^*$ is some positive constant and $\partial
W_\infty^{(j)}, j\ge1$, are independent copies of $\partial W_\infty$,
and $\partial W_\infty$ is the limit of the so-called derivative
martingale in the critical case. Therefore,
\[
\#\cL[0] \approx\sum_{j=1}^{H_g(L)} { e^{ V(u^{(j)})} \over
V(u^{(j)})} B^{(j)} \approx c^* \sum_{j=1}^{H_g(L)} { e^{ V(u^{(j)})}
\over V(u^{(j)})} \partial W_\infty^{(j)}.
\]

Remark that $V(u^{(j)}) \sim L$. By Proposition \ref{pyaglom}, a
modified version of Theorem \ref{thmyaglom-critiq}, under $\mathbf
{P}$ and
conditioning on $\{ H(L) >0\}$,
\[
\sum_{j=1}^{H_g(L)} { e^{ V(u^{(j)})} \over V(u^{(j)})} \partial
W_\infty^{(j)} \sim{ e^L\over L} \sum_{j=1}^{H_g(L)} e^{ V(u^{(j)})-
L} \partial W_\infty^{(j)} \approx{e^{L} \over L} \sum_{i=1}^{
\hat
\zeta} e^{x_i} \partial W_\infty^{(i)},
\]
where $\sum_{i=1}^{ \hat\zeta} \delta_{\{x_i\}} $ denotes some point
process on $(0, \infty)$ defined under $\Q$ and independent of
$(\partial W_\infty^{(i)}, i\ge1)$ which are i.i.d. and are
distributed as $\partial W_\infty$ under $\mathbf{P}$. Then by letting
$L= \log
n + \log\log n - A$ with a large $A$,
\[
\mathbf{P}\bigl( \# \cL[0] > n \bigr) \sim\Q\Biggl( \sum_{i=1}^{ \hat\zeta} e^{x_i}
\partial W_\infty^{(i)} > { e^A\over c^*} \Biggr) \mathbf{P}\bigl( H(L)
>0\bigr)
.
\]

By Theorem \ref{thmyaglom-critiq}(i), there exists some constant
$\hat c>0$ such that $\mathbf{P}( H(L) >0) \sim\hat c {e^{-L}
\over L}
\sim\hat c { e^{A} \over n ( \log n)^2}$ as $n \to\infty$. Then as
$n \to\infty$,
%
\begin{equation}\label{sk1} n (\log n)^2 \mathbf{P}\bigl( \# \cL[0] > n \bigr)
\sim
\hat c e^A \Q\Biggl( \sum_{i=1}^{ \hat\zeta} e^{x_i} \partial
W_\infty^{(i)} > { e^A\over c^*} \Biggr).
\end{equation}

Rigorously speaking, instead of the above equivalence as $n \to\infty
$ in \eqref{sk1}, we have to deal with $\limsup_{ n \to\infty} $ and
$\liminf_{ n \to\infty}$ on the left-hand side of \eqref{sk1}, and we
get an upper bound and a lower bound on the right-hand side of \eqref
{sk1} with an extra term $o_A(1)$ which goes to $0$ when $A \to\infty
$. It turns out that $ \sum_{i=1}^{ \hat\zeta} e^{x_i} \partial
W_\infty^{(i)} $ has a Cauchy-law tail; see Lemma \ref{Ltail9}. (The
point process $\hat\zeta$ may depend on some parameter after the
truncation argument.) Then
by letting $A \to\infty$ on the right-hand side of \eqref{sk1}, we
can obtain Theorem \ref{tleaves}(i) for the critical case.

\item[(3)] The subcritical case in Theorem \ref{tleaves} will be
proved in a similar way: By Nerman \cite{nerman81}, if we denote by
$B^{(j)}:= \#\cL^{(j)}[0] e^{ - \varrho_- V(u^{(j)})}$, then under
$\mathbf{P}$, conditioning on $\{V(u^{(j)}), 1\le j\le H(L)\}$ and letting
$L=L(n) \to\infty$, $B^{(j)}$ converges in law to $ c_{\mathrm{sub}}^*
M^{(\varrho_-, j)}_\infty$ where $c_{\mathrm{sub}}^*$ is some positive constant
and $M^{(\varrho_-, j)}_\infty, j\ge1$, are independent copies of
$M^{(\varrho_-)}_\infty$, and $M^{(\varrho_- )}_\infty$ is the
limit of
some positive martingale and has a power-law tail; see \eqref
{liuvarrho-}. As in the critical case, we get that under $\mathbf{P}$ and
conditioning on $\{ H(L) >0\}$,
\[
\#\cL[0] \approx\sum_{j=1}^{H_g(L)} e^{ \varrho_- V(u^{(j)})} B^{(j)}
\approx c^*_{\mathrm{sub}} e^{\varrho_- L} \sum_{i=1}^{\hat\zeta} e^{
\varrho_-
x_i} M^{(\varrho_-, i)}_\infty,
\]
with some point process $\sum_{i=1}^{\hat\zeta} \delta_{\{
x_i\}} $ on $(0, \infty)$ (this point process has of course nothing to
do with that in the critical case). Under $\Q$, $(M^{(\varrho_-,
i)}_\infty, i\ge1)$ are i.i.d., independent of $\sum_{i=1}^{\hat
\zeta
} \delta_{\{x_i\}} $ and distributed as $M^{(\varrho_-)}_\infty$ under
$\mathbf{P}$. By Theorem \ref{thmyaglom-critiq}(ii), $\mathbf{P}(
H(L) >0)
\sim\hat c_{\mathrm{sub}} e^{- \varrho_+ L}$ with some positive constant $\hat
c_{\mathrm{sub}}$. Let $L:= { \log n \over\varrho_-} - A$ with a large $A>0$.
It follows that as $n \to\infty$,
\begin{eqnarray*} n^{ \varrho_+/\varrho_-} \mathbf{P}\bigl( \#\cL[0]
> n
\bigr) & \sim& n^{ \varrho_+/\varrho_-} \Q\Biggl( c^*_{\mathrm{sub}} e^{\varrho_-
L} \sum_{i=1}^{\hat\zeta} e^{ \varrho_- x_i} M^{(\varrho_-,
i)}_\infty
> n \Biggr) \mathbf{P}\bigl( H(L) >0\bigr)
\\
&\sim& \hat c_{\mathrm{sub}} e^{ \varrho_+ A} \Q\Biggl( \sum
_{i=1}^{\hat\zeta} e^{ \varrho_- x_i} M^{(\varrho_-, i)}_\infty> {
1\over c^*_{\mathrm{sub}} } e^{ \varrho_- A}\Biggr),
\end{eqnarray*}
yielding the part (ii) in Theorem \ref{tleaves} by letting $A \to
\infty$.
\end{longlist}

\subsection{\texorpdfstring{Sketch of the proofs of Theorem \protect\ref{thmyaglom-critiq} and
Proposition \protect\ref{pyaglom}}
{Sketch of the proofs of Theorem 3 and Proposition 1}}

By the spine decomposition (see Proposition \ref{pchange-proba}), the
process $(S_k, k\ge0)$ formed by the positions of the spine $(\omega_k,
k\ge0)$ is a random walk under the probability $\Q$. Moreover $S$ has
zero mean in the critical case and positive mean in the subcritical
case. Let $\tau_t^+:= \inf\{ k\ge0\dvtx S_k > t\}$ and denote by $T^+_t:=
S_{ \tau^+_t} - t$ the overshoot. Then by the spine decomposition,
%
\begin{equation} \label{decmut} \mu_t:= \sum_{ u \in\H(t)} \delta
_{\{
V(u)- t\}} = \delta_{\{ T^+_t\}} + \sum_{k=1}^{ \tau_t^+}
\sum_{ u \in\I_k} \mu^{(u)}_t,
\end{equation}
where $\I_k$ denotes the set of brothers of $\omega_k$ at
$k$th generation [see \eqref{Ck}], and the point process $\mu_t^{(u)}$
is associated to the subtree $\cT^{(u)}$ (rooted at $u$) of $\cT$:
$\mu_{t}^{(u)}:=
\sum_{v \in\cT^{(u)} \cap\H(t)} \delta_{\{V(v)- t\}}$.

Consider a new probability $\Q^+$ defined in \eqref{Q+}.
Under $\Q^+$, $S$ is a random walk conditioned on staying nonnegative.
By \eqref{c31bis}, for any $f$ a nonnegative measurable function,
\begin{eqnarray}
\e\bigl[ e^{- \langle f, \mu_t\rangle} 1_{\{H(t) >0\}} \bigr] =
\Q
^+\biggl[ { e^{- \langle f, \mu_t\rangle} \over M^*_{ \cc_t}}\biggr]
\sim{1\over C_R {\mathscr R}(t)} e^{-\varrho t} \Q^+ \biggl[ { e^{-
\langle f, \mu_t\rangle} \over\int_\r e^{ \varrho z} \mu_t(dz) }
\biggr],\nonumber \\
\eqntext{t \to\infty,}
\end{eqnarray}
where $C_R $ denotes some positive constant and ${\mathscr
R}(t)$ and $\varrho$ are given in \eqref{def-Rt}. Therefore to prove
Theorem \ref{thmyaglom-critiq}, it is enough to check the convergence
in law of the point process $ \mu_t $ under $\Q^+$.

To this end, we first check that in the sum $\sum_{k=1}^{ \tau^+_t} $
in \eqref{decmut}, only those terms with $k$ near to $\tau^+_t$
contribute (see Lemma \ref{Lhut0} for the precise statement), and that
we may replace $\mu_t^{(u)}$ by $\widetilde\mu_t^{(u)}$ a point
process defined by some branching random walk starting from $V(u)$
without killing at $0$; see Lemma \ref{lemnonkilled}. Then by using
the convergence in law (Lemma \ref{lemlaw-convergence-reverse}) for
the time reversal random walk combined with the overshoot $\{ T_t^+,
S_{\tau_t^+}- S_{\tau_t^+- k}, 1\le k \le\tau_t^+\}$, we can obtain
the convergence of $\mu_t$ under $\Q^+$ and prove Theorem \ref
{thmyaglom-critiq}. Proposition \ref{pyaglom} will be proved in a
similar way.

\section{\texorpdfstring{From the number of leaves to the total progeny of the killed
branching walk: Proof of Lemma \protect\ref{lcom}}
{From the number of leaves to the total progeny of the killed
branching walk: Proof of Lemma 2}}
\label{scomparison}

We recall that our branching random walk starts from $x\ge
0$. We introduced for $u\in\cT$, $ \tau_a^-(u):= \inf\{0\le
k\le|u| \dvtx V(u_k)< a\}$ and
\[
\cL[a]:= \bigl\{ u\in\cT\dvtx|u|=\tau_a^-(u) \bigr\},\qquad a \le x.
\]

\begin{pf*}{Proof of Lemma \ref{lcom}} We equip the tree $\cT$
with the lexicographical order.
Let $U_k$ be the
$k$th vertex for this order in the set $\zz$ of the living particles.
It is well defined until $k=Z$ when all living particles have been
explored. For $k\in[1,Z]$, we introduce
\[
Y_k:= 1+ \sum_{i=1}^k \bigl(\nu(U_i)-1\bigr),
\]
and we notice that $Y_Z=\#\cL[0]$. (This can be easily checked
by using an argument of recurrence on the maximal generation of the
individuals of $\zz$.) We extend the definition
of $Y_k$ to $k>Z$, by $Y_{k+1}:=Y_k + \nu_k-1$ where $\nu_k$ is taken
from a family $\{\nu_{i},
i\ge1\}$ of i.i.d. random variables distributed as
$\nu(\varnothing)$ and independent of our branching random walk. We
claim that $(Y_k, k\ge1)$ is a random walk. To see this, observe that
we can construct the killed branching random walk in the following way.
Let $(\L^{(c)}_i, i\ge1)$ be i.i.d. copies of $\L$. At step 1, the
root $\varnothing=:U_1$ located at $x$ generates the point process $\L
^{(c)}_1$. If all the children are killed, we stop the construction.
Otherwise, we call $U_2$ the first vertex for the lexicographical order
that is alive. Then, $U_2$ generates the point process $\L_2^{(c)}$,
and we continue similarly. The process that we get has the law of the
killed branching random walk. In particular, if $ \nu_i^{(c)}$ denotes
the number of points of $\L_i^{(c)}$, then $(Y_k, k\ge1)$ has the law
of $(\sum_{i=1}^k (\nu_i^{(c)}-1), k\ge1)$ which is a random walk by
construction. This proves the claim. We suppose that Theorem \ref
{tleaves} holds and we want to deduce Theorem \ref{tmain}. Let us
look at the upper bound of $\mathbf{P}_x(Z>n)$. Let $m:= \e[\nu]>1$
and take
$\varepsilon\in(0,m-1)$. We have
\begin{eqnarray*}
\mathbf{P}_x\bigl(\#\cL[0]\le(m-1-\varepsilon)n, Z>n\bigr)
&=&
\mathbf{P}_x\bigl(Y_Z \le(m-1-\varepsilon)n, Z>n\bigr)\\
&=& \sum_{k> n} \mathbf{P}_x\bigl(Y_k \le(m-1-\varepsilon)n, Z=k \bigr)\\
&\le&
\sum_{k> n} \mathbf{P}_x\bigl(Y_k \le(m-1-\varepsilon)k\bigr),
\end{eqnarray*}
which is exponentially small by Cram\'er's bound. By Theorem
\ref{tleaves}, \mbox{$\mathbf{P}_x(\#\cL[0]>n)$} decreases polynomially. Therefore,
\begin{eqnarray*}
\mathbf{P}_x(Z>n) &\le& \mathbf{P}_x\bigl(\#\cL[0]>(m-1-\varepsilon)n\bigr)
+ \mathbf{P}_x\bigl(\#\cL[0] \le
(m-1-\varepsilon)n, Z>n\bigr)\\
&=&
\mathbf{P}_x\bigl(\#\cL[0]>(m-1-\varepsilon)n\bigr)\bigl(1+o(1)\bigr).
\end{eqnarray*}

Letting $n$ go to $\infty$, then $\varepsilon\to0$ yields
the upper bound. For the lower bound, we take $\varepsilon>0$, and we
observe that
\begin{eqnarray*}
\mathbf{P}_x\bigl(\#\cL[0] >(m-1+\varepsilon)n, Z\le n\bigr)
&=&
\mathbf{P}_x\bigl(Y_Z>(m-1+\varepsilon)n, Z\le n\bigr)\\
&\le&
\mathbf{P}_x\Bigl( \max_{1\le k \le n}\bigl (Y_{k}-(m-1)k\bigr) >\varepsilon n \Bigr).
\end{eqnarray*}

Let $\alpha> 2 $ in the critical case and $\alpha> 2
{\varrho_+ /\varrho_-}$ in the subcritical case. By Doob's $L^p$-inequality,
\[
\e\Bigl| \max_{1\le k \le n} \bigl(Y_{k}-(m-1)k\bigr)\Bigr|^\alpha\le{\alpha
^\alpha\over
(\alpha-1)^\alpha} \e\bigl| \bigl(Y_n-(m-1)n\bigr)\bigr|^\alpha,
\]
which according to Theorem 2.10 in Petrov \cite{petrov}, page 62, is
less than
\[
c(\alpha) n^{{\alpha/2} -1}
\e\sum_{i=1}^n |\nu_i- m|^\alpha= c(\alpha) n^{\alpha/2} \e
|\nu
-m|^\alpha,
\]
with some constant $c(\alpha)>0$. It follows that
\[
\mathbf{P}_x\bigl(\#\cL[0] >(m-1+\varepsilon)n,
Z\le n\bigr) \le{c(\alpha) \e|\nu-m|^\alpha\over\varepsilon^\alpha}
n^{-\alpha/2}.
\]
Therefore,
\begin{eqnarray*}
\mathbf{P}_x(Z>n) \ge\mathbf{P}_x\bigl(\#\cL[0]>(m-1+\varepsilon)n\bigr) -
{c(\alpha
) \e|\nu
-m|^\alpha\over\varepsilon^\alpha}
n^{-\alpha/2},
\end{eqnarray*}
which proves the lower bound by taking $n\to\infty$ then
$\varepsilon\to0$.
\end{pf*}

\section{One-dimensional real-valued random walks}\label{1dRW}

In this section we collect some preliminary results for a
one-dimensional random walk $(S_n)_{n\ge0}$ on some probability space
$(\Omega, {\mathscr F}, \ppp)$. Most of the results in this section
will be applied to the random walk $S$ defined in (\ref{S}) under $\Q$
in Section~\ref{secspinal-decompo}. For the sake of clarity of
presentation, the technical proofs are postponed to Section~\ref{proofsRW}.

\subsection{Time-reversal random walks}
\label{subsecretournement}

Let $((S_n)_n, \ppp_x)$ be a real-valued random walk starting from $x
\in\r$. We write $\ppp= \ppp_0$. Assume that $\eee[S_1] \ge0$ and
$\eee[|S_1|^{3+\delta}]<\infty$ for some $\delta>0$. In words, we
consider random walks that do not drift to $-\infty$. Moreover we
assume that the distribution of $S_1$ is nonarithmetic. Define
%
\begin{equation} \label{deftau}
\tau_a^+:= \inf\{k\ge0 \dvtx S_k> a\},\qquad
\tau_a^-:= \inf\{k\ge0 \dvtx S_k< a\},
\end{equation}
and the overshoot/undershoot
\[
T_a^+:= S_{\tau_a^+} -a > 0,\qquad T_a^-:= a - S_{\tau_a^-} > 0.
\]

Let $R(\cdot)$ be the renewal function of $(S_n)_{n\ge0}$ under
$\ppp$, that is, with $\tau^*:=\inf\{j\ge1\dvtx S_j \ge0\}$,
\[
R(x):= \eee\Biggl[\sum_{j=0}^{\tau^*-1} 1_{\{-x \le S_j\}} \Biggr]\qquad
\forall x >0,
\]
and $R(0)=1$.

Following \cite{bertoin-doney}, we introduce the law of the random
walk conditioned to stay nonnegative. To this aim, we see $(S_n)_{n\ge
0}$ under $\ppp_x$ as a Markov chain with transition function
$\mu(y,{d}z):= \ppp(y+ S_1 \in{d}z)$. We denote by
$\ppp^{+}_x$ the
$h$-transform of $\ppp_x$ by the function $R$. That is, $\ppp^{+}_x$
is a probability measure under which $(S_n)_{n\ge0}$ is a
homogeneous Markov chain on the nonnegative real numbers, with
transition function
%
\begin{equation}
\mu_R(y,{d}z):= \frac{R(z)}{R(y)} \mu(y,{d}z),\qquad
y,z \ge0.
\end{equation}
It is well known that $\ppp^{+}$-almost surely $S_n \to\infty$ when
$n \to\infty$. When $(S_n)_{n\ge0}$ drifts to $\infty$ (i.e., when
$\eee[S_1]>0$), $\ppp^{+}$ is the law of the random walk conditioned
to stay nonnegative in the usual sense, that is, $\ppp^{+}(\cdot
)=\ppp
(\cdot
\vert S_1\ge0,\ldots,S_n \ge0,\ldots)$.

We denote by $(\sigma_n, H_n)_{n\ge0}$ the
strict ascending ladder epochs and ladder heights of $S$ defined by
$(\sigma_0, H_0)=(0,0)$ and otherwise for $n \ge1$ by
\[
\sigma_n:= \cases{
\min\{ k > \sigma_{n-1} \dvtx S_{k} > H_{n-1} \}, & \quad$\mbox{if
$\sigma
_{n-1}<\infty$,}$ \vspace*{2pt}\cr
\infty, & \quad $\mbox{if $\sigma_{n-1}=\infty$}$}
\]
and
\[
H_n:= \cases{
S_{\sigma_n}, & \quad$\mbox{if $\sigma_n<\infty$,}$ \vspace*{2pt}\cr
\infty, & \quad$ \mbox{if $\sigma_n=\infty$.}$}
\]

Some results from random walk theory are important in the
proofs presented here and recorded in the following lemma.
%
\begin{lemma}
\label{lemRW} Assume that $\eee[S_1] \ge0$,
$\eee[|S_1|^{3+\delta}]<\infty$ for some $\delta>0$ and that the
distribution of $S_1$ is nonarithmetic. Then:

\begin{enumerate}[(iii)]
\item[(i)] $T_t^+$ converges in law to a finite random variable when
$t$ tends to infinity.
\item[(ii)] $(T^+_t, t\ge0)$ is bounded
in $L^p$ for all $1<p< 1+ \delta$.
\item[(iii)] $S_{ \tau_t^+}/t$ converges in probability to $1$ when $t$
tends to infinity.
\item[(iv)]
\begin{enumerate}
\item[$\bullet$] If $\eee[S_1]=0$, there exists a constant $C_R \in(0,\infty)$
such that $R(x)/x \longrightarrow C_R$ when $x \to\infty$. In this
case, $C_R= {1 \over\eee[T_0^{-}]}= {1 \over\eee[- S_{\tau_0^{-}}] }$.
\item[$\bullet$] If $\eee[S_1]>0$, there exists a constant $C_R \in(0,\infty)$
such that $R(x) \longrightarrow C_R$ when $x \to\infty$. In this case,
$C_R={1\over\ppp( \tau^-_0= \infty)}$.
\end{enumerate}
\item[(v)]
\begin{enumerate}
\item[$\bullet$] If $\eee[S_1]=0$, then $\ppp(\tau_{t}^+<\tau_0^-) \sim
{1\over
C_R t} $ when $t \to\infty$.
\item[$\bullet$] If $\eee[S_1]>0$, then $ \ppp(\tau_{t}^+<\tau_0^-) \to
{1\over
C_R}$ when $t \to\infty$.
\end{enumerate}
\end{enumerate}
\end{lemma}

\begin{pf}
Notice that $T_t^+$ is also the overshoot of the
random walk $(H_n)$ above the level $t$. In the case $\eee[S_1]=0$,
Doney \cite{D80} implies that $H_1$ has a finite $(2+\delta)$-moment
which in view of Lorden \cite{L70}, Theorem 3, applied to $(H_n)$,
implies that $(T^+_t, t\ge0)$ is bounded in $L^p$ for all $1<p< 1+
\delta$. In the case $\eee[S_1]>0$, again by Lorden \cite{L70}, Theorem
3, applied to $(S_n)$, $(T^+_t, t\ge0)$ is bounded in $L^p$ for all
$1<p< 2+ \delta$. This provides Part (ii)
of the lemma. Moreover, we see that in both cases, $H_1= T^+_0$ has a
finite expectation and obviously is nonarithmetic, then a refinement
of the renewal theorem gives part (i) of the lemma (Feller \cite
{feller}, page~370, equation (4.10)). For both cases, part (iii) is a
consequence of part (ii). To show (iv), we recall the duality lemma
(Feller \cite{feller}, page 395),
\[
R(x)=1+ \sum_{n=1}^\infty\ppp\bigl( H_n^{-} \le x\bigr), \qquad x>0,
\]
where $(H_n^{-}, n\ge0)$ denotes the (strict) ascending ladder heights
of $-S$ (in particular, $H^-_1= T^-_0$ the undershoot at $0$). In the
case $\eee[S_1]=0$, part (iv) is a
consequence of the renewal theorem (see Feller \cite{feller}, page
360) with $C_R= {1 \over\eee[T_0^{-}]}$, while part (v) is obtained by
applying the optional stopping theorem to the martingale $(S_k, 0\le k
\le\tau^+_t\wedge\tau_0^-)$ (the uniform integrability is guaranteed
by (ii); see \cite{aidekon}, Lemma 2.2).
In the case $\eee[S_1]>0$, parts (iv) and (v) follow from the duality
lemma, $C_R=\eee[ \tau^*]= \lim_{x\to\infty} R(x)= 1+ \sum
_{n=1}^\infty\ppp( H_n^- < \infty) = 1+ \sum_{n=1}^\infty\ppp(
\tau
^-_0 <\infty)^n = {1\over\ppp( \tau^-_0= \infty)}$.
\end{pf}

\begin{figure}

\includegraphics{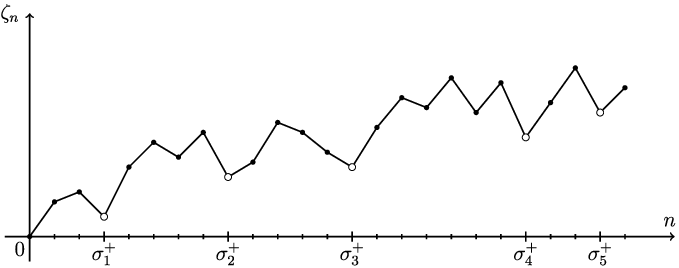}

\caption{Tanaka's construction.}\label{fig1}
\end{figure}

We recall now Tanaka's construction (see \cite{tanaka} and Figure~\ref{fig1}) of the random walk conditioned to
stay positive. Let us recall that $(\sigma_n, H_n)_{n\ge0}$ are the
strict ascending ladder epochs and ladder heights of $S$, and let
$(w_i)_{i \ge1}$ be independent copies of the segment of the random
walk $(S_n)_{n\ge0}$ up to time $\sigma:=\sigma_1$ viewed from
$(\sigma,S_\sigma)$ in reversed time and reflected in the $y$-axis;
that is, $(w_i)_{i \ge0}$ are independent copies of
%
\begin{equation}
(0,S_\sigma-S_{\sigma-1},S_{\sigma}-S_{\sigma-2},\ldots,S_{\sigma
}-S_{1},S_{\sigma}).
\end{equation}
We write now $w_i=(w_i(\ell); \ell= 0, 1, 2, \ldots
\sigma^{(i)})$ to identify the components of $w_i$. In
\cite{tanaka}, Tanaka shows that the random walk conditioned to stay
positive can be constructed by gluing the $w_i$'s together, each
starting from the end of the previous one. More formally, let
$(\sigma^+_n, H^+_n)_{n\ge1}$ be the renewal process formed from the
independent random variables $(\sigma^{(i)},w_i(\sigma^{(i)}))$,
that is,
%
\begin{equation}
\qquad\bigl(\sigma^+_n,H^+_n\bigr)=\bigl(\sigma^{(1)}+\cdots+\sigma^{(n)},w_1\bigl(\sigma
^{(1)}\bigr)+\cdots+w_n\bigl(\sigma^{(n)}\bigr)\bigr),\qquad
n \ge1.
\end{equation}
Then, Tanaka's result says that the random walk conditioned to stay
positive can be constructed via the process $(\zeta_n)_{n \ge0}$
given by
%
\begin{equation}
\label{eqtanakaconstruction}
\zeta_n=H_k^+ + w_{k+1}\bigl(n-\sigma^+_k\bigr),\qquad \sigma^+_k < n \le
\sigma^+_{k+1}.
\end{equation}


%
%
%
%
%
%
%
%
%
%
%
%
%
%
%
%
%
%
%
%


Finally we introduce a process $(\hat S_n)_{n \ge0}$ (obtained by
modifying slightly the random walk conditioned to stay positive)
which will be the limit process that appears in the following lemma.
Let $\widetilde\sigma:= \sup\{n \ge1 \dvtx\zeta_n=\min_{1 \le i
\le n}
\zeta_i\}$ and observe that $\widetilde\sigma$ is almost surely finite
since $\zeta_n \to\infty$. Then $(\hat S_n)_{n \ge0}$
is defined by
%
\begin{equation}
\label{eqdef-xi}
\eee\bigl[F\bigl((\hat S_n)_{n \ge0}\bigr) \bigr]= \frac{1}{\eee[ H_1
]} \eee\bigl[ \zeta_{\widetilde\sigma} F\bigl((\zeta_n)_{n
\ge
0}\bigr) \bigr]
\end{equation}
for any test function $F$. Observe that Tanaka's construction implies
$\eee[ \zeta_{\widetilde\sigma}] =\eee[ H_1 ]$. Moreover
we introduce $\hat\sigma:= \sup\{n \ge1 \dvtx\hat S_n=\min_{1
\le i \le n} \hat S_i\}$ which is almost surely finite since $\hat S_n
\to\infty$.

\begin{lemma}
\label{lemlaw-convergence-reverse} Assume that $\eee[S_1] \ge0$,
$\eee[|S_1|^{3+\delta}]<\infty$ for some $\delta>0$ and that the
distribution of $S_1$ is nonarithmetic. Recall (\ref{deftau}), and fix
an arbitrary integer $K \ge1$. Let $F\dvtx\r^*_+ \times\r_+^{K }
\to
\r$
be a bounded and measurable function. Suppose that for any $z \in\r
_+^K$, the set $\{x\in\r_+^*\dvtx F(\cdot, z) \mbox{ is not
continuous at
} x \}$ is at most countable [which may depend on $z$]. Then
\begin{longlist}[(ii)]
\item[(i)]
\[
\lim_{t \to\infty} \eee\bigl[ F\bigl( T_t^+,(S_{\tau_t^+}-S_{\tau
_t^+-j})_{1\le j \le K}\bigr) | \tau_t^+ > K \bigr] = \eee\bigl[ F\bigl(U
\hat S_{\hat\sigma}, (\hat S_j)_{1 \le j \le K}\bigr) \bigr],
\]
where $(\hat S_n)_{n \ge0}$ is the process defined by (\ref
{eqdef-xi}) and $U$ is a uniform random variable on $[0,1]$
independent of $(\hat S_n)_{n \ge0}$.\vspace*{1pt}

\item[(ii)]
\[
\lim_{t \to\infty} \eee^{+} \bigl[ F\bigl( T_t^+,(S_{\tau_t^+}-S_{\tau
_t^+-j})_{1 \le j \le K}\bigr) | \tau_t^+ > K \bigr] = \eee
\bigl[ F\bigl(
U \hat S_{\hat\sigma}, (\hat S_j)_{1 \le j \le K}\bigr) \bigr],
\]
where $\eee^{+}$ denotes the expectation with respect to the
probability measure $\ppp^{+}$.
\end{longlist}
\end{lemma}

As a consequence, under $\ppp(\cdot \vert \tau_t^+ > K)$ or under
$\ppp^{+}(\cdot \vert \tau_t^+ > K)$, the random vector
$(T_t^+,(S_{\tau_t^+}-S_{\tau_t^+-j})_{1 \le j \le K})$ converges in
distribution toward $(U \hat S_{\hat\sigma},\break (\hat S_j)_{1 \le j \le
K})$ when $t \to\infty$.
We also note that the conditioning with respect to the event $\{\tau
_t^+ > K\}$ is just technical since this event is asymptotically
typical (indeed almost surely $\tau_t^+ \to\infty$ when $t \to
\infty$).

\begin{pf*}{Proof of Lemma \ref{lemlaw-convergence-reverse}}
See Section~\ref{proofsRW}.
\end{pf*}

We end this subsection with an estimate on a random walk with positive drift:

\begin{lemma}
\label{lpositivewalk} Assume that $\eee[S_1] > 0$, $\eee[
S_1^2]<\infty
$. Let $(a_i, S_i- S_{i-1})_{i\ge1}$ be an i.i.d. sequence such that
$a_i\ge0$ almost surely. For any $p \ge1$ such that $ \eee[ a_1^p] <
\infty$ and for any $\kappa>0$, there exists some constant $c_{p,
\kappa
}>0$ such that
%
\begin{equation}
\eee_x\Biggl[ \sum_{k=0}^{\tau_t^+ -1} a_{k+1} e^{ \kappa (S_{k}-t)
} \Biggr]^p \le c_{p, \kappa}\qquad \forall t>0, \forall x\le t.
\label{zb2new}
\end{equation}
\end{lemma}

\begin{pf}
See Section~\ref{proofsRW}.
\end{pf}

\subsection{Centered random walks}

Let $((S_n)_{n\ge0}, \ppp_x)$ be a real-valued random walk starting
from $x \in\r$. We write $\ppp=\ppp_0$. Assume that
%
\begin{equation} \label{Hyp-expmoment} \qquad\eee[S_1]=0,\qquad \operatorname{Var}(S_1)>0,\qquad
\eee\bigl[e^{u S_1} \bigr]< \infty \qquad\forall u \in
\bigl(-(1+\eta),\eta\bigr)
\end{equation}
for some $\eta>0$. Recall that $ \ppp(\tau_L^+<\tau_0^-) $ is of order
$ {1\over L}$ as $ L\to\infty$ (cf. Lemma~\ref{lemRW}). For $a \in
\r$, recall that $T_a^+:=S_{\tau_a^+}-a>0$ (resp., $T^-_a:= a-
S_{\tau
_a^-}>0$)
denotes the overshoot (resp., undershoot) at level $a$.

We have the following estimate.

\begin{lemma} \label{lemestcrit} Under (\ref{Hyp-expmoment}). For any
$0<\delta<\eta$, there exist some constants $c>1$ and $c'=c'(\delta
)>1$ such that for all $b\ge a \in\r$ and $ x>0$,
%
\begin{eqnarray} \ppp_a\bigl( T_b^+ >x\bigr) &\le& c' e^{ - \delta x}, \label
{ta+} \\
\ppp_b\bigl( T_a^- >x\bigr) &\le& c ' e^{ -
(1+\delta) x} \label{ta-}.
\end{eqnarray}

Moreover, for all $L\ge1$, $0 \le a \le L$,
%
\begin{eqnarray}
\ppp_a \bigl( \tau_0^- <
\tau_L^+\bigr) & \le& { L -a + c \over L}, \label{pal1} \\
\ppp_a \bigl( \tau_0^- >
\tau_L^+\bigr) & \le& { a + c \over L}, \label{pal2}
\\ \eee_a \bigl[ e^{- S_{\tau_0^-}} 1_{\{\tau_0^- < \tau_L^+\}
}\bigr]
&\le& c { L -a +1 \over L}, \label{pal}
\\
\eee_a \Biggl[
\sum_{j =0}^{\tau_L^+-1} e^{- \delta(L-S_j)}\Biggr] + \eee_a \Biggl[
\sum_{j =0}^{\tau_0^--1} e^{- \delta S_j }\Biggr] &\le& c',
\label{foot3}
\\
\eee_a \bigl[ e^{S_{\tau_{0}^- -1} - S_{\tau_0^-}}\bigr]
&\le&
c,
\label{diftau0}\\
\eee_a \biggl[ \sum_{0\le j < \tau_0^- \wedge\tau_L^+} e^{- \delta S_j}
\biggr] & \le& c' { L-a +1\over L},
\label{L-a+1}
\\
\eee_a
\biggl[ \sum_{0\le j < \tau_0^- \wedge\tau_L^+} e^{- \delta(L-S_j)}
\biggr] & \le& c' { a +1\over L},
\label{a+1}
\\
\eee_a \biggl[
e^{- S_{\tau_0^-}} 1_{\{ \tau_0^- < \tau_L^+\}} \sum_{ 0\le j <
\tau_0^-} e^{ - \delta(L-S_j)}\biggr] &\le& c' { a +1
\over L^2}.
\label{eal}
\end{eqnarray}
\end{lemma}

\begin{remark*}
A weaker assumption $\sup_{ - \eta\le u \le
\eta} \eee[e^{u S_1} ]< \infty$ is enough to get~(\ref
{foot3}), (\ref{diftau0}), (\ref{L-a+1}) and (\ref{a+1}).
\end{remark*}

\begin{pf*}{Proof of Lemma \ref{lemestcrit}}
See Section~\ref{proofsRW}.
\end{pf*}

\subsection{Random walks with negative drift}

In this subsection, we give estimates on transient random walks. We
take again $((S_n)_{n\ge0}, \ppp_x)$ a random walk, but we suppose
now that $\eee[S_1]<0$, hence the random walk drifts to $-\infty$. We
suppose that there exist $\gamma, \eta_1, \eta_2>0$ such that
%
\begin{equation}\label{nega} \eee\bigl[e^{\gamma S_1}\bigr]=1,\qquad \eee\bigl[e^{u
S_1}\bigr]<\infty\qquad \forall u \in(-\eta_1,\gamma+\eta_2).
\end{equation}
Then
%
\begin{equation} \label{subtau-}
\ppp\bigl(\tau_a^-<\tau_0^+\bigr) \to\ppp\bigl(\tau_0^+=\infty\bigr)>0,\qquad a\to
-\infty.
\end{equation}

By Theorem 1 of \cite{iglehart}, if $S_1$ is nonarithmetic, then
%
\begin{equation}\label{subtau+}
\ppp\bigl(\tau_a^+<\tau_0^-\bigr) \sim c(\gamma) e^{-\gamma a},\qquad a
\to
+\infty
\end{equation}
for some constant $c(\gamma) >0$. We end this section with two lemmas:

\begin{lemma}\label{lemestsub} Under (\ref{nega}). For any $r >0$,
we can find some positive constants $c, c', c''$ such that for any
$a\ge0$, $L>1$,
%
\begin{eqnarray}\label{est1sub}
\eee_a\bigl[ e^{-r S_{\tau_0^-}} \bigr] &\le& c(r),\qquad \mbox{if }
r < \eta_1,
\\
\label{est3sub}\eee_a\biggl[ \sum_{0\le\ell< \tau_L^+}(1+L-S_{\ell})^{\alpha
}e^{r
S_{\ell}} \biggr]
&\le&
c'(r,\alpha)e^{\gamma(a-L)}e^{r L},
\nonumber
\\[-8pt]
\\[-8pt]
\eqntext{\mbox{if } r > \gamma
,
\alpha\ge0.}
\\
\label{est7sub}\eee_a\Biggl[ \sum_{\ell=0}^{\min(\tau_0^-,\tau_L^+)} (1+L- S_\ell
)^\alpha e^{\gamma S_{\ell}} \Biggr]
&\le&
c'' e^{\gamma a} (1+L-a)^{1+\alpha},
\nonumber
\\[-8pt]
\\[-8pt]
\eqntext{ a\in[0,L], \alpha\ge0.}
\end{eqnarray}
\end{lemma}

\begin{pf}
See Section~\ref{proofsRW}.
\end{pf}

\begin{lemma}\label{Laux1} Under (\ref{nega}). Fix some $ 0\le\eta<
\eta_1$ and $\alpha\ge0$. Assume that $(S_n-S_{n-1}, a_n)_{n\ge1}$ are
i.i.d. with $a_1\ge0$ almost surely.
\begin{longlist}[(ii)]
\item[(i)] Assume $b>0$, $0\le p < \gamma/b$ and $a_1$ are such that $ \eee
[ (1+ 1_{\{S_1<0\}} e^{ - \eta S_1}) \times  a_1^p ] < \infty$.
There exists some constant $c_p=c_p(b, \eta, \alpha)>0$ such that for
all $x \ge0$,
%
\begin{equation}\label{est8sub} \eee_x \Biggl[ e^{ - \eta S_{ \tau
^-_0}} \Biggl( \sum_{\ell=1}^{ \tau^-_0} e^{ b S_{\ell-1}} a_\ell
\Biggr)^p \Biggr] \le c_p e^{ b p x}.
\end{equation}

\item[(ii)] Assume $b>0$, $p \ge1 $ and $a_1$ are such that $ \eee[ (1+
1_{\{S_1<0\}} e^{ - \eta S_1}) a_1^p ] < \infty$ and $\eee[ e^{
p b S_1} ] < \infty$. There exists some constant $c_p=c_p(b, \eta,
\alpha)>0$ such that for all $L>0$ and $0\le x \le L$,
%
\begin{eqnarray}\label{est9sub} && \eee_x \Biggl[ e^{ - \eta S_{ \tau
^-_0 }} \Biggl( \sum_{\ell=1}^{ \min(\tau^-_0, \tau_L^+)}
(1+L-S_{\ell
-1})^\alpha e^{ b S_{\ell-1}} a_\ell\Biggr)^p \Biggr]
\nonumber
\\[-8pt]
\\[-8pt]
\nonumber
&&\qquad\le
c_p \times \cases{
(1+L-x )^{\alpha p } e^{ p b x}, & \quad$\mbox{if $ p < \gamma/b$,}$
\vspace*{2pt}\cr
e^{\gamma x} (1+ L-x)^{ 1+ \alpha p }, &\quad $\mbox{if $ p = \gamma/b$,}$
\vspace*{2pt}\cr
e^{\gamma(x-L) + p b L}, &\quad$ \mbox{if $ p > \gamma/b$.}$}
\end{eqnarray}
\end{longlist}
\end{lemma}

\begin{pf}
See Section~\ref{proofsRW}.
\end{pf}

\section{Spinal decomposition}
\label{secspinal-decompo}

\subsection{Spinal decomposition of a branching random walk (without killing)}
We begin with a general formalism of the spinal decomposition for a
branching random walk. This decomposition has already been used in the
literature
by many authors in various forms; see, for example, Lyons, Pemantle and
Peres~\cite{lyons-pemantle-peres},
Lyons~\cite{lyons97} and Biggins and Kyprianou \cite{BigKyp04}.

There is a one-to-one correspondence between the branching random
walk $(V(u)_{u \in\cT})$ and a marked tree $\{(u, V(u))\dvtx u \in
\cT\}$. For $n \ge1$, let $\F_n$ be the sigma-algebra generated by
the branching random walk in the first $n$ generations. For any $u
\in\cT\setminus\{\varnothing\}$, denote by ${\buildrel
\leftarrow\over u}$
the parent of $u$. Write as before $[\varnothing, u]=\{
u_0:=\varnothing
, u_1, \ldots, u_{|u|}\}$
the shortest path from the root $\varnothing$ to $u$(with
$|u_i|=i$ for any $0\le i \le|u|$).

Let $h\dvtx\cT\to[0, \infty)$ be measurable such that
$h(\varnothing)
>0$ and for any $x \in\r$, $v \in\cT$ with $|v|=n \ge0$,
%
\begin{equation}\label{h}
\e_x \biggl[ \sum_{{\buildrel\leftarrow\over u}=v} h( u )
\Big| \F
_n\biggr]= \lambda
h(v),
\end{equation}
where $\lambda>0$ is some positive constant. Let
$\H_+:=\{u \in\cT\dvtx h(u) >0\}$. In our examples of $h$ in this paper,
$\lambda=1$, $h(u)= f(V(u))$ or $h(u)= f(V(u_1), \ldots,\break V(u_{|u|}))$ for
some nonrandom function $f$, and $\H_+$ equals either $\cT$ or $\zz$
the set of progeny of the
killed branching walk.

Define
\[
W_n:= {1\over h(\varnothing) \lambda^n} \sum_{|u|=n} h( u),\qquad
n\ge0.
\]

Fix $x\in\r$. Clearly by (\ref{h}), $(W_n)$ is a $(\mathbf{P}_x,
(\F_n))$-martingale.


On the enlarged probability space formed by
marked trees with distinguished rays, we may construct a
probability $\Q_x^{(h)}$ and an infinite ray $\{
w_0= \varnothing, w_1, w_2, \ldots\}$ such that for any $n\ge1$,
${\buildrel\leftarrow\over w_n}= w_{n-1}$, and
%
\begin{equation}
\Q_x^{(h)}( w_n = u |
\F_n ) =
{ h(u) \over h(\varnothing) \lambda^n W_n }\qquad \forall |u|=n
\label{u*}
\end{equation}
and
%
\begin{equation} {d \Q_x^{(h)}\over d \mathbf{P}_x}\bigg|_{\F_n} = W_n.
\label
{Q}
\end{equation}


To construct $ \Q_x^{(h)}$, we follow Lyons \cite{lyons97} under a
slightly more general framework: Let $\L:= \sum_{|u|=1}
\delta_{\{V(u)\}}$. For any $y \in\H_+$, denote by $\widetilde
\L_y$ a random variable whose law has the Radon--Nikodym density $W_1$
with respect to the law of $\L$ under $\mathbf{P}_y$. Put one particle
$w_0= \varnothing$ at $x \in\H_+$. Generate offsprings and
displacements according to an independent copy of $\widetilde
\L_x$. Let $\{|u|=1\}$ be the set of the children of~$w_0$. We
choose $w_1=u$ according to the probability $ {h(u)\over h(w_0)
\lambda W_1 }$.\vspace*{1pt} All children $u \neq w_1$ give rise to
independent branching random walks of law $\mathbf{P}_{V(u)}$, while
conditioned on $V(w_1)=y$, $w_1$ gives offsprings and displacements
according to an independent copy of $\widetilde\L_y$. We choose
$w_2$ among the children of $w_1$ in the same size-biased way, and so
on. Denote by $ \Q_x^{(h)}$ the joint law of the marked tree
$(V(u))_{|u|\ge0}$ and the infinite ray $\{ w_0= \varnothing, w_1,
\ldots, w_n, \ldots\}$. Then $\Q_x^{(h)}$ satisfies~(\ref{Q}) and
(\ref{u*}), which can be checked in the same way as in Lyons
\cite{lyons97}.




Under $Q^{(h)}_x$, we write, for $k \ge1$,
%
\begin{equation}
\I_k
:=
\bigl\{ u \dvtx|u| =k,
{\buildrel\leftarrow\over u} = w_{k-1},
u \neq w_k \bigr\}.
\label{Ck}
\end{equation}

In words, $\I_k $ is the set of children of $w_{k-1} $
except $w_k $, or equivalently, the set of the brothers of $w_k$,
and is possibly empty. Define $S_0:=V(\varnothing)$ and
%
\begin{equation} \label{Theta} S_n:= V(w_n),\qquad
\Theta_n:= \sum_{ u
\in\I_n} \delta_{\{\Delta V(u)\}},\qquad n\ge1,
\end{equation}
where we recall that $
\Delta V(u):= V(u) - V({\buildrel\leftarrow\over u})$. Finally, let
us introduce the following
sigma-field:
%
\begin{equation}
\G_n:=
\sigma\bigl\{
\bigl(\Delta V(u), u\in\I_k\bigr),
V(w_k), w_k, \I_k, 1\le k\le n\bigr\}.
\label{Gn}
\end{equation}

Then $\G_\infty$ is the sigma-field generated by all random variables
related to the spine $\{ w_k, k\ge0\}$. Let us write $v< u$ if $v$ is
an ancestor of $u$ (then $v\le u$ if $v<u$ or $v=u$). By the standard
``\textit{words}''-representation in a tree, $u<v$ if and only if the word
$v$ is a concatenation of the word $u$ with some word $s$, namely $v =
us$ with $|s|\ge1$.

The promised spinal decomposition is as follows. Since it
differs only slightly from the spinal decomposition presented in
Lyons~\cite{lyons97} and Biggins and Kyprianou~\cite{BigKyp04}, we
feel free to omit the proof.

\begin{proposition}
\label{pchange-proba}
Assume (\ref{h}) and fix $x \in\H_+$.
Under probability $\Q_x^{(h)}$:
\begin{longlist}[(ii)]
\item[(i)] for each $n \ge1$, conditionally on $\G_{n-1}$ and on $\{
S_{n-1}=y\}$, the point process $(V(u), {\buildrel\leftarrow\over
u}= w_{n-1})$ is distributed as $\widetilde\L_y$. In particular, the
process $(S_n, \Theta_n)_{ n\ge0}$ is Markovian. Moreover,
$(S_n)_{n\ge
0}$ is also a Markov chain and satisfies
\[
\Q_x^{(h)} \bigl( f(S_n) |
S_{n-1}=y, \G_{n-1}\bigr)= {1\over\lambda} \e_y \biggl[ \sum_{|u|=1}
f\bigl(V(u)\bigr) {h( u ) \over h(\varnothing)} \biggr]
\]
for any nonnegative measurable function $f$, $ n\ge1$ and $ y\in\H_+$.

\item[(ii)] Conditionally on $\G_\infty$, the shifted branching
random walks $\{ V(vu) - V(v)\}_{|u|\ge0} $, for all
$v\in\bigcup_{k=1}^\infty\I_k$, are independent, and have the
same law as $\{V(u)\}_{|u|\ge0}$ under $\mathbf{P}_0$.
\end{longlist}
\end{proposition}


Remark that under $\Q_x^{(h)}$, $\{w_n, n\ge0\}$ lives
in $\H_+$ with probability one.
We can extend Proposition \ref{pchange-proba} to the so-called
stopping lines. Recall (\ref{deftau+}) and (\ref{deftau-}).
For $0\le x < t$, we consider the stopping line
%
\begin{equation}
\label{defcct}
\cc_t:=\bigl\{ u \in
\cT\dvtx \tau_t^+(u) = |u| \bigr\}.
\end{equation}


%
%
%
%
%
%
%
%
%
%
%
%
%
%
%
%
%
%
%
%
%
%
%
%
%
%
%
%
%
%
%
%
%
%
%
%
%
%
%
%
%
%
%
%
%
%
%
%
%
%
%
%
%
%
%
%
%
%
%
%
%
%
%
%
%
%
%
%
%
%
%
%
%
%
%
%
%
%
%
%
%
%
%
%
%
%
%
%
%
%
%
%
%
%
%
%
%
%
%
%
%
%
%
%
%
%
%
%
%
%
%
%
%
%
%
%
%
%
%
%
%
%
%
%
%
%
%
%
%
%
%
%
%
%
%
%
%
%

\begin{figure}

\includegraphics{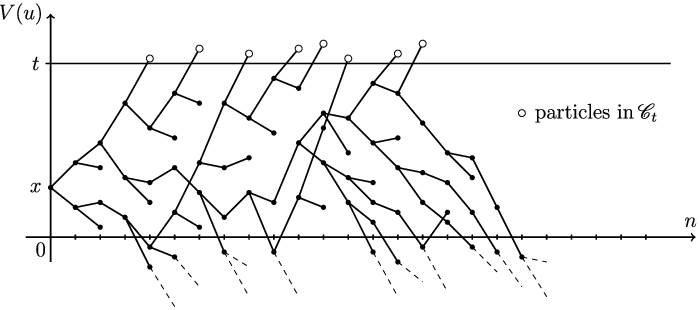}

\caption{The set $\cc_t$.}\label{fig2}
\end{figure}


Note that for any $v\in\cT$, $ |v|<
\tau_t^+(v)$ means that $ \sup_{0\le i \le|v|} V(v_i) \le t$; see
Figure~\ref{fig2}.
The process $\{V(u)\}_{|u|\le\tau_t^+(u)}$ can be interpreted as
the branching random walk stopped by the line $\cc_t$. Recalling (\ref
{defhl2}), we remark that $ \cc_t\cap\zz= \H(t)$, where as before
$\zz
$ denotes the set of progeny of the killed branching random walk.

Let $\F_{\cc_t}:= \sigma\{(u,V(u))\dvtx u \in\cT, |u| \le
\tau_t^+(u)\}$ be the $\sigma$-field generated by the branching walk
$V$ up to the stopping line $\cc_t$. Assuming (\ref{h}), we define
\[
W_{\cc_t}:=
{1 \over h(\varnothing)} \sum_{ u \in\cc_t } h(u) \lambda^{-
\vert u
\vert}.
\]

Define two families of stopping times for the process ($S_n:= V(w_n),
n\ge0$) under $\Q_x^{(h)}$,
%
\begin{equation}\label{tauS}\qquad  \tau^+_a:= \inf\{k\ge0\dvtx S_k >a\},\qquad
\tau_a^-:= \inf\{k\ge0\dvtx S_k < a\}\qquad \forall a\in\r,
\end{equation}
with the usual convention $\inf\varnothing=\infty$ and the corresponding
overshoot and undershoot processes
%
\begin{equation}\label{defovershoot} T_a^+:= S_{\tau_a^+} -a,\qquad
T_a^-:= a-S_{\tau_a^-}\qquad \forall a \in\r.
\end{equation}

Analogously to (\ref{Gn}), we introduce the sigma-field
%
\begin{equation}\label{Gct} \G_{\cc_t}:=
\sigma\bigl\{
\bigl(\Delta V(u), u\in\I_k\bigr),
V(w_k), w_k, \I_k, 1\le k \le\tau_t^+, \tau^+_t
\bigr\},
\end{equation}
generated by all information related to the spine $[\varnothing,
w(\tau
_t^+ ) ]$. Similarly, we recall $\cL[a]$ in (\ref{l[a]}) and define
$\F
_{\cL[a]}$, $W_{\cL[a]}$, $\G_{\cL[a]}$ as before. The next result
describes the decomposition along the spine $[\varnothing, w(\tau_t^+ )
]$ (resp., $[\varnothing, w(\tau_a^- ) ]$).

\begin{proposition}\label{pchange-proba2}
Assume (\ref{h}), and let $x\in\H_+$. Take $t\ge x$. Suppose that $h$
is such that $\Q_x^{(h)}(\tau_t^+<\infty)=1$. Then
%
\begin{equation} \label{c31} {d \Q_x^{(h)} \over d \mathbf{P}_x}
\bigg|_{\F_{\cc_t}}= W_{\cc_t}.
\end{equation}
\begin{longlist}[(iii)]
\item[(i)] Under probability $\Q_x^{(h)}$, conditionally on $\G_{\cc_t}$
and on $\{V(v)= x_v, v\in\bigcup_{k=1}^{\tau_t^+ } \I_k\}$, the shifted
branching
random walks $\{ V(vu) - V(v)\}_{u \dvtx|vu|\le\tau_t^+(vu)} $,
stopped by the line $\cc_t$, are independent, and have the
same law as\break $\{V(u)\}_{|u|\le\tau_{t-x_v}^+(u)}$ under $\mathbf{P}_0$,
stopped by the line $\cc_{t-x_v}$.

\item[(ii)] The distribution of the spine within $\cc_t$ is given by
\[
\Q^{(h)}_x( w_{\tau^+_t} = u | \F_{\cc_t} ) = {h(u)
\lambda
^{- \vert u\vert} \over h(\varnothing) W_{\cc_t} }\qquad \forall
u\in\cc_t.
\]

\item[(iii)] For any bounded measurable function $f\dvtx\r^\n\to\r$ and
for any bounded $\F_{\cc_t}$-measurable random variable $\Phi_t$,
\[
\e_x\biggl[ \sum_{ u \in\cc_t} {h (u) \over h(\varnothing) \lambda
^{|u|}} f\bigl( V(u_i), 0\le i \le\vert u \vert\bigr) \Phi_t \biggr] =
\Q
^{(h)}_x \bigl[ f\bigl( S_i, 0\le i \le\tau^+_t\bigr) \Phi_t \bigr].
\]

Similarly, take $a\le x$ and assume that $h$ is such that $\Q
_x^{(h)}(\tau_a^-<\infty)=1$. Then the analog holds for $\cc_t$
replaced by $\cL[a]$ (and $\tau_t^+$ by $\tau_a^-$).
\end{longlist}
\end{proposition}

\begin{remark} If $\Q_x^{(h)}(\tau_t^+<\infty) =1$ for all $t$, then
$W_{\cc_t}$ is a $(\mathbf{P}_x,
\F_{\cc_t})$-martingale by Lemma
6.1 and Theorem 6.1 in \cite{BigKyp04}. The equivalent holds for~$\cL
[a]$. 
\end{remark}

\begin{pf*}{Proof of Proposition \ref{pchange-proba2}} It is
enough to prove that for any $g\dvtx\cT\to\r$ measurable and bounded,
%
\begin{eqnarray}\label{eq3a}&&\e_x\biggl[ \sum_{ u \in\cc_t} {h (u)
\over
h(\varnothing) \lambda^{|u|}} f\bigl( V(u_i), 0\le i \le\vert u \vert
\bigr) g(u) \Phi_t \biggr]
\nonumber
\\[-8pt]
\\[-8pt]
\nonumber
&&\qquad= \Q^{(h)}_x \bigl[ f\bigl( S_i, 0\le i \le
\tau
^+_t\bigr) g(w_{\tau^+_t}) \Phi_t \bigr].
\end{eqnarray}

In fact, part (iii) follows from (\ref{eq3a}), and by taking
$f\equiv g \equiv1$ in (\ref{eq3a}) we get~(\ref{c31}); Taking
$f\equiv1$ in (\ref{eq3a}) and using (\ref{c31}), we get part (ii);
Finally since $\tau_t^+$ is a stopping time for $(S_k)_k$, the part (i)
follows easily from Proposition
\ref{pchange-proba}.

To check
(\ref{eq3a}), it is enough to show that for any $N \ge1$, (\ref{eq3a})
holds for all $\Phi_t$ of form $\Phi_{t, N}:= F(u,V(u), u \in\cT, |u|
\le\tau_t^+(u) \wedge N)$ with a bounded measurable function $F$.
Notice that the left-hand side of (\ref{eq3a}) equals
%
\begin{eqnarray} \label{eq3b} &&\sum_{n =0}^\infty\e_x\biggl[ \sum
_{|u|=n} 1_{\{
\tau_t^+(u)=n\}} {h (u) \over h(\varnothing) \lambda^{n}} f\bigl( V(u_i),
0\le i \le n \bigr) g(u) \Phi_{t, N} \biggr]
\nonumber
\\[-8pt]
\\[-8pt]
\nonumber
&&\qquad:= \sum_{ n=0}^\infty
(\ref{eq3b})_n,
\end{eqnarray}
with obvious definition of $(\ref{eq3b})_n$. If $n \ge N$,
since $ \Phi_{t, N} $ is measurable with respect to $\F_N$, we deduce
from (\ref{u*}) and the absolute continuity (\ref{Q}) that
\[
(\ref{eq3b})_n = \Q_x^{(h)} \bigl[ 1_{\{ \tau_t^+ =n\} } f(S_i, 0\le i
\le n) g(w_n) \Phi_{t,N}\bigr].
\]

For $n < N$, we deduce from the branching property
along the stopping line $\cc_t$ (see Jagers \cite{J89}) that
\begin{eqnarray*} (\ref{eq3b})_n &=& \e_x \biggl[ \sum_{|u|=n}
1_{\{ \tau_t^+(u)=n\}} f\bigl( V(u_i), 0\le i \le n \bigr) g(u) \Phi_{t,N}
\sum
_{|v|=N, u < v} {h(v) \over h(\varnothing) \lambda^N } \biggr] \\
&=& \e_x \biggl[ \sum_{|v|=N}
1_{\{ \tau_t^+(v)=n\}} f\bigl( V(v_i), 0\le i \le n \bigr) g(v_n) \Phi_{t,N}
{h(v) \over h(\varnothing) \lambda^N } \biggr]\\
&=& \Q_x^{(h)} \bigl[ 1_{\{
\tau_t^+ =n\}} f(S_i, 0\le i \le n) g(w_n) \Phi_{t,N}\bigr],
\end{eqnarray*}
by using again (\ref{u*}) and the absolute continuity (\ref
{Q}) at $N$. Noting that $f(S_i, 0\le i \le n) g(w_n) = f( S_i, 0\le i
\le\tau^+_t) g(w_{\tau^+_t})$ on\vspace*{1pt} $\{\tau_t^+=n\}$, we take the sum
of $(\ref{eq3b})_n$ over all $n$ and obtain (\ref{eq3a}). The proof for
$\cL[a]$ works by analogy.
\end{pf*}

Let us present below a particular example of $h$ and the corresponding
laws of $(\Theta_n, S_n)_{n\ge0}$. Recall (\ref{t*}). Define
%
\begin{equation} \label{defhu12} h(u):= \cases{
e^{\varrho_* V(u)}, &\quad $\mbox{if $\psi'(\varrho_*)=0$,}$
\vspace*{2pt}\cr
e^{\varrho_+ V(u)}, & \quad $\mbox{if $\psi'(\varrho_*) < 0$,}$}\qquad  u\in\cT.
\end{equation}

Since $\psi(\varrho_*) =0$ in the critical case and
$\psi(\varrho_+)=0$ in the subcritical case, the function $h $
satisfies (\ref{h}) with $\lambda=1$ and $\H_+=\cT$. We mention
that in
the subcritical case, since $\psi(\varrho_-)=0$, the function $u \to
e^{ \varrho_- V(u)} $ also satisfies (\ref{h}) with $\lambda=1$. This
fact will be explored in Section~\ref{prsubcritical} for the definition
of $\Q^{(\varrho_-)}$, the measure satisfying (\ref{Q}) with $h(u)=e^{
\varrho_- V(u)} $.

Write for any
$x \in\r$, $\Q_x \equiv\Q_x^{(h)}$ the probability with the choice of
$h $ given in~(\ref{defhu12}). For simplification,
let
%
\begin{equation}\label{kappa} \varrho:= \cases{
\varrho_*, & \quad $\mbox{if $\psi'(\varrho_*)=0$ (critical case);}$
\vspace*{2pt}\cr
\varrho_+, &\quad $\mbox{if $\psi'(\varrho_*)<0$ (subcritical case).}$}
\end{equation}

Then for any $x\in\r$, $\Q_x$ satisfies
%
\begin{equation}
{d \Q_x \over d \mathbf{P}_x}\bigg|_{\F_n} =
e^{ - \varrho x} \sum_{|u|=n} e^{\varrho V(u)}.
\label{Q1}
\end{equation}

We shall write $\Q\equiv\Q_0$ when $x=0$. The following
description of the law of $(S_n, \Theta_n)_{n\ge0}$ under $\Q_x$ is an
easy consequence of Proposition
\ref{pchange-proba}(i).

\begin{corollary} \label{cmany1} Recall (\ref{kappa}) and
(\ref{Theta}). Fix $x \in\r$.
\begin{longlist}[(iii)]
\item[(i)] Under
$\Q_x$, $(S_n-S_{n-1}, \Theta_n)_{n\ge1}$ are i.i.d. under $\Q_x$
whose common law is determined by
\[
\Q_x \bigl[ f(S_n-S_{n-1}) e^{ - \langle g, \Theta_n \rangle}
\bigr]= \e
\biggl[ \sum_{|u|=1} e^{\varrho V(u)} f\bigl(V(u)\bigr) e^{ - \sum_{v \neq u,
|v|=1} g(V(v))}\biggr]
\]
for any $n \ge1$, any measurable functions $f, g\dvtx\r\to\r_+$. In
particular, the process $(S_n)_{n\ge0 }$ is a random walk on $\r$,
starting from $S_0=x$, with step distribution given by
%
\begin{equation} \label{S} \Q_x \bigl[ f(S_n -S_{n-1}) \bigr]= \e
\biggl[
\sum_{|u|=1}
f\bigl(V(u)\bigr) e^{\varrho V(u)} \biggr], \qquad n\ge1.
\end{equation}

\item[(ii)] For any $n \ge1$ and any measurable function $F\dvtx\r^{n+1}
\to\r
_+$,
\[
\e_x \biggl[ \sum_{ |u|=n} F\bigl( V(u_i), 0\le i \le n\bigr) \biggr]=e^{\varrho
x} \Q_x\bigl [ e^{- \varrho S_n} F( S_i, 0\le i \le n)
\bigr].
\]

\item[(iii)] For any $n \ge1$, and any $|u|=n$,
\[
\Q_x(w_n=u | \F_n) = {e^{\varrho V(u)} \over\sum_{|v|=n}
e^{\varrho V(v)} }.
\]
\end{longlist}
\end{corollary}


Remark that by (\ref{S}), $\Q[S_1]= 0$ and
$ \Q[S_1^2]= \psi^{\prime\prime}(\varrho_*)>0$ in the critical case,
while $\Q[S_1]= \psi'(\varrho_+)>0$ in the subcritical
case.

\subsection{Spinal decomposition for a killed branching random walk}

Before introducing a change of measure related to the killed
branching walk, we recall some elementary facts on the Palm
distribution of the point process $\L= \sum_{|u|=1}
\delta_{\{V(u)\}}$ under $\mathbf{P}$. Let $ \e( \L(dx))$ be the
intensity measure of $\L$, namely for any measurable function $f\dvtx
\r\to\r_+$,
\[
\int_\r f(x)\e\bigl( \L(dx)\bigr) = \e\biggl[ \sum_{|u|=1}
f\bigl(V(u)\bigr)\biggr].
\]
Clearly $ \e( \L(dx))$ is $\sigma$-finite since $\psi$ is well
defined on some interval. Then
there exists a family $(\Xi_x, x\in\r)$, called reduced Palm
distributions, of distributions of random point measures on $\r$
such that
%
\begin{eqnarray}\label{palm}
\int_{\Omega_f} F(x, \theta) \Xi_x(d \theta) = { \e[ F(x, \L-
\delta_{\{x\}}) \L(dx ) ] \over\e( \L(dx))}
\nonumber
\\[-8pt]
\\[-8pt]
\eqntext{\mbox{for } \e\bigl( \L(dx)\bigr) \mbox{-a.e. } x}
\end{eqnarray}
for any measurable $F\dvtx\r\times{\Omega_f}(\r) \to
\r_+$, and where ${\Omega_f} $ denotes the set of $\sigma$-finite
measures on $\r$; see Kallenberg \cite{K76}, Chapter~10 for more
details. Roughly said, $\Xi_x$ is the distribution of $\L-
\delta_{\{x\}}$ conditioned on that $\L$ charges $x$.

In this subsection, let $((S_n), \Q_x)$ be as in Corollary \ref
{cmany1} and (\ref{Q1}). Based on
Corollary \ref{cmany1}(i) (with $n=1$ and $x=0$), elementary
computations give that for any measurable $f, g\dvtx\r\to\r_+$,
\[
\Q\bigl[ f(S_1) e^{ - \langle g, \Theta_1 \rangle} \bigr]= \int_\r
\e\bigl( \L(dx)\bigr) e^{\varrho x} f(x) \int_{\Omega_f} e^{ -
\langle
g, \theta\rangle} \Xi_x(d \theta).
\]

It follows immediately from (\ref{S}) that the law of
$S_1$ under $\Q$ is given by $\Q(S_1 \in dx)=\e( \L(dx))
e^{\varrho x}$. Hence for any measurable $f, g\dvtx\r\to\r_+$,
%
\begin{equation}\label{palm2} \Q\bigl[ f(S_1) e^{ - \langle g,
\Theta_1
\rangle} \bigr]= \int_\r\Q(S_1 \in dx ) f(x)
\int_{\Omega_f} e^{ - \langle g, \theta\rangle} \Xi_x(d \theta).
\end{equation}

In words, $\Xi_x$ is the law of $\Theta_1$ conditioned on
$\{S_1=x\}$ under $\Q$.

Now, we are interested in a change of measure in the killed
branching random walk. To introduce the corresponding density, we consider
$R(\cdot)$ the renewal function of the random walk
$(S_n)_{n\ge0}$ under $\Q$. More precisely, for $x>0$, $R(x)$ is
defined by the expected
number (under $\Q$) of visits to $[-x,0]$ before first returning to
$[0,\infty)$, that is, $R(0)=1$, and
%
\begin{equation}\label{Rx}
R(x):= \Q\Biggl[\sum_{j=0}^{\tau^*-1} 1_{\{-x \le S_j\}} \Biggr]\qquad
\forall x >0,
\end{equation}
with $\tau^*:=\inf\{j\ge1\dvtx S_j \ge0\}$. We extend the
definition of $R$ on the whole real line by letting $R(x)=0$ for all
$x <0$.

Recall that $\Q[S_1]= 0$ in the critical case and $\Q[S_1] >0$ in the
subcritical
case. It is known (see Lemma \ref{lemRW}) that the following limits
exist and is equal to some positive constants:
%
\begin{equation}\label{defcr} C_R:= \cases{
\displaystyle\lim_{x \to\infty} { R(x) \over x} = {1\over\Q[ - S_{\tau^-_0}]}
, &\quad
$\mbox{if $\psi'(\varrho_*)=0$ (critical case),}$ \vspace*{2pt}\cr
\displaystyle\lim_{x\to\infty} R(x)= {1\over\Q( \tau^-_0=\infty)}, & \quad$\mbox{if
$\psi'(\varrho_*)<0$ (subcritical case),}$}\hspace*{-40pt}
\end{equation}
with $\tau^-_0$ defined in (\ref{tauS}). Recall (\ref
{kappa}). Define
\[
h_+(u):= R\bigl(V(u)\bigr) e^{\varrho
V(u)} 1_{\{V(u_1)\ge0, \ldots, V(u_n) \ge0\}},\qquad |u|=n, u
\in\cT.
\]

It is well known that $(R(S_n)1_{\{\tau_0^->n\}}, n\ge0)$ is
a $\Q_x$-martingale for any $x\ge0$; see, for example, \cite
{bertoin-doney}. Then $h_+$ satisfies (\ref{h}) with $\lambda=1$. Note
that in this case,
$\H_+=\{u \in\cT\dvtx V(u_0)\ge0, \ldots, V(u_{|u|}) \ge0\}= \zz
$ is
exactly the set of progeny of the killed branching walk.

Let $\Q_x^{+}$ be the probability satisfying (\ref{Q}) and
(\ref{u*}) with $h=h_+$,
%
\begin{equation} \label{Q+} {d \Q_x^{+}\over d \mathbf{P}_x}\bigg|_{\F
_n}:=
{e^{-\varrho x} \over R(x)}
\sum_{|u|=n, u \in\zz} R\bigl(V(u)\bigr)
e^{\varrho V(u)} =: M^*_n,\qquad x \ge0, n\ge1\hspace*{-35pt}
\end{equation}
with $M^*_0:=1$. Write for simplification
$\Q^{+}=\Q_0^{+}$. Recalling (\ref{Theta}), we deduce from
Proposition \ref{pchange-proba} the following result; see Figure~\ref{fig5} below.

\begin{figure}

\includegraphics{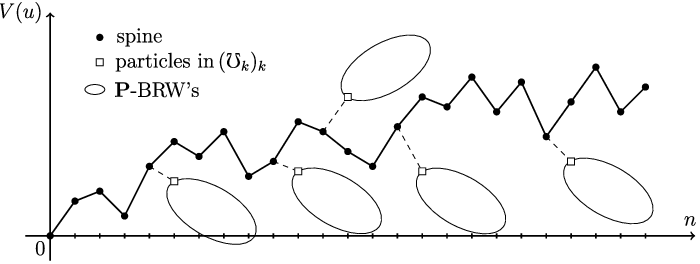}

\caption{Spinal decomposition under $\Q_0^+$.}\label{fig5}
\end{figure}

\begin{corollary} \label{cq+} Recall (\ref{kappa}). Fix $x\ge0$.
Under $\Q_x^{+}$:
\begin{longlist}[(a)]
\item[(a)] $(S_n)_{n\ge0}$ is a $(\G_n)$-Markov chain: for any $n \ge1$,
$y>0$ and any measurable function $f\dvtx\r_+ \to\r_+$,
\[
\Q_x^{+}\bigl[f(S_n) | \G_{n-1}, S_{n-1} =y \bigr] = \Q_y
\biggl[ {
R(S_1)\over R(y)} f(S_1) 1_{\{S_1\ge0\}}\biggr].
\]
In words, under $\Q_x^{+}$, the process $(S_n,
n\ge0)$ has the same law as the random walk $(S_n, n\ge0)$ under
$\mathbf{P}_x$, conditioned to stay nonnegative.

\item[(b)] Conditioned on $(S_n)_{n\ge0}$, the point processes $(\Theta_n)_{n
\ge1}$ are
independent, and each $\Theta_n$ is distributed as $\Xi_{ S_n-
S_{n-1}}$.

\item[(c)] For any nonnegative function $F$ and any $n \ge0$,
\[
\e_x \biggl[ \sum_{ u \in\zz, |u|=n} F\bigl( V(u_i), 0\le i \le|u|\bigr)
\biggr]=R(x) e^{\varrho x} \Q^{+}_x\biggl [ { e^{- \varrho
S_{n}}\over R(S_{n})} F( S_i, 0\le i \le n)
\biggr].
\]
\end{longlist}
\end{corollary}

\begin{pf}
The formula many-to-one
(c) is routine. Let us only check (a) and (b): By Proposition
\ref{pchange-proba}(i), we get that for any $n\ge1$,
%
\begin{eqnarray}&& \Q_x^{+}\bigl [ e^{- \langle g, \Theta_n
\rangle} f(S_n) | \G_{n-1}, S_{n-1} =y \bigr] \nonumber\\
&&\qquad= \e_y \biggl[ \sum_{ |u|=1} { 1 \over R(y) e^{\varrho y}}
e^{\varrho V(u)} R\bigl(V(u)\bigr)
\nonumber
\\[-8pt]
\\[-8pt]
\nonumber
&&\hspace*{71pt}{}\times 1_{\{V(u) \ge0\}} f\bigl(V(u)\bigr) e^{ - \sum_{v \not
=u, |v|=1}
g(V(v) -y)}\biggr] \label{c21} \\
&&\qquad= \Q_y \biggl[ {R(S_1) \over R(y)} 1_{\{S_1\ge0\}} f(S_1) e^{-
\langle g, \Theta_1 \rangle}\biggr] \label{c22} \\
&&\qquad= \Q\biggl[ {R(y+S_1) \over R(y)} 1_{\{y+S_1\ge0\}} f(y+S_1) e^{-
\langle g, \Theta_1 \rangle}\biggr], \label{c23}
\end{eqnarray}
by using Corollary \ref{cmany1}(i). Taking $g=0$ in
(\ref{c22}) yields the assertions in (a). Taking $n=1$ gives the joint
law of $(S_1, \Theta_1)$ under $\Q_x^{+}$. Let $p(dz)= \Q(S_1 \in
dz)$ be the law of $S_1$ under~$\Q$. Recall that $\Xi_z$ is the law of
$\Theta_1$ conditioning on $\{S_1=z\}$ under~$\Q$. Then for any event
$A \in
\G_{n-1}$, we deduce from (\ref{c23}) that
\begin{eqnarray*} &&\hspace*{-4pt}
\Q_x^{+} \bigl[ e^{- \langle g, \Theta_n \rangle} f(S_n) 1_A\bigr]
\\
&&\hspace*{-4pt}\quad= \Q_x^{+} \biggl[ 1_A \int_\r p(dz) {R(S_{n-1}+ z) \over
R(S_{n-1})} 1_{\{S_{n-1}+z \ge0\}} f(S_{n-1}+z) \int_{\Omega_f}
\Xi_{ z} (d \theta) e^{- \langle g, \theta\rangle} \biggr] \\
&&\hspace*{-4pt}\quad= \Q_x^{+} \biggl[ 1_A f(S_n) \int_{\Omega_f} \Xi_{S_n-
S_{n-1}} (d \theta) e^{- \langle g, \theta\rangle} \biggr],
\end{eqnarray*}
by using (a) for the last equality. This together with
Markov's property of $(S_n)$ with respect to $(\G_n)$, imply that
for any $n \ge1$ and $g\dvtx\r\to\r_+$,
\[
\Q_x^{+}\bigl [ e^{-
\langle g, \Theta_n \rangle} | \G_{n-1}, (S_j)_{j\ge0}
\bigr] =
\int_{\Omega_f} \Xi_{ S_n- S_{n-1}} (d \theta) e^{- \langle g,
\theta
\rangle},
\]
proving (b).
\end{pf}


%
%
%
%
%
%
%
%
%
%
%
%
%
%
%
%
%
%
%
%
%
%
%
%
%
%
%
%
%
%
%
%
%
%
%
%
%
%
%
%
%
%
%
%
%
%
%
%
%
%
%
%
%


\begin{remark}\label{rq+} If we assume that $\L= \sum_{i=1}^\nu
\delta_{\{X_i\}}$ with $(X_i)_{i\ge1}$ a sequence of i.i.d.
real-valued variables of the same law as $X$, independent of $\nu$,
then the expectation in
(\ref{c21}) equals
\[
\sum_{k\ge0} \mathbf{P}( \nu=k) k \e\biggl[{
R(X+y) \over R(y)} e^{ \varrho X} 1_{\{X+y \ge0\}} f(X+y) \biggr]
\bigl( \e e^{- g(X)}\bigr)^{k-1},
\]
which implies that under $\Q_x^{+}$ for each $n \ge1$, conditionally on
$\G_{n-1}$ and on $\{S_{n-1} =y\}$, $S_n$ and $\Theta_n$ are
independent and
$\Theta_n$ is distributed as $ \sum_{i=1}^{\widetilde\nu-1}
\delta_{X_i}$, with $\widetilde\nu$ the size-biased of $\nu$, $\Q
_x^{+} (\widetilde\nu=k)= k
\mathbf{P}(\nu=k)/\e[\nu], k\ge1$, and independent
of $(X_i)_{i\ge1}$.
\end{remark}

We may extend Corollary \ref{cq+} to the stopping lines. Remark that
$ \cc_t\cap\zz= \H(t)$; see (\ref{defcct}) and (\ref{defhl2}). We
deduce from Proposition \ref{pchange-proba2} the following result:

\begin{corollary} \label{cmany2} Recall (\ref{kappa}) and (\ref
{tauS}). Fix $0\le x <
t$. We have
%
\begin{equation} \label{c31bis} {d \Q_x^{+} \over d \mathbf{P}_x}
\bigg|_{\F_{\cc_t}}= {e^{-\varrho x} \over R(x)}
\sum_{ u \in\H(t) } R\bigl(V(u)\bigr)
e^{\varrho V(u)}=: M^*_{\cc_t}.
\end{equation}
\begin{longlist}[(ii)]
\item[(i)] Under probability $\Q_x^{+}$,
conditionally on $\G_{\cc_t}$ and on $\{V(v)=x_v, v\in\bigcup
_{k=1}^{\tau_t^+ } \I_k\}$, the shifted branching
random walks $\{ V(vu) - V(v)\}_{u \dvtx|vu|\le\tau_t^+(vu)} $,
stopped by the line $\cc_t$, are independent, and have the
same law as\break $\{V(u)\}_{|u|\le\tau_{t- x_v}^+(u)}$ under $\mathbf{P}_0$,
stopped by the line $\cc_{t-x_v}$.

\item[(ii)] Moreover, for any measurable function $F\dvtx\r^{\n_+}
\to\r_+$,
\begin{eqnarray*}
&&\e_x \biggl[ \sum_{ u \in\cc_t \cap\zz} F\bigl( V(u_i), 0\le i \le|u|\bigr)
\biggr]\\
&&\qquad=R(x) e^{\varrho x} \Q^{+}_x \biggl[ { e^{- \varrho
S_{\tau_t^+}}\over R(S_{\tau_t^+})} 1_{\{ \tau_t^+< \tau_0^-\}} F\bigl(
S_i, 0\le i \le\tau_t^+\bigr)
\biggr].
\end{eqnarray*}
\end{longlist}
\end{corollary}

\section{\texorpdfstring{Maximum of the killed branching random walk: Proofs of Theorem
\protect\ref{thmyaglom-critiq} and Proposition \protect\ref{pyaglom}}
{Maximum of the killed branching random walk: Proofs of Theorem 3 and Proposition 1}}
\label{secyaglom}

Let us first recall the following criterion for convergence in
distribution of point processes which can be found in Resnick \cite
{resnick}; see page 153, Proposition 3.19.
Let $E$ be a polish space. Then, let us define the Laplace transform of
a point process $\theta$ with probability measure $\P$ by
%
\begin{equation}
\Psi_{\P}(f):=\int\exp\biggl\{- \int f \,{d}\theta
\biggr\} \,{d}\P
(\theta)=\int\exp\bigl\{- \langle f, \theta\rangle\bigr\}
\,{d}\P
(\theta),
\end{equation}
where $f$ is a positive measurable function from $E$ to $\r$. Let
$C_K^+(E)$ be the space of continuous functions from $E$ to $\r_+$ with
compact support. Then we have
%
\begin{equation}
\lim_{n \to\infty} \Psi_{\P_n}(f)=\Psi_{\P}(f)\qquad \forall f
\in C_K^+(E),
\end{equation}
if and only if
%
\begin{equation}
\P_n \stackrel{\mathrm{(vague)}}{\longrightarrow} \P,\qquad n \to
\infty,
\end{equation}
which is the same as the convergence in distribution of the point processes.

Recall the real-valued random walk $(S_n)$ defined in Corollary \ref{cmany1}.
In order to treat both critical and subcritical cases in the same
proof, we introduce the following function defined on $\r_+$ by
%
\begin{equation} \label{def-Rt}\quad  {\mathscr R}(t):= \cases{
t, & \quad$\mbox{if $\psi'(\varrho_*)=0$,}$ \vspace*{2pt}\cr
1, & \quad$\mbox{if $\psi'(\varrho_*)<0,$}$}
\qquad
\varrho:= \cases{
\varrho_*, & \quad$\mbox{if $\psi'(\varrho_*)=0$,}$ \vspace*{2pt}\cr
\varrho_+, & \quad $\mbox{if $\psi'(\varrho_*)<0,$}$}
\end{equation}
and observe that the renewal function $R(\cdot)$, associated with the
random walk $(S_n, \Q) $ defined by (\ref{Rx}), satisfies [see
(\ref
{defcr})]
\[
R(t) \sim C_R {\mathscr R}(t),\qquad t \to\infty.
\]

We take the notation of Theorem \ref{thmyaglom-critiq} and Proposition
\ref{pyaglom}. The key step is to prove that for any $f \in C_K^+(\r)$
and when $t \to\infty$, we have
%
\begin{equation}
\label{eqkey-Yaglom}
\e_x \bigl[ e^{- \langle f, \mu_{\B,t} \rangle} 1_{\{H(t)>0\}}
\bigr] \sim\frac{R(x) e^{\varrho x}}{C_R {\mathscr R}(t)e^{\varrho
t}} \Q\biggl[ \frac{e^{- \langle f, \mu_{\B,\infty} \rangle
}}{\Re}
\biggr].
\end{equation}
%

We recall from (\ref{c31bis}) that $ M^*_{\cc_t}= {e^{-\varrho x}
\over
R(x)}
\sum_{ u \in\H(t)} R(V(u))
e^{\varrho V(u)}$, where $\H(t)$ denotes the set of those $u\in\zz$
satisfying $\tau^+_t(u) = \vert u \vert$ [see (\ref{defhl2})]. Then
$H(t)>0$ if and only if $ M^*_{\cc_t}>0$. It follows that
%
\begin{eqnarray}
\e_x \bigl[ e^{- \langle f, \mu_{\B,t} \rangle} 1_{\{H(t)>0\}}
\bigr] &=& \e_x \biggl[\frac{M^*_{\cc_t}}{M^*_{\cc_t}} e^{- \langle f,
\mu_{\B
,t} \rangle} 1_{\{H(t)>0\}}\biggr]
\nonumber
\\[-8pt]
\\[-8pt]
\nonumber
&=& \Q_x^+ \biggl[{e^{- \langle f,
\mu
_{\B,t} \rangle} \over M^*_{\cc_t}}\biggr].
\end{eqnarray}
We will now use the so-called ``decomposition along the spine'' $(w_k)$
(under $\Q_x^+$). Recalling that
$\I_k =\{ u \dvtx|u| =k, {\buildrel\leftarrow\over u} = w_{k-1},
x\neq w_k \} $, we have
%
\begin{equation}\qquad
\langle f, \mu_{\B,t} \rangle= f\bigl(T_t^+\bigr) 1_{\{ \beta_t(w_{\tau
_t^+})=\infty\}} + \sum_{1\le k\le\tau_t^+} \sum_{u \in\I_k}
1_{\{
\beta_t(u)=\infty\}} \bigl\langle f, \mu_{\B,t}^{(u)} \bigr\rangle,
\end{equation}
where $T_t^+= S_{\tau_t^+} - t$ denotes the overshoot of $S$ above the
level $t$ [see (\ref{defovershoot})], and for any $u \in\cT$ the point
process $\mu_{\B,t}^{(u)}$ is associated to the subtree $\cT^{(u)}$
(rooted at~$u$) of $\cT$ and defined by
%
\begin{equation}\label{defmutvu}\quad
\mu_{\B,t}^{(u)}:=
\sum_{{v \in\cT^{(u)} \cap\H_\beta(t)}} \delta_{\{V(v)- t\}},\qquad
\mu_{t}^{(u)}:=
\sum_{{v \in\cT^{(u)} \cap\H(t)}} \delta_{\{V(v)- t\}}.
\end{equation}
%
Recall that $R(s) \sim C_R {\mathscr R}(s)$ when $s \to\infty$. Since
$V(u) >t $ for all $u \in\H(t)$, we get that, under $\Q_x^+$,
%
\begin{eqnarray}
M^*_{\cc_t} \sim{e^{-\varrho x} \over R(x)} {C_R {\mathscr R}(t)
e^{\varrho t}} \sum_{u \in\H(t)} {\mathscr R}\biggl( 1+\frac{V(u)-t}{t}
\biggr) e^{\varrho(V(u)-t)},
\nonumber
\\[-8pt]
\\[-8pt]
 \eqntext{t \to\infty.}
\end{eqnarray}
Then repeating the spinal decomposition arguments for the above sum\break
$\sum_{u \in\H(t)} $ we obtain
%
\begin{equation}
\e_x \bigl[ e^{- \langle f, \mu_{\B,t} \rangle} 1_{\{H(t)>0\}}
\bigr] \sim\frac{R(x) e^{\varrho x}}{C_R {\mathscr R}(t)e^{\varrho
t}} \Q_x^+ \biggl[ \frac{I_\beta(t)}{J(t)} \biggr],
\end{equation}
with
\begin{eqnarray*}
I_\beta(t) &:=& \exp\biggl\{ - f\bigl(T_t^+\bigr) 1_{\{ \beta_t(w_{\tau
_t^+})=\infty\}} - \sum_{1\le k\le\tau_t^+} \sum_{u \in\I_k}
1_{\{
\beta_t(u)=\infty\}} \bigl\langle f, \mu_{\B,t}^{(u)} \bigr\rangle\biggr\},
\\
J(t)&:=& {\mathscr R}\biggl( 1+\frac{T_t^+}{t} \biggr) e^{\varrho
T_t^+} + \sum_{1\le k\le\tau_t^+} \sum_{u \in\I_k} \int
{\mathscr R}\biggl( 1+\frac{z}{t} \biggr) e^{\varrho z}
\mu_t^{(u)} ({d}z).
\end{eqnarray*}
Therefore, to prove (\ref{eqkey-Yaglom}) we only have to show that
%
\begin{equation}
\label{eqyaglom-critiq-resu1}
\lim_{t \to\infty} \Q_x^+ \biggl[ \frac{I_\beta(t)}{J(t)}
\biggr] = \Q
\biggl[ \frac{e^{- \langle f, \mu_{\B, \infty} \rangle}}{\Re}
\biggr].
\end{equation}

Note that $I_\beta(t) \in[0,1]$ and $J(t) \ge1$, hence
${I_\beta(t)\over J(t) }\in[0,1]$. Recalling the
convergence in law of the process
$(t-S_{\tau_t^+-j})_{0\le j \le K}$ for any fixed $K \ge1$ (see Lemma~\ref{lemlaw-convergence-reverse}), we will restrict the sums over $k$
in $I_\beta(t)$ and $J(t)$ to $k$'s between $\tau_t^+-K$ and $\tau
_t^+$. To this aim let us introduce $H^{u}(t)$ the number of
descendants of $u$ that reach $t$ before $0$ [with the convention
$H^{u}(t)=1$ if $V(u) > t$]. The following lemma ensures that with
probability close to $1$, $\sum_{1\le k\le\tau_t^+-K} \sum_{u \in
\I_k}
H^{u}(t) =0$ (the sum is $0$ if $\tau_t^+ \le K$):

\begin{lemma}\label{Lhut0} We have:
\begin{longlist}[(ii)]
\item[(i)] $\limsup_{K\to\infty} \limsup_{t\to\infty} \Q_x^+
(\sum
_{k=1}^{\tau_t^+-K} \sum_{u \in\I_k} H^{u}(t) \ge1 )=0$;

\item[(ii)] $\limsup_{K\to\infty} \limsup_{t\to\infty} \Q_x^+ (
\beta_t
( w_{\tau_t^+}) \le\tau_t^+-K )=0$.
\end{longlist}
\end{lemma}

\begin{pf}
See Section~\ref{pr-abc}.
\end{pf}

Notice that $\lim_{t\to\infty}\Q_x^+(\tau^+_t >K)=1$ and that on
$\{
\beta_t ( w_{\tau_t^+}) > \tau_t^+-K,\break \tau_t^+ >K\}$,
\[
\beta_t(u)= \inf\bigl\{ \tau_t^+-K < j \le|u|\dvtx \B(u_j) > e^{t-
V(u_{j-1})}\bigr\}
=: \beta_t^K(u)
\]
for any $u= w_{\tau_t^+}$ or $u \in\I_k$ with $\tau_t^+- K < k \le
\tau
_t^+$. The advantage of $\beta^K_t(u)$ is that $\beta_t^K(u)$ only
locally depends on the spines around $\tau_t^+$. Therefore (\ref
{eqyaglom-critiq-resu1}) will be a consequence of
%
\begin{equation}
\lim_{K \to\infty} \lim_{t \to\infty} \Q_x^+ \biggl[ \frac
{I_\beta
'(t,K)}{J'(t,K)} 1_{\{\tau^+_t >K\}} \biggr] = \Q\biggl[ \frac{e^{-
\langle f, \mu_{\B,\infty} \rangle}}{\Re} \biggr],
\end{equation}
with
\begin{eqnarray*}
I'_\beta(t,K) &:=& \exp\biggl\{ - f\bigl(T_t^+\bigr)1_{\{ \beta^K_t(w_{\tau
_t^+})=\infty\}} -\!\! \sum_{\tau_t^+-K < k \le\tau_t^+} \sum_{u \in
\I_k}
1_{\{ \beta^K_t(u)=\infty\}} \bigl\langle f, \mu_{\B,t}^{(u)} \bigr\rangle
\biggr\},
\\
J'(t,K)&:=& {\mathscr R}\biggl( 1+\frac{T_t^+}{t} \biggr) e^{\varrho
T_t^+} + \sum_{\tau_t^+-K < k \le\tau_t^+} \sum_{u \in\I_k} \int
{\mathscr R}\biggl( 1+\frac{z}{t} \biggr) e^{\varrho z}
\mu_t^{(u)} ({d}z).
\end{eqnarray*}

Recall from (\ref{defmutvu}) that the measures $ \mu_{\B,t}^{(u)}$ in
the previous expressions are associated with the branching random walk
killed at $0$. Now, we want to replace the measures $ \mu_{\B,
t}^{(u)}$ by the same measures $\widetilde\mu_{\B,t}^{(u)}$ but
associated with the nonkilled branching random walk,
%
\begin{eqnarray}
\widetilde\mu_{\B,t}^{(u)}&:=&
\sum_{{v \in\cT^{(u)} \cap\cc_t }} 1_{\{ \beta_t(v)=\infty\}}
\delta
_{\{V(v)- t\}},
\nonumber
\\[-8pt]
\\[-8pt]
\nonumber
 \widetilde\mu_t^{(u)}&:=&
\sum_{{v \in\cT^{(u)} \cap\cc_t }} \delta_{\{V(v)- t\}},
\end{eqnarray}
where we recall that $v \in\cT^{(u)} \cap\cc_t $ if and only if $v$
is a descendant of $u$ and $\tau_t^+(v)= \vert v\vert$; see (\ref
{defcct}) for the definition of $\cc_t$. The following lemma confirms
that we can replace $ \mu^{(u)}$ by $\widetilde\mu^{(u)}$ with
probability close to $1$:

\begin{lemma}
\label{lemnonkilled}
Let us define for $t>0$ and $K \ge1$,
\begin{eqnarray*}
\Gamma(t,K)&:=& \bigl\{ \tau_t^+ > K\bigr\} \\
&&{}\cap\bigl\{\bigl( \mu_{\B
,t}^{(u)}, \mu_{ t}^{(u)}\bigr) =\bigl(\widetilde\mu_{\B,t}^{(u)}, \widetilde
\mu_{ t}^{(u)}\bigr), \forall u \in\I_k, \forall k \in\bigl(
\tau_t^+-K, \tau_t^+\bigr] \bigr\}.
\end{eqnarray*}
Then for any $K \ge1$, we have
\[
\lim_{t \to\infty} \Q_x^+\bigl( \Gamma(t,K)^c\bigr)=0.
\]
\end{lemma}

\begin{pf}
See Section~\ref{pr-abc}. 
\end{pf}

By Lemmas \ref{Lhut0} and \ref{lemnonkilled}, to prove (\ref
{eqkey-Yaglom}) it is enough to
show that
%
\begin{equation}
\label{eqyaglomcritiq-resu2}
\lim_{K \to\infty} \lim_{t \to\infty} \Q_x^+ \biggl[ \frac
{\widetilde
I_\beta(t,K)}{\widetilde J(t,K)} 1_{\{\tau^+_t >K\}} \biggr] = \Q
\biggl[ \frac{e^{- \langle f, \mu_{\B, \infty} \rangle}}{\Re} \biggr],
\end{equation}
where $ \widetilde I_\beta(t,K)$ and $\widetilde J(t,K)$ are as $
I'_\beta(t,K)$ and $ J'(t,K)$ but with $\widetilde\mu^{(u)}$ in lieu
of~$ \mu^{(u)}$,
\begin{eqnarray*}
\widetilde I_\beta(t,K) &:=& \exp\biggl\{ - f\bigl(T_t^+\bigr)1_{\{ \beta
^K_t(w_{\tau_t^+})=\infty\}} -\!\!\sum_{\tau_t^+-K < k \le\tau_t^+}
\sum
_{u \in\I_k} 1_{\{ \beta^K_t(u)=\infty\}} \bigl\langle f, \widetilde\mu
_{\B,t}^{(u)} \bigr\rangle\biggr\},
\\
\widetilde J(t,K)&:=& {\mathscr R}\biggl ( 1+\frac{T_t^+}{t} \biggr)
e^{\varrho T_t^+} + \sum_{\tau_t^+-K < k \le\tau_t^+} \sum_{u \in
\I_k} \int{\mathscr R}\biggl( 1+\frac{z}{t} \biggr) e^{\varrho z}
\widetilde\mu_t^{(u)} ({d}z).
\end{eqnarray*}

Let us now introduce a family of point processes denoted by $(\overline
\mu_{\B, y}, \overline\mu_y)_{y\in\r}$, which are associated to the
nonkilled branching random walk $V$ under $\mathbf{P}$ and are defined by
%
\begin{equation}
\label{eqdef-mubar2}
\overline\mu_{\B,y}: = \cases{
\displaystyle\sum_{v \in\cc_y} 1_{\{ \beta_y (v)=\infty\}} \delta_{\{V(v)-y\}
}, &\quad
$\mbox{if $y \ge0$},$
\vspace*{2pt}\cr
\delta_{\{- y\}}, & \quad$\mbox{if $y < 0$,}$}
\end{equation}
and
%
\begin{equation}
\label{eqdef-mubar}
\overline\mu_y: = \cases{
\displaystyle\sum_{v \in\cc_y} \delta_{\{V(v)-y\}}, & \quad$\mbox{if $y \ge0,$}$
\vspace*{2pt}\cr
\delta_{\{- y\}}, & \quad$\mbox{if $y < 0$,}$}
\end{equation}
where $\cc_y$ was defined in (\ref{defcct}); in particular,
$\{V(v) - y, v \in\cc_y\}$ denotes exactly the set of overshoots of
the (nonkilled) branching random walk $V$ above the level $y$. By part
(i) of Corollary \ref{cmany2}, under $\Q^+$, conditionally on $\G
_{\cc
_t}$ and on $\{ V(u)=x_u, u \in\I_k, 1\le k \le\tau^+_t\}$, the
family $\{ (\widetilde\mu_{\B,t}^{(u)}, \widetilde\mu_t^{(u)}), u
\in\I_k, 1\le k \le\tau^+_t\}$ is independent and satisfies
%
\begin{equation} \label{eq518}
\bigl( \bigl(\widetilde\mu_{\B,t}^{(u)}, \widetilde\mu_t^{(u)} \bigr),
\mbox
{ under } \Q_x^+\bigr) \stackrel{\mathrm{law}}{=} \bigl( (\overline\mu
_{\B
, t-x_u}, \overline\mu_{t-x_u}), \mbox{ under } \mathbf{P}\bigr).
\end{equation}
For convenience of notation, let us introduce
%
\begin{eqnarray}
S^{(t)}_{i}&:=& S_{\tau_t^+} - S_{\tau_t^+-i},\qquad 1 \le i \le
\tau_t^+,
\\
\Theta^{(t)}_{i}&:=&\Theta_{\tau_t^+-i+1},\qquad 1\le i \le\tau_t^+.
\end{eqnarray}
Recall that $T_t^+:= S_{\tau_t^+} -t$ denotes the overshoot of $S$ over
$t$. Thus, (\ref{eq518}) yields that on $\{ \tau^+_t >K\}$,
%
\begin{equation} \label{qx521}
\Q_x^+ \biggl[ \frac{\widetilde I_\beta(t,K)}{\widetilde J(t,K)}
\Big| \G_{\cc_t} \biggr]
\stackrel{\mathrm{a.s.}}{=} \varphi_{t,K}\bigl(T_t^+,S^{(t)}_{1},\ldots
,S^{(t)}_{K},\Theta^{(t)}_{1},\ldots,\Theta^{(t)}_{K} \bigr),
\end{equation}
where for any $t_0>0, s_1, \ldots, s_K>0$ and the point measures
$\theta
^{(i)}$, $1\le i \le K$, of form $\theta^{(i)}=\sum_{j=1}^{m^{(i)}}
\delta_{x_j^{(i)}}$, we define
\[
\D_{i, K}:= \bigcap_{j=i}^K \bigl\{ \B\bigl(\theta^{(j)}\bigr) \le e^{- t_0+
s_j}\bigr\},\qquad 1\le i \le K
\]
and
\begin{eqnarray*}
&& \varphi_{t,K}\bigl(t_0,s_1,\ldots,s_K,\theta^{(1)},\ldots,\theta
^{(K)}\bigr)
\\
&&\qquad:= \e\biggl[\frac{\exp\{ - f(t_0) 1_{\D_{1, K}} - \sum_{i=1}^K
1_{\D_{i, K}}
\sum_{j=1}^{m^{(i)}} \langle f, \overline
\mu_{\B, s_i-t_0-x_j^{(i)}}^{(i,j)} \rangle\}}{
{\mathscr R}( 1+{t_0}/{t} ) e^{ \varrho t_0}
+\sum_{i=1}^K \sum_{j=1}^{m^{(i)}} \int{\mathscr R} (
1+{z}/{t} ) e^{\varrho z} \overline
\mu_{s_i-t_0-x_j^{(i)}}^{(i,j)} ({d}z)} \biggr],
\end{eqnarray*}
and with (under $\mathbf{P}$) $((\overline\mu_{\B, y}^{(i,j)},
\overline\mu
_{y}^{(i,j)}), y \in\r)_{i, j \ge1} $ i.i.d. copies of $((\overline
\mu
_{\B, y},\overline\mu_{ y}),\break y \in\r)$.
Then, applying part (b) of Corollary \ref{cq+} to (\ref{qx521})
implies that on\break $\{ \tau^+_t >K\}$,
%
\begin{equation}
\label{eq-prop-yag-eq1}
\Q_x^+ \biggl[ \frac{\widetilde I(t,K)}{\widetilde J(t,K)} \Big|
S_k, 0\le k \le\tau_t^+, \tau_t^+ \biggr]
\stackrel{\mathrm{a.s.}}{=} \widetilde\varphi_{t,K}
\bigl(T_t^+,S^{(t)}_{1},\ldots,S^{(t)}_{K}\bigr),
\end{equation}
with
\begin{eqnarray*}
\widetilde\varphi_{t,K}(t_0,s_1,\ldots,s_K)
:= \int\prod_{i=1}^K \Xi_{s_i- s_{i-1}}\bigl (d \theta^{(i)}\bigr)
\varphi
_{t,K}\bigl(t_0,s_1,\ldots,s_K,\theta^{(1)},\ldots,\theta
^{(K)}\bigr),
\end{eqnarray*}
with $s_0:=0$ for notational convenience.
Now for any $(t_0,s_1,\ldots,s_K) \in\r_+^* \times\r_+^K$ and for any
family $(\theta^{(i)})_{1\le i \le K}$ of point processes $\theta
^{(i)}:=\sum_{j=1}^{m^{(i)}} \delta_{x_j^{(i)}}$, let us define
\begin{eqnarray*}
\varphi_{\infty,K}\bigl(t_0,s_1,\ldots,s_K,\theta^{(1)},\ldots
,\theta
^{(K)}\bigr)&:=&\lim_{t \to\infty} \varphi_{t,K}
\bigl(t_0,s_1,\ldots
,s_K,\theta^{(1)},\ldots,\theta^{(K)}\bigr),
\\
\widetilde\varphi_{\infty,K}(t_0,s_1,\ldots,s_K
)&:=&\lim_{t
\to\infty} \widetilde\varphi_{t,K}(t_0,s_1,\ldots,s_K),
\end{eqnarray*}\eject\noindent
and observe that these limits exist by the dominated convergence
theorem, which also yields that
%
\begin{eqnarray}\label{defphitilde}
&& \widetilde\varphi_{\infty,K} (t_0,s_1,\ldots,s_K)
\nonumber\\
&&\qquad= \int \prod_{i=1}^K \Xi_{s_i-s_{i-1}} \bigl(d
\theta^{(i)}\bigr) \varphi_{\infty,K}\bigl(t_0,s_1,\ldots,s_K,\theta
^{(1)},\ldots,\theta^{(K)}\bigr)
\nonumber
\\[-8pt]
\\[-8pt]
\nonumber
&&\qquad= \int \prod_{i=1}^K \Xi_{s_i-s_{i-1}} \bigl(d \theta^{(i)}\bigr) \e
\\
&&\hspace*{34pt}\qquad{}\times\biggl[\frac{\exp\{ - f(t_0) 1_{\D_{1,K}}- \sum_{i=1}^K 1_{\D
_{i, K}}
\sum_{j=1}^{m^{(i)}} \langle f, \overline
\mu_{\B, s_i-t_0-x_j^{(i)}}^{(i,j)} \rangle\}}{
e^{ \varrho t_0}
+\sum_{i=1}^K \sum_{j=1}^{m^{(i)}} \int e^{\varrho z} \overline
\mu_{s_i-t_0-x_j^{(i)}}^{(i,j)} ({d}z)} \biggr]. \nonumber
\end{eqnarray}

The next step is to replace $\widetilde\varphi_{t,K}$ by $\widetilde
\varphi_{\infty,K}$:
%
\begin{lemma}
\label{lemvarphi-infinity} Fix $K \ge1$.
Then we have
%
\begin{eqnarray}
\label{eqlemvarinfinityeq0}
&&\lim_{t \to\infty} \Q_x^+ \bigl[ \bigl\vert\widetilde\varphi_{t,K}
\bigl(T_t^+,S^{(t)}_{1},\ldots,S^{(t)}_{K}\bigr)
\nonumber
\\[-8pt]
\\[-8pt]
\nonumber
 &&\hspace*{44pt}{}- \widetilde\varphi_{\infty
,K}\bigl(T_t^+,S^{(t)}_{1},\ldots,S^{(t)}_{K}\bigr) \bigr\vert1_{\{\tau^+_t >K\}}
\bigr] = 0.
\end{eqnarray}
\end{lemma}

\begin{pf}
See Section~\ref{pr-abc}.
\end{pf}

Note that since $\widetilde\varphi_{t,K} (\cdot)$ and $ \widetilde
\varphi_{\infty,K} (\cdot)$ differ only if $\psi'(\varrho_*)=0$, the
previous result is not trivial only in the critical case.

Finally thanks to (\ref{eq-prop-yag-eq1}) and Lemma \ref
{lemvarphi-infinity}, the double limits (\ref{eqyaglomcritiq-resu2})
will be a consequence of the following lemma.
%
\begin{lemma}
\label{lemcontinuity} We have
\[
\lim_{K \to\infty} \lim_{t \to\infty} \Q_x^+\bigl[\widetilde
\varphi
_{\infty,K}\bigl(T_t^+,S^{(t)}_{1},\ldots,S^{(t)}_{K}\bigr) 1_{\{\tau^+_t >K\}}
\bigr]
=\Q\biggl[ \frac{e^{- \langle f, \mu_{\B, \infty} \rangle
}}{\Re}
\biggr],
\]
where
%
\begin{eqnarray}\qquad
\mu_{\B, \infty}&:=& \delta_{U \hat S_{\hat\sigma}} 1_{\D_1}+
\sum
_{i= 1}^\infty1_{\D_i} \sum_{j=1}^{\widetilde\nu_i} \overline\mu
_{\B
, \hat S_{i} - U \hat S_{\hat\sigma}-\widetilde X_j^{(i)}}^{(i,j)},
\label{eqref-muinfty2}\\
\D_i&:=& \bigcap_{j=i}^\infty\bigl\{ \B(\widetilde\Theta_j) \le
e^{- U \hat S_{\hat\sigma} + \hat S_j } \bigr\}\qquad \forall
i\ge1,\\
\mu_{\infty}&:=& \delta_{U \hat S_{\hat\sigma}} + \sum_{i=
1}^\infty\sum_{j=1}^{\widetilde\nu_i} \overline\mu_{\hat S_{i} - U
\hat S_{\hat\sigma}-\widetilde X_j^{(i)}}^{(i,j)}, \label
{eqref-muinfty} \\
\Re&:=& e^{\varrho U \hat S_{\hat\sigma}} +\sum_{i= 1}^\infty
\sum_{j=1}^{\widetilde\nu_i} \int e^{\varrho z} \overline\mu
_{\hat S_{i} - U \hat S_{\hat\sigma}-\widetilde X_j^{(i)}}^{(i,j)}
({d}
z)=\int e^{\varrho z} \mu_{\infty}({d}z) \label{eqref-R}
\end{eqnarray}
and $\varrho= \varrho_*$ if $\psi'(\varrho_*)=0$, $\varrho=
\varrho_+$
if $\psi'(\varrho_*)<0$, and under $\Q$:
\begin{itemize}
\item the $((\overline\mu_{\B,y}^{(i,j)},\overline\mu
_{y}^{(i,j)}), y
\in\r)_{i, j \ge1} $ are i.i.d. with common distribution that of
$((\overline\mu_{\B, y},\overline\mu_{y}), y \in\r) $ under
$\mathbf{P}$
[see (\ref{eqdef-mubar})], and are independent of everything else; %
%
\item the process $(\hat S_n)_n$ (as well as the associated random time
$\hat\sigma$) and the random variable $U$ are introduced in Lemma
\ref
{lemlaw-convergence-reverse} (see Section~\ref{subsecretournement});
\item conditionally on $\{\hat S_n, n \ge0\}$, the random point
processes $\widetilde\Theta_i:=\sum_{j=1}^{\widetilde\nu_i} \delta
_{\widetilde X_j^{(i)}}$ for $i\ge1$ are independent, and $\widetilde
\Theta_i$ is distributed as $\Xi_{\hat S_{i-1}-\hat S_i}$; see (\ref
{palm}) and Corollary \ref{cq+} for the Palm measures $(\Xi_z, z\in
\r)$.
\end{itemize}
\end{lemma}

\begin{pf}
See Section~\ref{pr-abc}.
\end{pf}

\begin{pf*}{Proof of Theorem \ref{thmyaglom-critiq} and Proposition
\ref{pyaglom}}
Assembling (\ref{eq-prop-yag-eq1}), Lemmas~\ref{lemvarphi-infinity}
and~\ref{lemcontinuity} imply (\ref{eqyaglomcritiq-resu2}), hence
(\ref{eqkey-Yaglom}): namely for any $f \in C_K^+(\r)$ and when $t
\to
\infty$, we have
\begin{eqnarray*}
\e_x \bigl[ e^{- \langle f, \mu_{\B,t} \rangle} 1_{\{H(t)>0\}}
\bigr] &\sim&\frac{R(x) e^{\varrho x}}{C_R {\mathscr R}(t)e^{\varrho
t}} \Q\biggl[ \frac{e^{- \langle f, \mu_{\B,\infty} \rangle
}}{\Re
} \biggr] \\
&=& \frac{R(x) e^{\varrho x}}{C_R {\mathscr R}(t)e^{\varrho t}}
\Q\bigl[\Re^{-1}\bigr] \Q\bigl[ e^{- \langle f, \widehat\mu_{\B
,\infty
} \rangle} \bigr],
\end{eqnarray*}
by the definition of $\widehat\mu_{\B, \infty}$. Taking $f = 0$ in the
above asymptotical equivalence yields parts~(i) and~(ii) of Theorem
\ref{thmyaglom-critiq} while Proposition \ref{pyaglom} is a
consequence of parts~(i) and~(ii) together with (\ref
{eqkey-Yaglom}). Finally, taking $\B\equiv0$ in Proposition \ref
{pyaglom} gives part (iii), which completes the proof of Theorem \ref
{thmyaglom-critiq}.
\end{pf*}

\section{\texorpdfstring{Proof of Theorem \protect\ref{tleaves}: The critical case}
{Proof of Theorem 2: The critical case}}\label{prcritical}

We look at the critical case $\psi'(\varrho_*)=0$. By linear
transformation on $V$, we may assume that $\varrho_*=1$ in the whole
section without any loss of generality. We investigate the tail
distribution of the number of leaves $\#\cL[0]$; see (\ref{l[a]}) for
the definition. We will see that when $\cL[0]$ is large, the main
contribution comes from particles that reached a critical height $L$.
For integrability reasons, we will also restrict to \textit{good}
particles whose brothers do not display atypical jumps, and are not too
many. We denote by $\I(v):= \{ u\in\cT \dvtx{\buildrel\leftarrow
\over u}={\buildrel\leftarrow\over v}; u \neq v \}$ the set of
brothers of $v$ (${\buildrel\leftarrow\over v}$ denotes as before
the parent of $v$ in the tree $\cT$). For $\lambda>1, L>1$ (typically
$\lambda$ is a large constant whereas $L\to\infty$), we say that
%
\begin{eqnarray}\label{defbad} u \in\bb(L, \lambda) \mbox{ if there
exists some } 1\le j \le|u|\dvtx
\nonumber
\\[-8pt]
\\[-8pt]
\eqntext{\displaystyle\sum_{v \in\I(u_j) } \bigl(1+ e^{
\Delta V(v)
}\bigr) > \lambda e^{ {(L- V(u_{j-1}))/2}}}
\end{eqnarray}
and $u \in\gg(L, \lambda)$ if such $j$ does not exist. In
words, $\gg(L, \lambda)$ collects \textit{good} particles in the sense
that their large moments are finite; however, $\bb(L, \lambda)$ is a
set of \textit{bad} particles for which only low moments exist. Recall
from (\ref{zl}) that $Z[0, L]= \sum_u 1_{\{ \tau_0^-(u) = |u| < \tau
_L^+(u)\}}$ counts the number of leaves in the killed branching random
walk that did not touch the level $L$. Let us decompose $Z[0, L]$ as
the sums over good particles and bad particles,
\[
Z[0, L]= Z_g[0, L]+ Z_b[0, L]
\]
with
%
\begin{eqnarray}\label{defzgb} Z_g[0, L]&:=& \sum_{u \in\gg(L,
\lambda)}
1_{\{ \tau_0^-(u) = |u| < \tau_L^+(u)\}},
\nonumber
\\[-8pt]
\\[-8pt]
\nonumber
 Z_b[0, L]&:=& \sum
_{u\in
\bb(L, \lambda) } 1_{\{ \tau_0^-(u) = |u| < \tau_L^+(u)\}}.
\end{eqnarray}
The following lemma shows that we can neglect the number of \textit
{bad} particles.

\begin{lemma}\label{lbad}
For $\delta>0$ small enough, there exist constants $c=c(\delta)>0$ and
$c'=c'(\delta)>0$ such that for $x\ge0$, $\lambda\ge1$ and $L\ge1$,
%
\begin{equation}\label{eqbadL[0]}
\e_x\bigl[Z_b[0,L]\bigr] \le c\lambda^{-\delta} {(1+x)e^x \over
L^2} +
c e^x e^{-c' L}.
\end{equation}

For $\delta>0$ small enough, there exists a constant
$c=c(\delta)>0$ such that for $x\ge0$, $\lambda\ge1$, $L\ge1$ and
$B\ge0$,
%
\begin{equation}
\label{eqbadHL}\quad
\e_x\biggl[\sum_{u\in\H(L)} 1_{\{u \in\bb(L,\lambda)\}}Z^{(u)}[0,L
+B] \biggr] \le c\lambda^{-\delta}{1+B\over L+B}{(1+x)e^x\over L},
\end{equation}
where $Z^{(u)}[0,L+B]$ is the number of leaves in $\cL[0]$ that are
descendants of $u$ and did not cross level $L+B$.
\end{lemma}
\begin{pf}
We prove separately (\ref{eqbadL[0]}) and
(\ref
{eqbadHL}). The notation $c$ denotes a constant that can change
value from line to line.

\begin{pf*}{Proof of equation (\ref{eqbadL[0]})}
Mentioning here that (\ref{h}) holds with $\lambda=1$ (because $\psi
(\varrho_*) = \varrho_* \psi'(\varrho_*) =0$), Proposition \ref
{pchange-proba2} (applied to $\cL[0]$ and $h(u):=e^{V(u)}$) implies that
\begin{eqnarray*}
\e_x\bigl[Z_b[0,L]\bigr]
&=&
e^{x}\Q_x\biggl[{1\over\sum_{u\in\cL[0]} e^{V(u)}}Z_b[0,L]
\biggr]\\
&=&
e^{x}\Q_x\biggl[\sum_{u\in\cL[0]}{e^{V(u)}\over\sum_{u\in\cL[0]}
e^{V(u)}}e^{-V(u)}1_{\{\tau_0^-(u)<\tau_{L}^+(u)\}} 1_{\{ u \in\bb
(L,\lambda)\}}\biggr].
\end{eqnarray*}

The weight ${e^{V(u)}\over\sum_{u\in\cL[0]} e^{V(u)}}$ is
the probability that the vertex $u$ is the spine; see Proposition \ref
{pchange-proba2}. Therefore,
\[
\e_x\bigl[Z_b[0,L]\bigr]
=
e^x \Q_x\bigl[ e^{-S_{\tau_0^-} } 1_{\{\tau_0^-<\tau_{L}^+\}}1_{\{
w_{\tau_0^-}\in\bb(L,\lambda) \}} \bigr],
\]
where $\tau_0^-$ (resp., $\tau_L^+$) is the hitting time of
$(-\infty,0)$ [resp., $(L,+\infty)$] by the random walk $S$. Let
$\delta
\in(0,1)$, and, for $k\ge1$, $a_k:=\sum_{u\in\I_k} \{1+ e^{\Delta
V(u)}\}$ [we recall that $\I_k:=\I(w_k)$]. From the definition of
$\bb
(L,\lambda)$, we observe that
\[
1_{\{w_{\tau_0^-} \in\bb(L,\lambda)\}}
\le
\sum_{k=1}^{\tau_0^-} 1_{\{a_k > \lambda e^{(L-S_{k-1})/2}\}}
\le
\sum_{k=1}^{\tau_0^-} \min\bigl(a_k^{\delta} \lambda^{-\delta}
e^{-\delta(L- S_{k-1}) / 2},1\bigr).
\]

It follows that
%
\begin{equation}\label{eqbadL[0]1}
\e_x\bigl[Z_b[0,L]\bigr]
\le
e^x \Q_x\Biggl[e^{-S_{\tau_0^-}}1_{\{\tau_0^-<\tau_{L}^+\}} \sum
_{k=1}^{\tau_0^-}\min\bigl(a_k^{\delta}\lambda^{-\delta}
e^{-\delta(L-
S_{k-1}) / 2},1\bigr) \Biggr].\hspace*{-35pt}
\end{equation}

We first consider the term corresponding to $k=\tau_0^-$,
that is,
\begin{eqnarray*}
&&\Q_x\bigl[e^{- S_{\tau_0^-}}1_{\{\tau_0^-<\tau_{L}^+\}}\min
\bigl(a_{\tau_0^-}^{\delta} \lambda^{-\delta} e^{-\delta(L- S_{\tau_0^--1})
/2},1\bigr)\bigr] \\
&&\qquad\le
\Q_x\bigl[e^{- S_{\tau_0^-} }\min\bigl(a_{\tau_0^-}^{\delta}
\lambda
^{-\delta} e^{-\delta(L- S_{\tau_0^--1}) / 2 },1\bigr)\bigr].
\end{eqnarray*}

We know that $(S_n)_n$ is under $\Q$ a centered random walk
[since $\psi'(1)=0$]. Assumption (\ref{hypcrit}) ensures that $\Q
[e^{-(1+\eta)S_1}]$ is finite if $\eta$ is small enough. In turn,
it implies [see (\ref{ta-})] that
\[
\Q_x\bigl[e^{-(1+\eta) S_{\tau_0^-}}\bigr]\le c
\]
for small $\eta>0$, and any $x\ge0$. We also have $\Q_x[e^{S_{\tau
_0^--1} - S_{\tau_0^-}}]\le c$ by (\ref{diftau0}). Then it is not hard
to see that, with $\mathcal E:=\{ S_{\tau_0^-} \ge- \delta L/8,
S_{\tau
_0^--1} \le L/2\}$, we have, for some constant $\eta'>0$,
\[
\Q_x\bigl[e^{- S_{\tau_0^-} }1_{\mathcal E^c}\bigr]\le c' e^{-\eta'\delta L}.
\]
Therefore, we can restrict to the event $\mathcal E$, on which $e^{-
S_{\tau_0^-} } \le e^{\delta L/8}$, and $e^{-\delta(L- S_{\tau
_0^--1})/ 2} \le e^{-\delta L/4}$. It yields that
\[
\Q_x\bigl[e^{- S_{\tau_0^-}}\min\bigl(a_{\tau_0^-}^{\delta}
\lambda
^{-\delta}e^{-\delta(L- S_{\tau_0^--1} )/2},1\bigr)\bigr]
\le
c' e^{-\eta'\delta L} + \lambda^{-\delta}e^{-\delta{L/8}}\Q_x\bigl[
a_{\tau_0^-}^{\delta}\bigr].
\]

Observe that
\[
\Q_x\bigl[ a_{\tau_0^-}^{\delta}\bigr]= \sum_{j=1}^\infty\Q_x\bigl[ 1_{\{
j-1 <
\tau_0^-\}} \Q_{S_{j-1}}\bigl [ 1_{\{S_1<0\}} a_{j}^{\delta}\bigr] \bigr],
\]
by Markov's property at $j-1$. For $y:= S_{j-1} \ge0$,
\[
\Q_{y} \bigl[ 1_{\{S_1<0\}} a_{j}^{\delta}\bigr] \le\Q_{y} \bigl[e^{- ({1/2}) S_1
} a_{j}^{\delta} \bigr] = e^{- {(1/2)} y} \Q\bigl[e^{- {(1/2)} S_1 }
a_{j}^{\delta} \bigr].
\]

By the Cauchy--Schwarz inequality and (\ref{hypcrit}), we
have $\Q[e^{- S_1/2 } a_{j}^{\delta} ]\le c$ if $\delta>0$ is chosen
small enough. Therefore,
\[
\Q_x\bigl[ a_{\tau_0^-}^{\delta}\bigr] \le c \sum_{j=1}^\infty\Q_x
\bigl[ 1_{\{
j-1 < \tau_0^-\}} e^{- {(1/2)} S_{j-1}} \bigr],
\]
which is uniformly bounded by (\ref{foot3}). Hence, we
showed that
%
\begin{equation}\label{eqbadL[0]2}
\Q_x\bigl[e^{-S_{\tau_0^-} }1_{\{\tau_0^-<\tau_{L}^+\}}a_{\tau
_0^-}^{\delta} \lambda^{-\delta} e^{-\delta(L- S_{\tau_0^--1} ) /
2}\bigr] \le c e^{-\eta'' \delta L}.
\end{equation}

We consider now the terms corresponding to $k<\tau_0^-$ in
(\ref{eqbadL[0]1}). By Markov's property at time $k$, we get
\begin{eqnarray*}
&& \Q_x\bigl[e^{-S_{\tau_0^-} }1_{\{k< \tau_0^- < \tau_{L}^+ \}
}a_k^{\delta}\lambda^{-\delta} e^{-\delta(L- S_{k-1})/ 2} \bigr]
\\
&&\qquad\le
\lambda^{-\delta} \Q_x\bigl[ 1_{\{k< \tau_0^- < \tau_{L}^+ \}}
a_k^{\delta} e^{-\delta(L- S_{k-1})/ 2}\bigr]\sup_{y\ge0}\Q
_y\bigl[e^{-S_{\tau_0^-} }\bigr] \\
&&\qquad= c\lambda^{-\delta} \Q_x\bigl[ 1_{\{k < \tau_0^- < \tau_{L}^+
\}}
a_k^{\delta} e^{-\delta(L- S_{k-1})/ 2}\bigr],
\end{eqnarray*}
again by (\ref{ta-}). By Markov's property at time $k-1$, we
observe that the last expectation is $\Q_x[ 1_{\{k < \tau_0^- <
\tau_{L}^+ \}} e^{-\delta(L- S_{k-1})/ 2}]\Q[a_1^{\delta}]$.
Summing over $k\ge1$, we deduce that
\[
\Q_x\Biggl[e^{- S_{\tau_0^-} }1_{\{\tau_0^-<\tau_{L}^+\}}\! \sum
_{k=1}^{\tau_0^- - 1}\! a_k^{\delta} e^{-\delta(L- S_{k-1})/ 2}\Biggr]
\le
c \Q_x\Biggl[1_{\{\tau_0^-<\tau_{L}^+\}} \!\sum_{k=1}^{\tau_0^- - 1}\!
e^{-\delta(L- S_{k-1})/ 2}\Biggr].
\]

By (\ref{eal}), we have $\Q_x[1_{\{\tau_0^-<\tau
_{L}^+\}
}\sum_{k=0}^{\tau_0^- - 1} e^{-\delta(L- S_{k-1}) / 2}] \le c
{1+x \over L^2}$ for some $c=c(\delta)$. We obtain that
\[
\lambda^{-\delta}\Q_x\Biggl[e^{- S_{\tau_0^-} }1_{\{\tau_0^-<\tau
_{L}^+\}
} \sum_{k=1}^{\tau_0^- - 1} a_k^{\delta} e^{-\delta(L- S_{k-1}) / 2}
\Biggr]
\le
c' \lambda^{-\delta}{1+x \over L^2}.
\]

Then (\ref{eqbadL[0]}) follows from equations (\ref
{eqbadL[0]1}) and (\ref{eqbadL[0]2}).
\end{pf*}

\begin{pf*}{Proof of equation (\ref{eqbadHL})}
By the branching property, we have
\[
\e_x\biggl[\sum_{u\in\H(L)} 1_{\{u \in\bb(L,\lambda)\}
}Z^{(u)}[0,L +B]
\biggr]
=
\e_x\biggl[\sum_{u\in\H(L)} 1_{\{u \in\bb(L,\lambda)\}}
f\bigl(V(u)\bigr)\biggr],
\]
with $f(y):=\e_{y}[Z[0,L +B]]$. Using the measure $\Q_y$,
Proposition \ref{pchange-proba2} implies that $f(y)=e^y\Q
_y[e^{-V(w_{\tau_0^-})}1_{\{\tau_{0}^-<\tau_{L+B}^+\}}]$ which is
smaller than $c {1+(L+B-y)_+\over L+B}e^y$ by~(\ref{pal}). It follows that
%
\begin{eqnarray}\label{eqbadHL1}
&&\e_x\biggl[\sum_{u\in\H(L)} 1_{\{u \in\bb(L,\lambda)\}}Z^{(u)}[0,L+B]
\biggr]
\nonumber
\\[-8pt]
\\[-8pt]
\nonumber
&&\qquad\le
{c(1+B)\over L+B}\e_x\biggl[\sum_{u\in\H(L)} 1_{\{u \in\bb
(L,\lambda)\}
}e^{V(u)} \biggr].
\end{eqnarray}

By Proposition \ref{pchange-proba2} with $\cc_{L}$ and
$h(y):=e^{y}$, we observe that
\begin{eqnarray*}
\e_x\biggl[ \sum_{u\in\H(L)} 1_{\{u \in\bb(L,\lambda)\}}e^{V(u)}
\biggr]
=
e^x \Q_x[1_{\{\tau_{L}^+<\tau_0^-\}}1_{\{w_{\tau_{L}^+} \in
\bb
(L,\lambda)\}} ]
.
\end{eqnarray*}
As before, we have for $\delta>0$,
\[
1_{\{w_{\tau_{L}^+} \in\bb(L,\lambda)\}}
\le\lambda^{-\delta}\sum_{k=1}^{\tau_{L}^+} a_k^{\delta}
e^{-\delta
(L- S_{k-1}) / 2},
\]
where $a_k:= \sum_{u\in\I_k} \{1+e^{\Delta V(u)}\}$. Hence,
\begin{eqnarray*}
\Q_x[1_{\{\tau_{L}^+<\tau_0^-\}}1_{\{w_{\tau_{L}^+} \in\bb
(L,\lambda)\}} ]
&\le&
\lambda^{-\delta}\Q_x\Biggl[ 1_{\{\tau_{L}^+<\tau_0^-\}}\sum
_{k=1}^{\tau
_{L}^+} a_k^{\delta}e^{-\delta(L- S_{k-1}) / 2}\Biggr]\\
&=&
\lambda^{-\delta}\sum_{k\ge1} \Q_x\bigl[ 1_{\{k\le\tau
_{L}^+<\tau
_0^-\}}a_k^{\delta} e^{-\delta(L- S_{k-1}) / 2}\bigr].
\end{eqnarray*}
Using Markov's property at time $k-1$, for every $k\ge1$, yields
\[
\Q_x[1_{\{\tau_{L}^+<\tau_0^-\}}1_{\{w_{\tau_{L}^+} \in\bb
(L,\lambda)\}} ]
\le
c'\lambda^{-\delta}\Q_x\Biggl[ 1_{\{\tau_{L}^+<\tau_0^-\}}\sum
_{k=1}^{\tau_{L}^+} e^{-\delta(L- S_{k-1}) / 2}\Biggr],
\]
with $c'=\Q[a_1^{\delta}]<\infty$ if $\delta>0$ is small
enough by (\ref{hypcrit}). We get by equation~(\ref{a+1})
\[
\Q_x[1_{\{\tau_{L}^+<\tau_0^-\}}1_{\{w_{\tau_{L}^+} \in\bb
(L,\lambda)\}} ]
\le
c \lambda^{-\delta} {1+x \over L}.
\]
Going back to (\ref{eqbadHL1}), we obtain
\begin{eqnarray*}
\e_x\biggl[\sum_{u\in\H(L)} 1_{\{u \in\bb(L,\lambda)\}
}Z^{(u)}[0,L +B]
\biggr]
\le
c\lambda^{-\delta}{(1+B)\over L+B}{(1+x)e^x \over L},
\end{eqnarray*}
proving (\ref{eqbadHL}).
\end{pf*}
\noqed\end{pf}

We are going to re-prove the following estimate of A\"{i}d\'ekon \cite
{aidekon} but in a more general setting.
We recall that $\cL[0,L]$ is the set of leaves in $\cL[0]$ which did
not hit $(L,+\infty)$, and $Z[0,L]:=\#\cL[0,L]$. We call similarly
$\cL
_{g}[0,L]$ the leaves in $\cL[0,L]$ which are in $\gg(L,\lambda)$,
hence we have $Z_{ g}[0,L]:=\#\cL_{g}[0,L]$ the number of such leaves.

\begin{lemma}\label{L2moment} Fix $\lambda\ge1$ and assume that
$\psi
'(\varrho_*)=0$ with $\varrho_*=1$. Under (\ref{hypcrit}), there exists
some constant $c>0$ such that for all $L\ge1$, and $0\le x \le L$,
\[
\e_x \bigl[\bigl (Z_{g}[0,L]\bigr)^2\bigr] \le c\lambda(1+x)e^x { e^{L}\over
{L}^3}.
\]
\end{lemma}

\begin{pf}
Writing $Z_{g}[0,L] = \sum_{v\in\cL
[0]}e^{V(v)}1_{\{\tau_{L}^+(v)>|v|\}}{e^{-V(v)}}1_{\{v \in\gg(L,
\lambda)\}}$, we deduce from Proposition \ref{pchange-proba2} (applied
to $\cL[0]$ and $h(u):=e^{V(u)}$) that
%
\begin{equation}\label{eqboundZOL}
\e_x\bigl[\bigl(Z_{g}[0,L]\bigr)^2\bigr]
= e^{x}\Q_x\bigl[Z_{g}[0,L]e^{-S_{\tau_0^-} } 1_{\{ \tau_0^-<\tau
_{L}^+\}}1_{\{w_{\tau_0^-} \in\gg(L,\lambda)\}}\bigr].
\end{equation}

We decompose $Z_{g}[0,L]$ along the spine $(w_n, n\ge0)$
as follows:
\[
Z_{g}[0,L] \le1 + \sum_{k=1}^{\tau_0^-} \sum_{u\in\I_k} Z^{(u)}[0,L],
\]
where $Z^{(u)}[0, L]:= \sum_{v \in\cT^{(u)}} 1_{\{\tau
_0^-(v) = |v| < \tau_{L}^+(v)\}} $ denotes the number of descendants of
$u$, touching
$0$ before $L$ ($\cT^{(u)}$ means as before the subtree rooted at~$u$).
We have an inequality here since we dropped the condition that the
particles must be in $\gg(L,\lambda)$. By Proposition \ref
{pchange-proba}, under $\Q_x$,
conditioned on $\G_\infty:=\sigma\{\omega_j, S_j, \I_j, (V(u),
u\in\I
_j), j\ge0\}$, $(Z^{(u)}[0, L])_{u \in\I_{j}, j \le\tau_0^-}$ are
independent, and each $Z^{(u)}[0, L]$ is distributed as $(Z[0, L],
\mathbf{P}_{V(u)})$. In particular,
\[
\Q_x\bigl[Z_{g}[0,L] | \G_{\infty} \bigr] \le1 + \sum_{k=1}^{\tau_0^-}
\sum
_{u\in\I_k} \e_{V(u)}\bigl[Z[0,L]\bigr].
\]
Proposition \ref{pchange-proba2} implies as well that for
any $z\in\r$,
\[
\e_z\bigl[Z[0,L]\bigr] = e^z \Q_z\bigl[e^{-S_{\tau_0^-} } 1_{\{ \tau
_0^-<\tau
_{L}^+\}}\bigr],
\]
which is zero if $z >L $ and is $1$ if $z <0$. By (\ref
{pal}), we get that
\[
\e_z\bigl[Z[0,L]\bigr] \le c e^z{L-z+1\over L}1_{\{z\in[0,L]\}} + 1_{\{z<0\}}.
\]
Hence,
\[
\Q_x\bigl[Z_{g}[0,L] | \G_\infty\bigr] \le1+ \sum
_{k=1}^{\tau
_0^-} \sum_{u\in\I_k}\biggl ( c e^{V(u)}{L-V(u)+1 \over L}1_{\{
V(u)\in
[0,L]\}} + 1_{\{V(u)<0\}}\biggr).
\]
For $k<\tau_{L}^+$, we observe that [recalling
$S_{k-1}=V(w_{k-1}) \le L$]
\begin{eqnarray*}
&& \sum_{u\in\I_k} e^{V(u)}{L-V(u)+1 \over L}1_{\{V(u)\in[0,L]\}}\\
&&\qquad=
e^{S_{k-1}}\sum_{u\in\I_k} e^{\Delta V(u) }{L-V(u)+1 \over L}1_{\{
V(u)\in[0,L]\}}\\
&&\qquad\le
{L-S_{k-1} +1 \over L}e^{S_{k-1}} a_{k},
\end{eqnarray*}
with $a_k:= \sum_{u\in\I_k} \{1+e^{\Delta V(u)}\}$. If
$w_{\tau_0^-}\in\gg(L,\lambda)$, it follows that for any $k<\tau_0^-$,
\[
 \sum_{u\in\I_k} e^{V(u)}{L-V(u)+1 \over L}1_{\{V(u)\in[0,L]\}}
\le
\lambda e^{L}{L- S_{k-1} +1 \over L} e^{( S_{k-1} - L)/ 2}.
\]
Similarly, we observe that $\sum_{u\in\I_k} 1_{\{V(u)<0\}}
\le
a_{k} \le\lambda e^{{L /2}}$. Therefore, if $w_{\tau_0^-}\in\gg
(L,\lambda)$, then
\[
\Q_x\bigl[Z_{g}[0,L] | \G_\infty\bigr] \le1+ c \lambda{e^{L}
\over L} \sum_{k=1}^{\tau_0^-} (L- S_{k-1} +1) e^{( S_{k-1} - L)/ 2}.
\]
The equality (\ref{eqboundZOL}) implies that
\begin{eqnarray*}
&&\e_x\bigl[\bigl(Z_{g}[0,L]\bigr)^2\bigr]\\
&&\qquad\le
e^x\Q_x\bigl[e^{-S_{\tau_0^-} } 1_{\{ \tau_0^-<\tau_{L}^+\}}\bigr]
\\
&&\quad\qquad{}+ c\lambda{e^{x+L}\over L}\Q_x\Biggl[e^{-S_{\tau_0^-} } 1_{\{ \tau
_0^-<\tau_{L}^+\}} \sum_{k=1}^{\tau_0^-} (L- S_{k-1} + 1)e^{( S_{k-1}
-L)/ 2} \Biggr].
\end{eqnarray*}
The right-hand side is smaller than $e^x (1+ c'(1+x)\lambda
{e^{L}\over L^3})$ by (\ref{eal}). It completes the proof of the lemma.
\end{pf}

We look now at the progeny of a particle which went far to the right.
Let the derivative martingale be defined by
\[
\partial W_n:= - \sum_{|u|=n}
V(u) e^{V(u)},\qquad n\ge0.
\]
According to Theorems 5.1 and 5.2 in Biggins and Kyprianou
\cite{BigKyp04}, under $\mathbf{P}$, $\partial W_n$ converges almost
surely to
$\partial
W_\infty$ which has infinite mean and is almost surely positive on
$\{\cT=\infty\}$.

\begin{lemma} \label{Lconjecture2} Assuming $\psi'(\varrho_*)=0$ with
$\varrho_*=1$. Under (\ref{hypcrit}), as $t \to\infty$, the law of
$\#
\cL[0]$ under
$\mathbf{P}_t$, normalized by
$e^t/t$ converges in distribution to $c^{*} \partial W_\infty$,
with
%
\begin{equation} \label{eqconstant1} c^{*}:= { \Q[ e^{-
S_{\tau_0^-}} -1 ] \over-\Q[ S_{\tau_0^-} ] }.
\end{equation}
\end{lemma}

\begin{pf}
By linear translation, it is equivalent to prove
that under
$\mathbf{P}_0$, $\#\cL[-t]$ normalized by
$e^t/t$ converges in law to $c^*
\partial W_\infty$. If we define
\[
\partial W_{ \cL[-t]}:= -
\sum_{u \in\cL[-t]} V(u) e^{V(u)},
\]
then $ \partial W_{\cL[-t]}$ converges almost surely to
$\partial W_\infty$; cf. Biggins and Kyprianou \cite{BigKyp04},
Theorem 5.3. We write
%
\begin{equation}\label{eqnermanderiv}
\partial W_{\cL[-t]} = te^{-t}\biggl ( \sum_{u \in\cL[-t]} e^{V(u)+t} +
{1\over t}\eta_t\biggr),
\end{equation}
with $\eta_t = -\sum_{u \in\cL[-t]}(V(u)+t) e^{V(u)+t}$. At
this stage, we may apply a result of Nerman \cite{nerman81} for
the asymptotic behavior of ${ 1\over\# \cL[-t]}\sum_{u \in
\cL[-t]} e^{V(u)+t} $:\vspace*{1pt} Let $\xi:= \sum_{u \in\cL[0]}
\delta_{\{-V(u)\}}$ be the point process formed by the (nonkilled)
branching walk $V$ stopped at the line $\cL[0]$. Generate a
branching random walk $(V_\xi(u), u \in\cT_\xi)$ from the point
process $\xi$, where $V_\xi, \cT_\xi$ are related to $\xi$ in the
same way as $ V, \cT$ are to~$\L$. Define $ \cL_\xi[a]:= \{ u\in
\cT_\xi\dvtx|u|=\tau_a^+(u) \}$ for all $a> 0$. Clearly $
\cL_\xi[t]= \cL[-t] $ and
\[
{ \sum_{u \in\cL[-t]} e^{V(u)+t}
\over\# \cL[-t]} = { \sum_{u \in\cT_\xi} \psi_u( t - \sigma_u)
\over\sum_{u \in\cT_\xi} \phi_u( t - \sigma_u)},
\]
where for any $u \in\cT_\xi$, $\sigma_u:= - V_\xi(u )$ and
\[
\psi_u(x):= 1_{\{ x\ge0\}} \sum_{ {\buildrel
\leftarrow\over v} =u} e^{x- (\sigma_v- \sigma_u)} 1_{\{\sigma_v-
\sigma_u > x\}}, \qquad\!\! \phi_u(x):= 1_{\{ x\ge0\}} \sum_{
{\buildrel
\leftarrow\over v} =u} 1_{\{\sigma_v- \sigma_u > x\}}.
\]
Applying Theorem 6.3 in Nerman \cite{nerman81} (with $\alpha=1$
and $\lambda_u=\infty$ there) gives that conditioned on $ \{ \cT
=\infty
\}$, almost surely, when $t$ tends to infinity
\[
{ \sum_{u \in\cT_\xi} \psi_u( t - \sigma_u) \over\sum_{u \in
\cT_\xi} \phi_u( t - \sigma_u)} \to { \e[ \sum
_{|v|=1,v\in\cT
_\xi}
e^{- \sigma_v} \sigma_v ] \over\e[ \sum_{|v|=1,v\in\cT
_\xi}
(1- e^{-
\sigma_v}) ]}.
\]

Observe that $ \e[ \sum_{|v|=1,v\in\cT_\xi} e^{-
\sigma_v}
\sigma_v ] = - \e[ \sum_{u \in\cL[0]} e^{V(u)} V(u)]
= - \Q[ S_{\tau_0^-} ]$, and similarly, $\e
[ \sum_{|v|=1,v\in\cT_\xi} (1- e^{- \sigma_v}) ]= \Q
[ e^{-
S_{\tau_0^-}} ] -1$. Therefore conditioned on $ \{ \cT=\infty\}$,
almost surely
\[
{ \sum_{u \in\cL[-t]} e^{V(u)+t}
\over\# \cL[-t]} \to {\Q[ S_{\tau_0^-} ] \over1-\Q
[ e^{-
S_{\tau_0^-}} ]},\qquad t \to\infty.
\]

On the other hand, following the same strategy, we get that
conditioned on $\{\cT= \infty\}$, we have almost surely
\[
{ \eta_t
\over\# \cL[-t]} \to {\Q[ (S_{\tau_0^-})^2/2 ] \over
\Q
[e^{-S_{\tau_0^-}} ] -1},\qquad t \to\infty.
\]

Dividing both sides of (\ref{eqnermanderiv}) by $\#\cL[-t]$,
and using the fact that $\partial W_{\cL[-t]} $ goes to $\partial
W_{\infty}$, we deduce the lemma.
\end{pf}

We also need the following simple technical lemma whose proof is
postponed until Section~\ref{proofsRW}:
%
\begin{lemma} \label{Lconvolution} On some probability space $(\Omega,
{\mathscr F}, \ppp)$, let $\sum_{i=1}^\xi\delta_{\{Y_i\}}$ be a point
process on $\r_+$. Let
$(\Gamma_i, i\ge1)$ be a sequence of i.i.d. random variables on $\r_+$,
independent of $\sigma\{\xi, Y_i, 1\le i \le\xi\}$. Assume that for
some
$p>0$ and $a>0$,
\[
\ppp( \Gamma_1 > t ) = \bigl(a+o(1)\bigr)
t^{-p},\qquad t \to\infty.\vspace*{-12pt}
\]
\begin{longlist}[(ii)]
\item[(i)] If $p=1$ and if there exists some $\delta>0$ such that $ \eee
[ \sum_{i=1}^\xi Y_i^{1+ \delta}]< \infty$, then
\[
\lim_{t\to\infty} t \ppp\Biggl(
\sum_{i=1}^\xi Y_i \Gamma_i > t \Biggr) = a \eee\Biggl[
\sum
_{i=1}^\xi Y_i \Biggr].
\]

\item[(ii)] If $p>1$ and if there exists some $\delta>0$ such that $ \eee
[
\sum_{i=1}^\xi(1+Y_i) ]^{ p+ \delta} < \infty, $ then
\[
\lim_{t \to\infty} t^p \ppp\Biggl(
\sum_{i=1}^\xi Y_i \Gamma_i > t \Biggr) = a \eee\Biggl[
\sum
_{i=1}^\xi Y_i^p\Biggr].
\]
\end{longlist}
\end{lemma}

In the critical case, the branching random walk goes to $-\infty$. In
particular, almost surely, $H(L)=0$ if $L$ is large enough. Fix
$\lambda
\ge1$. For $L\ge1$, let $ \mu_{\lambda, L}:=\sum_{u\in\H
(L)}\delta
_{\{V(u)-L\}}1_{\{u\in\gg(L,\lambda)\}}$. Then Proposition \ref
{pyaglom} implies that $\mu_{\lambda,L}$ under $\mathbf{P}(\cdot|
H(L)>0)$ converges when $L\to\infty$ to $\hat\mu_{\B,\infty}$ defined
in Proposition \ref{pyaglom} with $\B(u):=\lambda^{-2}(\sum_{v\in
\I
(u)} \{1+ e^{\Delta V(v)}\})^2$. We will write $\hat\mu_{\lambda
,\infty
}:=\sum_{i=1}^{\hat\zeta_{\lambda}}\delta_{x_i}$ instead of $\hat
\mu
_{\B,\infty}$. Since the measures $\hat\mu_{\lambda,\infty}$ are
increasing in $\lambda$, we can assume that the labeling $(x_i)_{i}$
does not depend on $\lambda\ge1$. We write similarly $\mu_{\lambda
,\infty}:=\sum_{i=1}^{\zeta_{\lambda}}\delta_{x_i}$ for the
measure $\mu
_{\B,\infty}$ given by Proposition \ref{pyaglom}, and we know that
the Radon--Nikodym derivative of $\hat\mu_{\lambda,\infty}$ with
respect to $\mu_{\lambda,\infty}$ is ${\Re^{-1} \over\Q[\Re^{-1}]}$.
Notice that if $\hat\zeta_{\lambda}=0$, then $\hat\mu_{\lambda
,\infty
}$ is the measure zero.

\begin{lemma}\label{Ltail9} Assuming $\psi'(\varrho_*)=0$ with
$\varrho
_*=1$ and (\ref{hypcrit}), fix $\lambda\ge1$ and let $\hat\mu
_{\lambda
,\infty}$ and $\mu_{\lambda,\infty}$ be as above. Under $\Q$, let
$(\partial
W^{(i)}_\infty, i\ge1)$ be a sequence of i.i.d. random variables,
independent of $\hat\mu_{\lambda,\infty}$ and of common law that of
$\partial
W_\infty$ under $\mathbf{P}$. For any $\lambda\ge1$, we have
%
\begin{equation}\label{last2} \lim_{t\to\infty}t\Q\Biggl(\sum
_{i=1}^{\hat\zeta_{\lambda}} e^{
x_i}
\partial W_\infty^{(i)}
> t\Biggr) = {\Q[\Re^{-1}\sum_{i=1}^{\zeta_{\lambda}}
e^{x_i}]\over\Q
[\Re^{-1}] }.
\end{equation}
Moreover, for any $c >0$,
%
\begin{equation} \label{last9}
\lim_{\lambda\to\infty} \lambda^2 \Q\Biggl(\sum_{i=1}^{\hat
\zeta
_{\lambda}} e^{ x_i}
\partial W_\infty^{(i)}
> c \lambda^2 \Biggr) = {1\over c \Q[\Re^{-1}] }.
\end{equation}
\end{lemma}

\begin{pf}
For any $i\ge1$, by
Theorem 1.2 in
\cite{Bura09},
%
\begin{equation} \label{eqLiu-result-critiq}
\Q\bigl( \partial W^{(i)}_\infty> t \bigr)= \mathbf{P}( \partial
W_\infty>
t ) \sim
{1\over t},\qquad t \to\infty.
\end{equation}

In order to prove (\ref{last2}), we shall apply Lemma \ref
{Lconvolution}(i) and it is enough to show that there exists some
$\delta>0$ such that $\Q[ \sum_{i=1}^{\hat\zeta_{\lambda}}
(1+e^{x_i})^{1+\delta}] < \infty$.
Remark that $\hat\mu_{\lambda,\infty}$ has the support contained in
$\r
_+$, hence for $\delta>0$,
$\Q[ \sum_{i=1}^{\hat\zeta_{\lambda}} (1+e^{x_i})^{1+\delta
}
] \le2^{1+\delta} \Q[ \sum_{i=1}^{\hat\zeta_{\lambda}}
e^{(1+\delta)x_i} ]$. We are going to prove a stronger statement:
for $\widehat\mu_\infty$ the point process defined in Theorem \ref
{thmyaglom-critiq}(iii), we have
%
\begin{equation}\label{last6}
\Q\biggl[ \int\widehat\mu_\infty(dx) e^{(1+\delta) x} \biggr] <
\infty.
\end{equation}

The statement (\ref{last6}) implies the corresponding integrability for
$\hat\mu_{\lambda, \infty}$ since $\hat\mu_{\lambda, \infty}$ is
stochastically dominated by $\widehat\mu_\infty$. To prove (\ref
{last6}), we consider $\chi(L):= \sum_{ u \in\H(L)} e^{(1+\delta)
(V(u) -L)}$ and prove first that,
under $\mathbf{P}( \cdot| \H(L)\neq \varnothing)$, $\chi(L)$
converges in law to $\int\widehat\mu_\infty(dx) e^{(1+\delta) x} $.
In order to apply the convergence in law of Theorem \ref
{thmyaglom-critiq}(iii), we need some tightness result. We claim that
%
\begin{equation}
\label{eqcompact1}\quad
\sup_{L\ge1} \mathbf{P}_x\bigl(\exists i \in\bigl[\bigl\vert1, H(L)\bigr\vert\bigr] \dvtx
V\bigl(u^{(i)}\bigr)-L > K | H(L) > 0 \bigr) = o_K(1),
\end{equation}
where we order
the set of particles in $\H(L)$ (eventually empty) in an
arbitrary way: $\H(L)= \{u^{(i)}, 1\le i \le H(L)\}$. Markov's
inequality yields that the probability term in (\ref{eqcompact1}) is
smaller than
%
\begin{eqnarray}
\label{eqcompact12}
&&e^{-K} e^{-L} \e_x\biggl[\sum_{u\in\H(L)} e^{V(u)}\biggr] \mathbf{P}_x\bigl( H(L) > 0\bigr)^{-1}
\nonumber
\\[-8pt]
\\[-8pt]
\nonumber
&&\qquad\le c e^{-K} L \e_x\biggl[\sum_{u\in\H(L)} e^{V(u)}\biggr],
\end{eqnarray}
where the inequality is a consequence of Theorem \ref
{thmyaglom-critiq}(i). To prove the claimed tightness result it is
sufficient to show that
there exists some constant $c>0$ such that for any
$x\ge0$ and $L\ge\max(1,x)$ we have
%
\begin{equation}
\label{momentshl0} \e_x\biggl[\sum_{u\in\H(L)} e^{V(u)}\biggr] \le c
(1+x){e^{x}\over L}.
\end{equation}
To see it, a change of measure from $\mathbf{P}_x$ to $\Q_x$ by
Proposition \ref
{pchange-proba2} is applied to $\cc_{L}$ and $h(u):=e^{V(u)}$, and we
find that
\[
\e_x\biggl[\sum_{u\in\H(L)} e^{V(u)}\biggr]
=
e^{x} \Q_x\bigl(\tau_{L}^+ <\tau_{0}^-\bigr).
\]
Then (\ref{pal2}) implies (\ref{momentshl0}).
Assembling (\ref{eqcompact12}) and (\ref{momentshl0}) yields (\ref
{eqcompact1}) and allows us to apply Theorem \ref
{thmyaglom-critiq}(iii) to obtain the convergence in distribution,
under $\mathbf{P}(
\cdot|
\H(L)\neq \varnothing)$, of $\chi(L)$ toward $\int\widehat\mu
_\infty(dx) e^{(1+\delta) x}$.

Then (\ref{last6}) will hold once we have checked that $\e(\chi(L)
| \H(L) \neq \varnothing)$ is bounded on $L$. By Theorem \ref
{thmyaglom-critiq}(i) with $\varrho_*=1$, it is enough to show that
%
\begin{equation}\label{last8} \e\bigl[\chi(L) \bigr] \le c { e^{- L}
\over L}.
\end{equation}
But by the change of measure,
\[
\e\bigl[\chi(L) \bigr]= e^{-L} \Q\bigl[ e^{ \delta(S_{\tau_L^+}
-L)},
\tau_L^+ < \tau_0^-\bigr].
\]

The above expectation $\Q[\cdot]$ is less than ${c \over L}$ by
applying (\ref{pal}) to the random walk $( \delta(L-S_j)) _{j\ge0}$
(the integrability is guaranteed if $\delta$ is sufficiently small).
This proves (\ref{last8}) and a fortiori (\ref{last2}).

Recall that by (\ref{last6}) and Lemma \ref{Lconvolution}(i), if we
write $\widehat\mu_\infty= \sum_{i=1}^{\widehat\zeta} \delta_{\{
x_i\}
}$, then
\[
\Q\Biggl(\sum_{i=1}^{\hat\zeta} e^{
x_i} \partial W_\infty^{(i)}
> t\Biggr) \sim{\Q[\Re^{-1}\sum_{i=1}^{\zeta} e^{x_i}]\over\Q
[\Re
^{-1}] } {1 \over t} = {1\over\Q[\Re^{-1}] } {1 \over t},\qquad t
\to\infty
\]
since $\Re=\sum_{i=1}^{\zeta} e^{x_i}$ by definition; see (\ref
{eqref-R}). We have already observed that $\widehat\mu_{\lambda,
\infty
}$ is stochastically nondecreasing in $\lambda$ and is dominated by
$\widehat\mu_\infty$ ($\widehat\mu_\infty$ corresponds to
$\widehat\mu
_{\lambda, \infty}$ with $\lambda=\infty$). Then $\limsup_{\lambda
\to
\infty} \lambda^2 \Q(\sum_{i=1}^{\hat\zeta_{\lambda}} e^{ x_i}
\partial W_\infty^{(i)} > c \lambda^2 ) \le\limsup
_{\lambda
\to\infty} \lambda^2 \Q(\sum_{j=1}^{\hat\zeta} e^{
x_i} \partial W_\infty^{(i)} > c \lambda^2 )$ which is
${1\over c \Q[\Re^{-1}] }, $ yielding the upper bound in (\ref{last9}).

For the lower bound, let $\lambda_0>1$ and by the monotonicity in
$\hat
\mu_\lambda$,
\begin{eqnarray*}
\liminf_{\lambda\to\infty} \lambda^2 \Q\Biggl(\sum_{i=1}^{\hat
\zeta
_{\lambda}} e^{ x_i}
\partial W_\infty^{(i)} > c \lambda^2 \Biggr) &\ge& \liminf
_{\lambda\to\infty} \lambda^2 \Q\Biggl(\sum_{i=1}^{\hat\zeta
_{\lambda
_0}} e^{ x_i}
\partial W_\infty^{(i)} > c \lambda^2 \Biggr) \\
&=& {\Q[\Re^{-1}\sum_{i=1}^{\zeta_{\lambda_0}} e^{x_i}]\over c \Q
[\Re
^{-1}] },
\end{eqnarray*}
by applying (\ref{last2}) to $\widehat\mu_{\lambda_0, \infty} $.
Letting $\lambda_0\to\infty$ and noting that $\sum_{i=1}^{\zeta
_{\lambda
_0}} e^{x_i}= \int e^x \mu_{\lambda_0, \infty} (dx) \to\Re$, this
gives the lower bound of (\ref{last9}).
\end{pf}

We now have all the ingredients to prove Theorem \ref
{tleaves} in the critical case.

\begin{pf*}{Proof of Theorem \ref{tleaves}(i), (critical case)}
\mbox{}

\textit{Lower bound of Theorem \ref{tleaves}}(i). Recall that we have
assumed $\varrho_*=1$ by linear transformation.
Fix a constant $A>0$. Consider $n
\to\infty$, and let $L_{n,A}:= \log n +\log\log n -A$. We recall from
(\ref{hl}) that $H(L_{n,A})= \#\H(L_{n,A})$ is the number of
particles that
hit level $L_{n,A}$ before touching $0$. Recall \eqref{defbad}. We call
$H_{g}(L_{n,A}):=\#\H_{g}(L_{n,A})$ the number of particles in $\H
(L_{n,A})$\vspace*{1pt} which are in $\gg(L_n,\lambda)$ with $\lambda:=e^{A/
2}$,
%
\begin{equation} \label{defhgl} \H_g(L_{n, A}):= \H(L_{n, A}) \cap
\gg
\bigl(L_n, e^{A/2}\bigr).
\end{equation}

Let us order
the set of particles in $\H_{g}(L_{n,A})$ (eventually empty) in an
arbitrary way, $\H_{g}(L_{n,A})= \{u^{(i)}, 1\le i \le H_{
g}(L_{n,A})\}
$. Denote by
$\#\cL^{(i)}[0]$ the number of descendants of the $i$th particle
$u^{(i)} $ which are absorbed at $0$. Then
%
\begin{eqnarray}\label{eq1critlow}
&&\mathbf{P}_x\bigl(\#\cL[0]>n \bigr)\nonumber\\[-1pt]
&&\qquad \ge \mathbf{P}_x\Biggl(\sum
_{i=1}^{H_{g}(L_{n,A})}\#\cL^{(i)}[0]>n \Biggr)
\\[-1pt]
&&\qquad=
\mathbf{P}_x\bigl( H(L_{n,A}) >0 \bigr) \mathbf{P}_x\Biggl(\sum_{i=1}^{H_{
g}(L_{n,A})}\#\cL^{(i)}[0]>n
\Big| H(L_{n,A}) > 0 \Biggr).\nonumber
\end{eqnarray}

By Theorem \ref{thmyaglom-critiq}(i), $\mathbf{P}_x( H(L_{n,A})
> 0 ) \sim \frac{ \Q[\Re^{-1}]}{C_R} R(x)
e^x
\frac{e^{-L_{n,A}}}{L_{n,A}}$ as $n \to\infty$. On the other hand,
conditioned on $\H_{ g}(L_{n,A})$ and on $\{V(u^{(i)}),1\le i \le H_{
g}(L_{n,A})\}$, $(\#\cL^{(i)}[0])_{1\le i \le H_{ g}(L_{n,A})}$ are
independent, and each $\#\cL^{(i)}[0]$ is distributed as $\#\cL[0]$
under $\mathbf{P}_{V(u^{(i)})}$.

By Lemma \ref{Lconjecture2}, if we denote by $B^{(i)}:= \#\cL
^{(i)}[0] e^{-V(u^{(i)})} V(u^{(i)})$, then
conditioned on $\H_{ g}(L_{n,A})$ and on $\{V(u^{(i)}),1\le i \le H_{
g}(L_{n,A})\}$,
for each $i$, $B^{(i)}$ converges in law to $c^{*}
\partial W^{(i)}_\infty$ as $n \to\infty$, where $\partial
W^{(i)}_\infty, i\ge1$, is a sequence of i.i.d. random variables of
common law that of $(\partial
W_\infty, \mathbf{P})$, and independent of $\mu_{L_{n,A}}$. We may
assume by
Skorohod's representation theorem that for each $i$, $B^{(i)}$
converges almost surely to $c^*\partial W_\infty^{(i)}$.

Let $\varepsilon\in(0,1)$. First, we want to show that we can restrict
to the event $E(L_{n,A}):= \{ B^{(i)}>(1-\varepsilon) c^*
\partial
W^{(i)}_\infty; \forall i \dvtx1 \le i \le H_{ g}(L_{n,A})\}$. We have
\begin{eqnarray*}
&&\mathbf{P}_x\bigl(E(L_{n,A})^c | H(L_{n,A}) > 0\bigr)\\
&&\qquad\le
\e_x\bigl[H_{g}(L_{n,A}) | H(L_{n,A}) > 0\bigr]\sup_{z\ge L_{n,A}}
\mathbf{P}
_z\bigl( ze^{-z}\#\cL[0]<(1-\varepsilon)c^{*}\partial W_\infty\bigr)\\
&&\qquad=:\e_x\bigl[H_{ g}(L_{n,A}) | H(L_{n,A})> 0\bigr] \eta_{L_{n,A}}.
\end{eqnarray*}

The term $\eta_{L_{n,A}}$ goes to zero as $n\to\infty$ by
Lemma \ref{Lconjecture2}. By (\ref{momentshl0}) and Theorem~\ref
{thmyaglom-critiq}(i), we have
\begin{eqnarray*}
&&\e_x\bigl[H_{ g}(L_{n,A})| H(L_{n,A})
> 0\bigr]\\
&&\qquad\le e^{-L_{n,A}} \e_x\biggl[\sum_{u \in\H(L_{n,A})} e^{V(u)} \Big|
H(L_{n,A})> 0\biggr] \le c
\end{eqnarray*}
for some positive constant $c=c(x)$ which
depends on $x$. Hence, $\mathbf{P}_x(E(L_{n,A})^c
|\break  H(L_{n,A}) > 0)
=o_{L_{n,A}}(1)$, where $o_{L_{n,A}}(1)\to0$ as ${L_{n,A}}\to\infty
$.
We have
%
\begin{eqnarray}\label{eq2critlow}
&& \mathbf{P}_x\Biggl(\sum_{i=1}^{H_{g}(L_{n,A})}\#\cL^{(i)}[0]>n
\Big| H(L_{n,A}) > 0 \Biggr)\nonumber \\
&&\qquad=  \mathbf{P}_x\Biggl(\sum_{i=1}^{H_{g}(L_{n,A})} {e^{V(u^{(i)})}\over
V(u^{(i)})} B^{(i)} >
n
\Big| H(L_{n,A}) > 0 \Biggr) \\
&&\qquad\ge
\mathbf{P}_x\Biggl(\sum_{i=1}^{H_{g}(L_{n,A})} {e^{V(u^{(i)})}\over V(u^{(i)})}
B^{(i)} >
n, E(L_{n,A})
\Big| H(L_{n,A}) > 0 \Biggr).\nonumber
\end{eqnarray}

Observe that
%
\begin{eqnarray}\label{eq3critlow}\quad
&& \mathbf{P}_x\Biggl(\sum_{i=1}^{H_{ g}(L_{n,A})} {e^{V(u^{(i)})}\over
V(u^{(i)})} B^{(i)} >
n, E(L_{n,A})
\Big| H(L_{n,A}) > 0 \Biggr)\nonumber\\
&&\qquad\ge
\mathbf{P}_x\Biggl(\sum_{i=1}^{H_{g}(L_{n,A})} {e^{V(u^{(i)})}\over V(u^{(i)})}
\partial W^{(i)}_\infty> {n \over c^{*} (1-\varepsilon)}
, E(L_{n,A})
\Big| H(L_{n,A}) > 0 \Biggr)\\
&&\qquad\ge
\mathbf{P}_x\Biggl(\sum_{i=1}^{H_{g}(L_{n,A})} {e^{V(u^{(i)})}\over V(u^{(i)})}
\partial W^{(i)}_\infty> {n \over c^{*} (1-\varepsilon)} \Big|
H(L_{n,A}) > 0 \Biggr) + o_{L_{n,A}}(1).\nonumber
\end{eqnarray}

In order to apply the convergence in law of Proposition \ref
{pyaglom}, we need some tightness results. Recalling (\ref
{eqcompact1}), it is sufficient to show that
\[
\sup_{L\ge1} \mathbf{P}_x\bigl(\exists i \in\bigl[\bigl\vert1, H(L)\bigr\vert\bigr] \dvtx
\partial W^{(i)}_\infty> K | H(L) > 0 \bigr) = o_K(1).
\]
Since the $\partial
W^{(i)}_\infty$'s are i.i.d. copies of $\partial
W_\infty$ and independent of $\mu_{L_{n,A}}$, Markov's inequality
yields that the probability term in the previous equation is smaller than
\[
K^{-1/2} \e_x \bigl[H(L) | H(L) > 0\bigr] \e[\sqrt{\partial W_\infty}] =
O\bigl( K^{-1/2}\bigr),
\]
by using (\ref{momentshl0}), Theorem \ref{thmyaglom-critiq}(i) and
(\ref{eqLiu-result-critiq}).
This yields the claimed tightness and allows us to apply Proposition
\ref{pyaglom} to get
%
\begin{eqnarray}\label{eq4critlow}
&&\lim_{n\to\infty} \mathbf{P}_x\Biggl(\sum_{i=1}^{H_{g}(L_{n,A})}
{e^{V(u^{(i)})}\over V(u^{(i)})} \partial W^{(i)}_\infty>{n \over
c^{*} (1-\varepsilon)} \Big| H(L_{n,A}) > 0 \Biggr)
\nonumber
\\[-8pt]
\\[-8pt]
\nonumber
&&\qquad= \Q\Biggl(\sum_{i=1}^{\hat\zeta_{\lambda}} e^{ x_i}
\partial W_\infty^{(i)} >
{e^A \over c^{*} (1-\varepsilon)} \Biggr),
\end{eqnarray}
where $\hat\mu_{\lambda,\infty}:=\sum_{i=1}^{\hat\zeta
_{\lambda}} \delta_{x_i} $ is
the point process defined before Lemma \ref{Ltail9}, and we recall
that $\lambda:=e^{A/2}$. By (\ref{eq1critlow})--(\ref{eq4critlow}) and the
definition of $L_{n,A}$, we deduce that for
any $A>0$,
\begin{eqnarray*}
&&\liminf_{n\to\infty}n (\log n)^2 \mathbf{P}_x\bigl(\cL[0]>n
\bigr) \\
&&\qquad\ge\frac
{ \Q[\Re^{-1} ]}{C_R} R(x) e^x e^{A} \Q
\Biggl(\sum
_{i=1}^{\hat\zeta_{\lambda}} e^{ x_i} \partial W_\infty^{(i)} >
{\lambda^2 \over c^{*} (1-\varepsilon)} \Biggr).
\end{eqnarray*}

We let $\varepsilon\to0$ to get
\[
\liminf_{n\to\infty}n (\log n)^2 \mathbf{P}_x\bigl(\#\cL[0] > n
\bigr) \ge R(x)e^xC(A),
\]
with $C(A):=\frac{ \Q[\Re^{-1} ]}{C_R} e^{A}
\Q(\sum_{i=1}^{\hat\zeta_{\lambda}} e^{ x_i}
c^{*} \partial W_\infty^{(i)} > \lambda^2)$.

By Lemma \ref{Ltail9}, we have $C(A) \to{ c^*\over C_R}$ as $A \to
\infty$, which leads to
\[
\liminf_{n\to\infty}n (\log n)^2 \mathbf{P}_x\bigl(\#\cL[0] > n
\bigr) \ge
R(x)e^x\frac{ c^* }{C_R}.
\]

\textit{Upper bound of Theorem \ref{tleaves}}(i).
We notice that we showed in fact that, for any $A>0$,
\[
\liminf_{n\to\infty}n (\log n)^2 \mathbf{P}_x\Biggl(\sum
_{i=1}^{H_{g}(L_{n,A})}\#\cL^{(i)} [0] >
n \Biggr) \ge R(x)e^xC(A).
\]

Repeating the same argument with this time $E'(L_{n,A}):= \{
B^{(i)}<(1+\varepsilon)\partial W^{(i)}_\infty;  \forall i \dvtx1
\le i \le H_{g}(L_{n,A}) \}$ yields that $C(A)$ is also a limsup. Therefore,
%
\begin{equation}\label{57}
\lim_{n\to\infty}n (\log n)^2 \mathbf{P}_x\Biggl(\sum
_{i=1}^{H_{g}(L_{n,A})}\#
\cL^{(i)}[0]>n \Biggr) = R(x)e^x C(A),
\end{equation}
with $C(A) \to{ c^*\over C_R}$ as $A \to\infty$.

Then, let $\eta>0$ and $\varepsilon>0$. We take again $L_{n,A}:= \log n
+ \log\log n -A$ and $\lambda:=e^{A/2}$. Markov's inequality with
(\ref{eqbadL[0]}) implies that if $A$ is taken large enough,
\[
\limsup_{n\to\infty}n(\log n)^2\mathbf{P}_x\bigl(Z_b[0,L_{n,A}]>\eta
n\bigr)
\le\varepsilon.
\]

By Theorem \ref{thmyaglom-critiq}(i), we can choose $B>0$
large enough such that
%
\begin{equation}
\limsup_{n} n(\log n)^2 \mathbf{P}_x\bigl(H(L_{n}+B)>0\bigr) \le\varepsilon.
\label
{eqbadHLB}
\end{equation}

On the other hand, by (\ref{eqbadHL}) and Markov's
inequality, we obtain that for $A$ large enough,
%
\begin{eqnarray}\label{eqbadHL0}
&& \limsup_{n} n(\log n)^2\mathbf{P}_x\biggl(\sum_{u \in\H
(L_{n,A})}\!\! 1_{\{u
\in\bb(L_{n,A},\lambda)\}}\#\cL^{(u)}[0] > \eta n, H(L_n +B) =0
\biggr)
\nonumber\hspace*{-25pt}\\
&&\qquad\le
\limsup_{n} n(\log n)^2{1\over\eta n} \e_x\biggl[\sum_{u\in\H
(L_{n,A})} 1_{\{u \in\bb(L_{n,A},\lambda)\}}Z^{(u)}[0,L_n +B]
\biggr]\hspace*{-25pt}\\
&&\qquad\le
\varepsilon,\nonumber\hspace*{-25pt}
\end{eqnarray}
where the notation $Z^{(u)}[,]$ was introduced in Lemma \ref
{lbad}. Finally, it yields that
%
\begin{equation}\qquad
\limsup_{n} n(\log n)^2\mathbf{P}_x\biggl(\sum_{u \in\H(L_{n,A})}
1_{\{u \in
\bb(L_{n,A},\lambda)\}}\#\cL^{(u)}[0] > \eta n \biggr)
\le
2\varepsilon.
\end{equation}

We now show that the ``good particles'' which never touch
$L_{n,A}$ are negligible when $A$ is large. We recall that
$Z_{g}(0,L_{n,A})$ is the number of particles in $\gg(L_n,\lambda)$
that touch $0$
before $L_{n,A}$. By Lemma \ref{L2moment},
\[
\e_x\bigl[Z_{g}(0,L_{n,A})^2\bigr] \le c (1+x) e^x \lambda
{e^{L_{n,A}} \over L_{n,A}^3}.
\]

Therefore, by the choice of $L_{n,A}$ and $\lambda$ we have
that for any fixed $\eta>0$,
\[
\label{58} \limsup_{n\to\infty} n (\log n)^2
\mathbf{P}_x\bigl(Z_{g}[0,L_{n,A}]>\eta n\bigr) \le{ c (1+x) e^x
e^{-{A/2}} \over\eta^2},
\]
which is less than $\varepsilon$ if $A$ is large enough. By
the triangle inequality, for any $0<\eta<1/3 $ and any $\varepsilon>0$,
we deduce that if $A$ is large enough,
\[
\mathbf{P}_x\bigl(\#\cL[0] >n\bigr) \le\mathbf{P}_x\Biggl( \sum
_{i=1}^{H_{g}(L_{n,A})}\#\cL
^{(i)}[0] > (1-3\eta) n\Biggr) + \frac{4 \varepsilon}{n (\log n)^2}.
\]

From this and (\ref{57}), by letting $A \to\infty$ and
$\eta
\to0$, we deduce the upper bound
\[
\limsup_{n\to\infty}n (\log n)^2 \mathbf{P}_x\bigl(\#\cL[0] > n
\bigr) \le
R(x)e^x\frac{ c^*}{C_R}.
\]

Thus we have
\[
\lim_{n\to\infty}n (\log n)^2 \mathbf{P}_x\bigl(\#\cL[0] > n \bigr)
= R(x)e^x c_{\mathrm{crit}}',
\]
with $c_{\mathrm{crit}}' = \frac{ c^* }{C_R}$. Finally, we recall that $C_R$ is
the limit of $R(x)/x$ as $x\to\infty$, $R(x)$ being the renewal
function for the descending ladder heights. The renewal theorem implies
that $C_R = \Q[-S_{\tau_0^-}]^{-1}$. Hence, from the value of $c^*$ in
(\ref{eqconstant1}), we end up with $c_{\mathrm{crit}}' = \Q[e^{-S_{\tau_0^-}}
- 1]$ indeed.
\end{pf*}

\section{\texorpdfstring{Proof of Theorem \protect\ref{tleaves}: The subcritical case}
{Proof of Theorem 2: The subcritical case}}\label{prsubcritical}

We treat here the subcritical case $\psi'(\varrho_*)<0$. Define a new
probability measure $\Q^{(\varrho_-)}$ by (\ref{Q}) with $h(u)=\ee
^{\varrho_- V(u)}$ for all $u \in\cT$. Then for any $x\in\r$,
\[
{d \Q_x^{(\varrho_-)}\over d \mathbf{P}_x}\bigg|_{\F_n} = e^{-\varrho
_- x} \sum
_{\vert u \vert=n} e^{ \varrho_- V(u)},\qquad n\ge0.
\]
We recall that $\Q$ satisfies (\ref{Q1}) with $\varrho=\varrho_+$.

Applying Proposition \ref{pchange-proba}, we see that the trajectory
of the spine $(S_n)$ is a random walk that drifts to $+\infty$ under
$\Q$, and drifts to $-\infty$ under $\Q^{(\varrho_-)}$, in fact,
$\Q
[S_1]= \psi'(\varrho_+)>0$ and $\Q^{(\varrho_-)}[S_1] = \psi
'(\varrho
_-) <0$. In particular [see (\ref{subtau-}) and~(\ref{subtau+}),
changing $S_1$ in $-S_1$ for $\Q^{(\varrho_-)}$], we deduce the
existence of $C_R^{(\varrho_-)}>0 $ such that
%
\begin{eqnarray}\label{est5sub}
\Q^{(\varrho_-)}\bigl(\tau_L^+<\tau_0^-\bigr)& \sim&{1\over C_R^{(\varrho_-)}}
e^{(\varrho_- - \varrho_+)L},
\nonumber
\\[-8pt]
\\[-8pt]
\nonumber
\Q\bigl(\tau_L^+<\tau_0^-\bigr) &\sim&{1\over C_R },\qquad L\to\infty,
\label{est6sub}
\end{eqnarray}
(the second equivalence follows from Lemma \ref{lemRW}).
The strategy of the proof of Theorem \ref{tleaves}(ii) is in the same
spirit as in the critical case (i). Recall (\ref{l[a]}) that $\cL[0]$
denotes the set of leaves of the killed branching random walk. We give
first an estimate on the moments of $\#\cL[0]$.

\begin{lemma} \label{lmomentsubL}
For any integer $k< {\varrho_+ \over
\varrho_-}$, there exists some constant $c_{k}>0$ such that for any
$x\ge0$,
\[
\e_x\bigl[\bigl(\#\cL[0]\bigr)^{k}\bigr] \le c_{k} e^{k \varrho_- x}.
\]
\end{lemma}

\begin{pf*}{Proof of Lemma \ref{lmomentsubL}} We give a proof
by induction on $k$. Changing
measure from $\mathbf{P}_x$ to $\Q_x^{(\varrho_-)}$ with Proposition
\ref
{pchange-proba2} (with $\cL[0]$ and $h(u)=e^{\varrho_- V(u)}$ for
$u\in\cT$) yields
the identity
%
\begin{equation}\label{eqmomentsub1}
\e_x\bigl[\bigl(\#\cL[0]\bigr)^{k }\bigr]
=
e^{\varrho_- x}\Q_x^{(\varrho_-)}\bigl[e^{-\varrho_- S_{\tau
_0^-}}
\bigl(\#\cL[0]\bigr)^{k-1}\bigr].
\end{equation}

By (\ref{est1sub}), the case $k=1$ holds. Suppose that it is
true for $k-1 \ge1$, and that $2\le k < {\varrho_+\over\varrho_-}$.
We decompose $\#\cL[0]$ along the spine
\[
\#\cL[0] = 1 + \sum_{\ell=1}^{\tau_0^-} \sum_{u\in\I_{\ell}} \#
\cL^{(u)}[0],
\]
where $\#\cL^{(u)}[0]$ is the number of particles descended
from $u$ absorbed at $0$. We mention that if $V(u)<0$, then $\#\cL^{(u)}[0]=1$.
Conditionally on $\mathcal{G}_\infty$, $(\#\cL^{(u)}[0])_{u\in\I
_{j+1}}$, $0\le j < \tau_0^-$, are independent and each $\#\cL
^{(u)}[0]$ is distributed as $(\#\cL[0],\mathbf{P}_{V(u)})$. By the
triangle inequality,
\begin{eqnarray*}
\bigl(\Q_x^{(\varrho_-)}\bigl[\bigl(\#\cL[0]\bigr)^{k-1} |
\mathcal
{G}_\infty\bigr]\bigr)^{1/(k-1)}
\le1+
\sum_{\ell=1}^{\tau_0^-}\sum_{u\in\I_\ell}\bigl( \Q
_x^{(\varrho
_-)}\bigl[\bigl( \#\cL^{(u)}[0] \bigr)^{k-1} | \mathcal
{G}_\infty
\bigr]\bigr)^{1/(k-1)}.
\end{eqnarray*}

For each $\ell$ and $u\in\I_\ell$, we have from our induction
assumption
%
\begin{eqnarray}
\label{induc}\quad  \Q_x^{(\varrho_-)}\bigl[\bigl( \#\cL^{(u)}[0]
\bigr)^{k-1} | \mathcal{G}_\infty\bigr]
&\le&
1_{\{V(u)<0\}}+1_{\{V(u)\ge0\}}c_{k-1} e^{\varrho_- (k-1) V(u)}
\nonumber
\\[-8pt]
\\[-8pt]
\nonumber
&\le& c\bigl(1+e^{\varrho_- V(u)}\bigr)^{k-1}.
\end{eqnarray}

Therefore we get
\[
\Q_x^{(\varrho_-)}\bigl[\bigl(\#\cL[0] \bigr)^{k-1} |
\mathcal
{G}_\infty\bigr]^{1/(k-1)}
\le1+
c' \sum_{\ell=1}^{\tau_0^-} \sum_{u\in\I_\ell} 1 + e^{\varrho_-
V(u)}.
\]

In view of (\ref{eqmomentsub1}), we deduce that
\begin{eqnarray*}
\e_x\bigl[\bigl(\#\cL[0]\bigr)^{k }\bigr]
&\le& c e^{\varrho_- x} +
c e^{\varrho_- x} \Q_x^{(\varrho_-)}\Biggl[e^{-\varrho_-
S_{\tau
_0^-}} \Biggl(\sum_{\ell= 1}^{\tau_0^-} \sum_{u\in\I_\ell}\bigl \{
1+e^{\varrho
_- V(u)}\bigr\}\Biggr)^{k-1}\Biggr] \\
&\le& c e^{\varrho_- x} +
c e^{\varrho_- x} \Q_x^{(\varrho_-)}\Biggl[e^{-\varrho_-
S_{\tau
_0^-}} \Biggl(\sum_{\ell= 1}^{\tau_0^-} e^{\varrho_- S_{\ell
-1}}a_{\ell
}\Biggr)^{k-1}\Biggr],
\end{eqnarray*}
where for any $\ell\ge1$, $a_\ell:= \sum_{u\in\I_{\ell
}} \{
1+e^{\varrho_- \Delta V(u)}\}$. Plainly Corollary \ref{cmany1} also
holds with $\varrho=\varrho_-$, which implies that under $\Q
_x^{(\varrho
_-)}$, the random variables $(S_\ell- S_{\ell-1}, a_{\ell})_{\ell
\ge
1}$ are i.i.d. (whose law does not depend on $x$). Moreover,
\[
\Q^{(\varrho_-)} \bigl[ \bigl(1+ 1_{\{S_1<0\}} e^{ - \varrho_- S_1}\bigr)
a_1^{k-1}\bigr] \le\e\Biggl[ \sum_{|u|=1}\bigl (1+e^{\varrho_- V(u)}\bigr) \Biggr]^k
< \infty,
\]
by (\ref{hypsub}). Applying (\ref{est8sub}) with $b=\varrho
_-$, $p=k-1$, $\gamma=\varrho_+-\varrho_-$ (recalling that $\varrho
_+/\varrho_- > k \ge2$), we get $\Q_x^{(\varrho_-)}
[e^{-\varrho_-
S_{\tau_0^-}} (\sum_{\ell= 1}^{\tau_0^-} e^{\varrho_-
S_{\ell
-1}}a_{\ell})^{k-1}] \le c e^{ (k-1) \varrho_- x}$,
proving the lemma.
\end{pf*}

We introduce the analog of \textit{good} and \textit{bad} particles
in the
subcritical case, and we feel free to use the same notation. For
$\lambda>1, L>1$, we say now that
%
\begin{eqnarray} \label{defbad2} u \in\bb(L, \lambda) \mbox{ if there
exists some } 1\le j \le|u|\dvtx
\nonumber
\\[-8pt]
\\[-8pt]
\eqntext{\displaystyle \sum_{v \dvtx{\buildrel
\leftarrow\over v} = u_{j-1} }\bigl (1+ e^{ \varrho_- \Delta V(v) }\bigr) >
\lambda e^{ \varrho_- ( L- V(u_{j-1}))},}
\end{eqnarray}
and $u \in\gg(L, \lambda)$ otherwise, and we define again
%
\begin{eqnarray} \label{defzgb2} Z_g[0, L]&:=& \sum_{u \in\gg(L,
\lambda
)} 1_{\{ \tau_0^-(u) = |u| < \tau_L^+(u)\}},
\nonumber
\\[-8pt]
\\[-8pt]
\nonumber
 Z_b[0, L]&:=& \sum
_{u\in\bb(L, \lambda) } 1_{\{ \tau_0^-(u) = |u| < \tau_L^+(u)\}}.
\end{eqnarray}
Recall the notation $\delta^*$ in (\ref{hypsub}).

\begin{lemma} \label{lmomentsubZ} Let $k^*:= \lfloor{ \varrho_+
\over
\varrho_- }\rfloor+1$ be the smallest integer such that
$k^*> {\varrho_+\over\varrho_-}$. Let $0< \delta_2 < \min( {\delta
^*\over2},k^*- { \varrho_+\over\varrho_ - } )$.
\begin{longlist}[(iii)]
\item[(i)] There exists some constant $c >0$ such that for any $L>x\ge0$,
\[
\e_x\bigl[Z_g[0,L]^{k^*}\bigr] \le c \lambda^{ k^*- { \varrho
_+/\varrho_ - } - \delta_2} e^{\varrho_+ x} e^{(\varrho_-
k^* -
\varrho_+)L}.
\]

\item[(ii)] For $q:= {\varrho_+\over\varrho_-} + \delta_2$, there exists some
constant $c':=c'(\lambda, q)>0$ such that for any $L>x\ge0$,
\[
\e_x\biggl[ \sum_{u \in
\H(L) \cap\gg(L,\lambda) } e^{\varrho_- V(u)} \biggr]^{q} \le c'
e^{\varrho_+ x} e^{ (q \varrho_- - \varrho_+)L}.
\]

\item[(iii)] If we assume (\ref{hyp-extra}), then
\[
\e_x\biggl[ \sum_{u \in
\H(L) } e^{\varrho_- V(u)} \biggr]^{k^*} \le c e^{\varrho_+ x}
e^{ (k^* \varrho_- - \varrho_+)L},\qquad 0\le x < L.
\]
\end{longlist}
\end{lemma}

\begin{pf*}{Proof of Lemma \ref{lmomentsubZ}}
(i) Let $k$ be an integer. By
changing of measure from $\mathbf{P}_x$ to $\Q_x^{(\varrho_-)}$, we obtain
%
\begin{eqnarray}\label{changsub1} &&\e_x\bigl[\bigl(Z_g[0,L]\bigr)^{k}\bigr]
\nonumber
\\[-8pt]
\\[-8pt]
\nonumber
&&\qquad
= e^{\varrho_- x} \Q_x^{(\varrho_-)}\bigl[e^{-\varrho_- S_{\tau_0^-}}
1_{\{ w_{\tau_0^-} \in\gg(L,\lambda) \}} \bigl(Z_g[0,L]\bigr)^{k-1}, \tau
_0^-<\tau_L^+\bigr].
\end{eqnarray}

By decomposing the tree $\cT$ along the spine $(w_\ell)$, we get that
%
\begin{equation}\label{zgL1} Z_g[0, L] \le Z[0, L]= 1 + \sum_{\ell
=1}^{\tau_0^-} \sum_{u \in\I_\ell} Z^{(u)}[0, L],
\end{equation}
where $Z^{(u)}[0, L]:= \sum_{v \in\cT^{(u)}} 1_{\{\tau
_0^-(v) = |v| < \tau_{L }^+(v)\}} $ denotes the number of descendants
of $u$, touching
$0$ before $L $ ($\cT^{(u)}$ means as before the subtree rooted at~$u$).
By Proposition \ref{pchange-proba}, under $\Q_x$,
conditioned on $\G_\infty=\sigma\{\omega_j, S_j, \I_j, (V(u),\break u\in \I
_j), j\ge0\}$, the random variables $(Z^{(u)}[0, L ])_{u \in\I_\ell,
\ell\le\tau_0^-}$ are
independent and each $Z^{(u)}[0, L]$ is distributed as $(Z[0, L],
\mathbf{P}_{V(u)})$. Conditioning and using the triangle inequality,
we have
%
\begin{eqnarray}\label{eqcondsub2}
&&\bigl(\Q_x^{(\varrho_-)}\bigl[\bigl(Z_g[0,L]\bigr)^{k-1} | \mathcal
{G}_\infty
\bigr]\bigr)^{1/(k-1)}
\nonumber
\\[-8pt]
\\[-8pt]
\nonumber
&&\qquad
\le1+
\sum_{\ell=1}^{\tau_0^- }\sum_{u\in\I_\ell}\bigl(\Q
_x^{(\varrho
_-)}\bigl[\bigl( Z^{(u)}[0,L] \bigr)^{k-1} | \mathcal
{G}_\infty
\bigr]\bigr)^{1/(k-1)}.
\end{eqnarray}

Assume $k< (\varrho_+/\varrho_-) +1$. From Lemma \ref{lmomentsubL}, since
$Z^{(u)}[0,L] \le\#\cL^{(u)}[0]$ and $k-1<\varrho_+/\varrho_-$, we
know that
\[
\bigl(\Q_x^{(\varrho_-)}\bigl[\bigl( Z^{(u)}[0,L] \bigr)^{k-1}
|
\mathcal{G}_\infty\bigr]\bigr)^{1/(k-1)}\le
c e^{\varrho_- V(u)} + 1_{\{V(u) <0\}},
\]\eject\noindent
where the indicator comes from $Z^{(u)}[0,L] =1$ if $V(u)
<0$. It follows that
\fontsize{10pt}{\baselineskip}\selectfont
\makeatletter
\def\tagform@#1{\normalsize\maketag@@@{(\ignorespaces#1\unskip\@@italiccorr)}}
\makeatother
\begin{eqnarray}\label{zg1}&& \e_x\bigl[\bigl(Z_g[0,L]\bigr)^{k}\bigr] \nonumber\\
&&\qquad\le c e^{\varrho_- x} \Q
_x^{(\varrho_-)}\bigl[e^{-\varrho_- S_{\tau_0^-}} \bigr]
\nonumber
\\[-8pt]
\\[-8pt]
\nonumber
&&\qquad\quad{}+ c e^{\varrho_- x} \Q_x^{(\varrho_-)}\Biggl[e^{-\varrho_-
S_{\tau
_0^-}} 1_{\{ w_{\tau_0^-} \in\gg(L,\lambda),\tau_0^-<\tau_L^+ \}}
\Biggl( \sum_{\ell=1}^{\tau_0^- }\sum_{u\in\I_\ell} \bigl(1+e^{\varrho_- V(u)}\bigr)
\Biggr)^{k-1}\Biggr] \\
&&\qquad=: c e^{\varrho_- x} \Q_x^{(\varrho_-)}\bigl[e^{-\varrho_-
S_{\tau
_0^-}} \bigr] + c e^{\varrho_- x} A_{\scriptsize(\ref{zg1})},\nonumber
\end{eqnarray}
\normalsize
with some larger constant $c >0$ and the obvious definition
of $A_{\scriptsize(\ref{zg1})}$ for the remaining expectation under $\Q
_x^{(\varrho
_-)}$. By (\ref{est1sub}), see also Theorem 4 in \cite{L70} applied to
$-S$ at $\tau_x^+$, $\Q_x^{(\varrho_-)}[e^{-\varrho_- S_{\tau
_0^-}} ] \le c$. Therefore we have shown that for all $k<
(\varrho
_+/\varrho_-) +1$,
%
\begin{equation} \label{zg11} \e_x\bigl[\bigl(Z_g[0,L]\bigr)^{k}\bigr] \le c' e^{\varrho_-
x} + c e^{\varrho_- x} A_{\scriptsize(\ref{zg1})}.
\end{equation}

To estimate $A_{\scriptsize(\ref{zg1})}$, let us adopt the notation $a_{\ell}$:
for any $\ell\ge1$, $a_\ell:= \sum_{u\in\I_{\ell}} (1+
e^{\varrho_-
\Delta V(u)})$, hence $\sum_{\ell=1}^{\tau_0^- }\sum_{u\in\I_\ell}
(1+e^{\varrho_- V(u)} )\le\sum_{\ell=1}^{\tau_0^- } e^{\varrho_-
S_{\ell-1}} a_\ell$. On $\{w_{\tau_0^-} \in \gg(L, \lambda)\}$,
$a_\ell\le\lambda^s e^{ s \varrho_- (L-S_{\ell-1})} a_\ell
^{1-s}$ for any $0< s <1$. It follows that
%
\begin{eqnarray}\label{zg2}\quad  A_{\scriptsize(\ref{zg1})} &\le&\lambda^{s( k-1)}
e^{ s \varrho_- ( k-1) L}
\nonumber
\\[-8pt]
\\[-8pt]
\nonumber
&&{}\times\Q_x^{(\varrho_-)}\Biggl[e^{-\varrho_-
S_{\tau_0^-}} \Biggl( \sum_{\ell=1}^{\tau_0^- } e^{ \varrho_- (1 - s)
S_{\ell-1}} a_\ell^{1-s} \Biggr)^{k-1}, \tau_0^-<\tau_L^+\Biggr],
\end{eqnarray}
for any $0< s <1$ and $ k < (\varrho_+/\varrho_-) +1$.

If $\varrho_+/\varrho_-$ is not an integer, then $k^* < { \varrho
_+\over\varrho_-} +1$ and (\ref{zg2}) holds for $k= k^*$. Take
%
\begin{equation}\label{defs}
s= { k^*- { \varrho_+/\varrho_ - } - \delta_2 \over k^*-1}.
\end{equation}

Notice that
\begin{eqnarray*}  \Q_x^{(\varrho_-)}\bigl[ \bigl(1+ 1_{\{S_1<0\}} e^{-
\varrho_- S_1}\bigr) a_1^{ (1-s) (k^*-1)}\bigr] &\le&\e\biggl[ \sum_{
\vert u \vert=1}\bigl (1+ e^{ \varrho_- V(u)} \bigr)\biggr]^{ (1-s) (k^*-1) +1}
\\
&<& \infty, \\
\Q^{(\varrho_-)}\bigl[ e^{ (1-s) \varrho_- (k^*-1) S_1}
\bigr]&=& e^{ \psi( \varrho_- ( 1+ (1-s) ( k^*-1)))} < \infty,
\end{eqnarray*}
by (\ref{hypsub}). Under $\Q^{(\varrho_-)}$, $(S_\ell-
S_{\ell-1}, a_\ell^{1-s})_{\ell\ge1}$ are i.i.d. Applying (\ref
{est9sub}) (with $\alpha=0$) to the expectation term $\Q_x^{ (\varrho
_-)}[\cdot]$ in (\ref{zg2}) with $ \gamma= \varrho_+ - \varrho_-, b=
\varrho_- (1-s), \eta= \varrho_-, p=k^*-1$ and noticing that $pb >
\gamma$, we get that if we take $k=k^*$ in (\ref{zg1}), then
\begin{eqnarray*}
A_{\scriptsize(\ref{zg1})} &\le& c \lambda^{s( k^*-1)} e^{ s \varrho_- (
k^*-1) L} e^{ (\varrho_+- \varrho_-) (x -L) + (k^*-1) (\varrho_- - s
\varrho_- )L} \\
&=& c \lambda^{s(k^*-1)} e^{ (\varrho_+- \varrho_-) (x -L) + (k^*-1)
\varrho_- L}.
\end{eqnarray*}
This estimate with (\ref{zg11}) proves (i) in the case that
$\varrho_+/\varrho_-$ is not an integer.

It remains to treat the case when $\varrho_+/\varrho_-$ is an integer.
Then $k^*= { \varrho_+\over\varrho_-} +1$. Applying (\ref{zg1}) to $k=
k^*-1$ (which is less than ${ \varrho_+\over\varrho_-} +1$), we have
that
\begin{eqnarray*}
&&\e_x\bigl[\bigl(Z_g[0,L]\bigr)^{k^*-1}\bigr] \\
&&\qquad\le c' e^{\varrho_- x} + c e^{\varrho_- x}
\Q_x^{(\varrho_-)}\Biggl[e^{-\varrho_- S_{\tau_0^-}}\Biggl ( \sum
_{\ell
=1}^{\tau_0^- } e^{\varrho_- S_{\ell-1}} a_\ell
\Biggr)^{k^*-2},
\tau_0^-<\tau_L^+\Biggr],
\end{eqnarray*}
which by an application of (\ref{est9sub}) with $\alpha=0,
\gamma= \varrho_+- \varrho_-, b= \varrho_-, p=k^*-2=\gamma/b $ [it is
easy to check the integrability hypothesis in Lemma \ref{Laux1}(ii)],
yields that
\[
\e_x\bigl[\bigl(Z_g[0,L]\bigr)^{k^*-1}\bigr] \le c ( 1+ L- x) e^{ \varrho_+ x},\qquad
0\le x \le L.
\]

Moreover, $\e_x[(Z_g[0,L])^{k^*-1}]$ is $1$ if $x<0$ and $0$ if $x>L$.
Going back to (\ref{eqcondsub2}) and (\ref{changsub1}) with now
$k=k^*$, we obtain that
\[
\e_x\bigl[\bigl(Z_g[0,L]\bigr)^{k^* }\bigr] \le
c e^{ \varrho_- x} \Q_x^{(\varrho_-)} \bigl[1+ e^{-\varrho_-
S_{\tau_0^-}}
1_{\{ w_{\tau_0^-} \in\gg(L,\lambda) \}} A^{k^*-1}, \tau_0^-<\tau
_L^+\bigr]
\]
with
\[
A:=\sum_{\ell=1}^{\tau_0^- }\sum_{u\in\I_\ell} \bigl(\bigl( 1+L- V(u)\bigr)^{
\varrho_-/\varrho_+} e^{\varrho_- V(u)}1_{\{V(u)\in[0,L]\}} +
1_{\{
V(u)<0\}}\bigr).
\]
Observe that on $\{ \ell\le\tau_0^-<\tau_L^+\}$, $S_{\ell
-1} \in[0, L]$. For any $u \in\I_\ell$ such that $V(u) \in[0, L]$,
either $\Delta V(u) \ge0$ then $(1+L- V(u))^{ \varrho_-/ \varrho_+}
\le
(1+L- S_{\ell-1})^{ \varrho_-/ \varrho_+} $, or $\Delta V(u)< 0$, then
$ ( 1+L- V(u))^{ \varrho_-/ \varrho_+} e^{\varrho_- V(u)} \le(1+L-
S_{\ell-1})^{ \varrho_-/ \varrho_+} e^{ \varrho_- V(u)} + \vert
\Delta
V(u) \vert^{ \varrho_-/ \varrho_+} e^{ \varrho_- V(u)} \le(1+L-
S_{\ell-1})^{ \varrho_-/ \varrho_+} e^{ \varrho_- V(u)} + c(\varrho
)
e^{ \varrho_- S_{\ell-1}} $, with $c(\varrho):= \sup_{y \le0}
\vert y
\vert^{ \varrho_-/ \varrho_+} e^y< \infty$. It follows that there
exists some $c>0$ such that
\begin{eqnarray*}
&& \sum_{u\in\I_\ell} \bigl( \bigl( 1+L- V(u)\bigr)^{ \varrho_-/
\varrho_+}
e^{\varrho_- V(u)} 1_{\{V(u)\in[0,L]\}} + 1_{\{V(u)<0\}} \bigr)\\
&&\qquad\le c (1+L- S_{\ell-1})^{ \varrho_-/\varrho_+} e^{\rho_-
S_{\ell-1}}\sum_{u\in\I_\ell} \bigl(1+e^{\varrho_- \Delta V(u)} \bigr),
\end{eqnarray*}
which in turn is bounded by $ c (1+L- S_{\ell-1})^{ \varrho
_-/\varrho_+}e^{\rho_- S_{\ell-1}} \lambda^s e^{s \varrho_-
(L-S_{\ell-1})} a_\ell^{1-s}$ since $w_{\tau_0^-} \in\gg(L,
\lambda)$,
where $0< s <1$ is as in (\ref{defs}). It follows that
\begin{eqnarray*} &&\e_x\bigl[\bigl(Z_g[0,L]\bigr)^{k^* }\bigr] \\
&&\qquad\le
c' \lambda^{s( k^*-1)} e^{ s\varrho_- (k^*-1) L} e^{ \varrho_- x}
\\
&&\qquad\quad{}\times \Q_x^{(\varrho_-)} \Biggl[ e^{-\varrho_- S_{\tau_0^-}}
\Biggl( \sum_{\ell=1}^{\tau_0^- } ( 1+L- S_{\ell-1})^{ \varrho_-/
\varrho
_+} a_\ell^{1-s} e^{ \varrho_- (1 - s) S_{\ell-1}}
\Biggr)^{k^*-1},\\
&&\hspace*{294pt}\tau_0^-<\tau_L^+\Biggr].
\end{eqnarray*}

Again, we apply (\ref{est9sub}) with $\alpha= \varrho_-/ \varrho_+$ to
$(S_\ell- S_{\ell-1}, a_\ell^{1-s})_{\ell\ge1}$ with $\gamma=
\varrho
_+ - \varrho_-, b= \varrho_- (1-s), \eta= \varrho_-, p=k^*-1>
\gamma/b$
(the integrability hypothesis can be easily checked as before), which
yields that $\e_x[(Z_g[0,L])^{k^* }] \le c' \lambda^{s(k^*-1)} e^{
\varrho_+ x + (k^* \varrho_- - \varrho_+) L} $, proving (i) in the case
that $\varrho_+/\varrho_-$ is an integer.

(ii) Write in this proof $\Lambda:= \sum_{u \in
\H(L) \cap\gg(L,\lambda) } e^{\varrho_- V(u)} $. Instead of $\Q
_x^{(\varrho_-)}$, we shall make use of the probability $\Q$ defined in
(\ref{Q1}) with $\varrho=\varrho_+$ for the change of measure. We
stress that under $\Q$, $(S_n)$ drifts to $+\infty$.

Firstly, we prove by induction on $k$ that for any $1\le k \le k^* -1$,
there exists some constant $c_k=c_k(\lambda)>0$ such that
%
\begin{equation}\label{toto10} \e_x\bigl[\Lambda^k\bigr] \le c_k
e^{\varrho_+ x} e^{ (k \varrho_- - \varrho_+)L}.
\end{equation}

By the change of measure, we get that for $k\ge1$,
%
\begin{eqnarray}\label{newsub1}\qquad  \e_x\bigl[\Lambda^k\bigr]
&= & e^{\varrho_+ x} \Q_x \bigl[e^{ (\varrho_- - \varrho_+)
S_{\tau
_L^+}} 1_{\{ w_{\tau_L^+} \in\gg(L,\lambda) \}} \Lambda^{k-1},
\tau
_L^+<\tau_0^-\bigr]
\nonumber
\\[-8pt]
\\[-8pt]
\nonumber
&=& e^{\varrho_+ x + (\varrho_-- \varrho_+) L} \Q_x \bigl[ e^{
(\varrho
_- - \varrho_+ )T_L^+} 1_{\{ w_{\tau_L^+} \in\gg(L,\lambda) \}}
\Lambda
^{k-1}, \tau_L^+<\tau_0^-\bigr],
\end{eqnarray}
where $T_L^+:=S_{\tau_L^+} - L>0$. This yields the case $k=1$
of (\ref{toto10}).

Assume $ 2\le k \le k^*-1 $ and that (\ref{toto10}) holds for $1,
\ldots, k-1$.
Exactly as before, we decompose $\Lambda$ along the spine up to $\tau
_L^+$, apply the triangular inequality and arrive at
\begin{eqnarray*}
&&\bigl( \Q_x\bigl [\Lambda^{k-1} | \G_\infty\bigr]\bigr)^{1/(k-1)}\\
&&\qquad \le e^{ \varrho_-
S_{\tau
_L^+}} + \sum_{\ell=1}^{\tau_L^+} \sum_{u \in\I_\ell} \bigl( \Q_x
\bigl[\bigl(\Lambda
^{(u)}\bigr)^{k-1} | \G_\infty\bigr]\bigr)^{1/(k-1)},
\end{eqnarray*}
where $\Lambda^{(u)}:= \sum_{v \in\cT^{(u)} \cap
\H(L) \cap\gg(L,\lambda) } e^{\varrho_- V(v)} $ with $\cT^{(u)}$
the subtree rooted~at~$u$. By Proposition \ref{pchange-proba}, under
$\Q_x$ and conditioning on $\G_\infty$, each $\Lambda^{(u)}$ is
distributed as $(\Lambda,
\mathbf{P}_{V(u)})$. Hence by induction assumption, $ ( \Q_x
[(\Lambda
^{(u)})^{k-1} | \G_\infty])^{1/(k-1)} \le c_{k-1}^{1/ (k-1)} e^{
\varrho_+ (V(u) - L)/ (k-1)} e^{\varrho_ - L}. $ Then
\begin{eqnarray*}
&&\bigl( \Q_x \bigl[\Lambda^{k-1} | \G_\infty\bigr]\bigr)^{1/(k-1)} \\
&&\qquad\le e^{ \varrho_-
S_{\tau
_L^+}} + c_{k-1}^{1/ (k-1)} e^{\varrho_ - L} \sum_{\ell=1}^{\tau
_L^+} \sum_{u \in\I_\ell} e^{ \varrho_+ \Delta V(u) / (k-1)} e^{
\varrho_+ (S_{\ell-1} -L)/ (k-1)}.
\end{eqnarray*}

Notice that $ {\varrho_+ \over k-1} \ge\varrho_-$ and that on $\{
w_{\tau_L^+} \in\gg(L,\lambda) \}$,
\begin{eqnarray*}
\sum_{u \in\I_\ell} e^{{ \varrho_+ / (k-1)} \Delta V(u)} &\le&
a_\ell
\max_{ u \in\I_\ell} e^{({ \varrho_+ / (k-1)} - \varrho_-)
\Delta
V(u) } \\
&\le&(a_\ell)^{1-s} \lambda^{ {( \varrho_+ /(\varrho_-(
k-1)))} - (1-s)} e^{ ( {\varrho_+ / (k-1)}- (1-s) \varrho_- )
(L-S_{\ell-1})},
\end{eqnarray*}
with $s:= { k^*- { \varrho_+/\varrho_ - } - \delta_2
\over k^*-1}$. We mention that the above inequality holds for \mbox{$k=k^*$}.

Going back to (\ref{newsub1}), we obtain that [we keep the density
there $e^{(\varrho_- -\varrho_+) T_L^+}$ only for $e^{\varrho_-
S_{\tau
_L^+}}$ and use the inequality $(x+y)^{k-1} \le2^{k-1} (x^{k-1}+ y^{k-1})$]
\begin{eqnarray*}  &&\e_x[\Lambda^{k}] \\
&&\quad\le c(\lambda)
e^{\varrho_+
x + (\varrho_- - \varrho_+) L}\\
&&\qquad{}\times e^{\varrho_ - (k-1) L} \Biggl( \Q_x
\bigl[e^{ ( k\varrho_- - \varrho_+)T_L^+ }\bigr]+ \Q_x\Biggl[\sum_{\ell
=1}^{\tau
_L^+} (a_\ell)^{1-s} e^{ (1-s) \varrho_- (S_{\ell-1} - L)}
\Biggr]^{k-1}\Biggr).
\end{eqnarray*}

Recall that $\Q_x[ e^{ ( k\varrho_- - \varrho_+) T_L^+}] =
\Q
[e^{ ( k\varrho_- - \varrho_+)T_{L-x}^+}]$ is bounded by some constant
since we have $\Q[ e^{ (k\varrho_- - \varrho_+ + \delta) S_1 }] =
\exp
\{ \psi( k\varrho_-+\delta) \} < \infty$ if $\delta>0$ is sufficiently
small [here we use the fact that $k\le k^*-1$]. By Lemma \ref
{lpositivewalk}, the above expectation $\Q_x [\cdot\cdot\cdot]^{k-1}$
is uniformly bounded, which proves (\ref{toto10}).\

To control $\e_x[\Lambda^q]$, we use the change of measure
\[
\e_x\bigl[\Lambda^q\bigr]=e^{\varrho_+ x + (\varrho_-- \varrho_+) L} \Q_x
\bigl[
e^{ (\varrho_- - \varrho_+ )T_L^+} 1_{\{ w_{\tau_L^+} \in\gg
(L,\lambda
) \}} \Lambda^{q-1}, \tau_L^+<\tau_0^-\bigr].
\]

Since $q < k^*$, $ ( \Q_x [\Lambda^{q-1} | \G_\infty
])^{1/(q-1)} \le( \Q_x [\Lambda^{k^*-1} | \G_\infty])^{1/(k^*-1)}$.
From (\ref{toto10}) with $k=k^*-1$ there, we use the same arguments as
before and get that
\begin{eqnarray*}  \e_x[\Lambda^{q}] &\le& c e^{\varrho_+ x +
(\varrho_- - \varrho_+) L} e^{\varrho_ - (q-1) L} \Biggl( \Q_x
\bigl[e^{ (
q \varrho_- - \varrho_+)T_L^+ } \bigr]\\
&&\hspace*{120pt}{}+ \Q_x\Biggl[\sum_{\ell=1}^{\tau_L^+}
(a_\ell)^{1-s} e^{ (1-s) \varrho_- (S_{\ell-1} - L)}
\Biggr]^{q-1}\Biggr).
\end{eqnarray*}

Again, $\Q_x [e^{ ( q\varrho_- - \varrho_+) T_L^+}]$ is bounded by
some constant since
\[
\Q\bigl[ e^{ ( q \varrho_- - \varrho_+ +\delta) S_1 }\bigr]
= \exp\bigl( \psi( q \varrho_-+\delta)\bigr) < \infty
\]
if $\delta>0$ is
sufficiently small. By Lemma \ref{lpositivewalk}, the above
expectation $\Q_x [\cdot\cdot\cdot]^{q-1}$ is uniformly bounded, which
proves (ii).

(iii) The proof follows in the same spirit as that of (i) and (ii):
Let $ \chi(L):= \sum_{u \in
\H(L) } e^{\varrho_- (V(u)-L) }$, and we prove by induction that for
any $1\le k \le k^*$,
%
\begin{equation}\label{toto20} \e_x\bigl[\chi(L)^k\bigr] \le c_k
e^{\varrho_+ (x - L) },\qquad x\in\r.
\end{equation}

The case $k=1$ is obvious by the change of measure. Assume (\ref
{toto20}) for $k-1$ and $2\le k \le k^*$. By repeating the same
arguments as in (ii), we get that
%
\begin{eqnarray}\label{toto25}
\e_x\bigl[ \chi(L)^k\bigr] &\le& c e^{ \varrho_- (x- L) }
\nonumber
\\
&&{}\times
\Biggl( \Q
^{(\varrho_-)}_x \bigl[e^{ (k-1) \varrho_- T_L^+}, \tau_L^+< \tau_0^-\bigr] \\
&&\hspace*{18pt}{}+
\Q
^{(\varrho_-)}_x \Biggl[ \Biggl( \sum_{ \ell=1}^{ \tau_L^+} \sum_{u
\in\I
_\ell} e^{ {(\varrho_+/ (k-1))} ( V(u) -L)} \Biggr)^{k-1}, \tau_L^+<
\tau
_0^- \Biggr] \Biggr).\nonumber
\end{eqnarray}

By the absolute continuity between $\Q_x^{(\varrho_-)}$ and
$\Q_x$,
\begin{eqnarray*}
\Q^{(\varrho_-)}_x \bigl[e^{ (k-1) \varrho_- T_L^+}, \tau_L^+< \tau
_0^-\bigr] &
=& e^{ (\varrho_+-\varrho_-) x - (k-1) \varrho_- L} \Q_x \bigl[e^{ (k
\varrho
_- - \varrho_+) S_{\tau_L^+}}, \tau_L^+< \tau_0^-\bigr] \\
&= & e^{ (\varrho_+-\varrho_-) (x - L)} \Q_x \bigl[e^{ (k \varrho_- -
\varrho_+) T_L^+}, \tau_L^+< \tau_0^-\bigr] \\
&\le& c e^{ (\varrho_+-\varrho_-) (x - L)},
\end{eqnarray*}
where the term $\Q_x [e^{( k \varrho_- - \varrho_+) T_L^+}]$
is uniformly bounded, since for $k\le k^*$ and sufficiently small
$\delta_4>0$, $\Q[ e^{ (k \varrho_- - \varrho_+ + \delta_4)S_1}] = e^{
\psi( k \varrho_- + \delta_4)} <\infty$ by (\ref{hyp-extra}).

It remains to control the second expectation term $\Q_x^{(\varrho_-)}$
in (\ref{toto25}). Let $b_\ell:= \sum_{u \in\I_\ell} e^{ {\varrho
_+/ (k-1)} \Delta V(u) }$, for $ \ell\ge1$.
Under $\Q^{(\varrho_-)}_x$, $( S _\ell- S_{ \ell-1}, b_\ell)_{\ell
\ge
1}$ are i.i.d. and
\begin{eqnarray*}
\Q^{(\varrho_-)}\bigl[b_1^{k-1}\bigr] &=& \e\biggl[ \biggl(\sum_{\vert u \vert=1} e^{
\varrho_- V(u)} \biggr) \biggl( \sum_{v \not=u} e^{ {(\varrho_+/( k-1))} V(v)}
\biggr)^{k-1} \biggr]\\
& \le&\e\biggl[ \sum_{\vert u \vert=1} e^{ \varrho_-
V(u)} \biggr] ^{ 1+ {\varrho_+/\varrho_-}},
\end{eqnarray*}
where the last inequality follows from the elementary
inequality: for any $n\ge1$ and $x_1, \ldots, x_n \in\r$, $(
\sum
_{i=1}^n e^{ \varrho_+ x_i/ (k-1)})^{k-1} \le( \sum
_{i=1}^n e^{ {\varrho_- /r\varrho_+} \varrho_+ x_i})^{
\varrho
_+/ \varrho_-}= \break ( \sum_{i=1}^n e^{ \varrho_- x_i})^{
\varrho
_+/ \varrho_-}$, since $k-1 < {\varrho_+\over\varrho_-}$. Then $\Q
^{(\varrho_-)}[b_1^{k-1}] < \infty$ by (\ref{hyp-extra}). Going back to
(\ref{toto25}), we see that the expectation term $\Q_x[(\cdot\cdot
)^{k-1}, \tau_L^+< \tau_0^-]$ equals
\[
\Q^{(\varrho_-)}_x \Biggl[ \Biggl( \sum_{\ell=1}^{ \tau_L^+} b_\ell
e^{
{\varrho_+/ (k-1)} ( S_{\ell-1} -L)} \Biggr)^{k-1}, \tau_L^+< \tau_0^-
\Biggr] \le c' e^{ (\varrho_+-\varrho_-) (x - L)},
\]
by applying (\ref{est9sub}) to $( S_\ell- S_{ \ell-1},
b_\ell
)_{\ell\ge1}$ with $\gamma= \varrho_+- \varrho_-$, $b= \varrho
_+/(k-1)$ and $p=k-1$. This proves (\ref{toto20}) hence (iii).
\end{pf*}

The next lemma controls the number of bad particles.

\begin{lemma}\label{lsub-bad} Let $ r= { \varrho_+ \over\varrho_- }
-1 + {\delta^*\over2}$ [with $\delta^*$ as in (\ref{hypsub})].
\begin{longlist}[(ii)]
\item[(i)] There exists some constant $c=c(r)>0$ such that
for all $0\le x \le L$,
\[
\e_x\bigl[Z_b[0,L] \bigr] \le c \lambda^{-r} e^{\varrho_+ x}
e^{(\varrho_- - \varrho_+)L}.
\]

\item[(ii)] Denote by $\cL_{b, L} [0]:=\{v \in\cL[0]\dvtx\exists u \in\H(L)
\cap\bb(L, \lambda) \mbox{ with } u < v\} $ the set of leaves which
are descendants of some element of $ \H(L) \cap\bb(L, \lambda)$. Then
for any $0\le x \le L$,
\[
\e_x \bigl[ \#\cL_{b, L}[0]\bigr] \le c \lambda^{-r}
e^{\varrho_+
x} e^{(\varrho_- - \varrho_+)L}.
\]
\end{longlist}
\end{lemma}

\begin{pf}
(i) By changing the measure from $\mathbf{P}_x$ to $\Q_x^{(\varrho_-)}$,
\[
\e_x\bigl[Z_b[0,L] \bigr] = e^{\varrho_- x} \Q_x^{(\varrho
_-)}
\bigl[ e^{ - \varrho_- S_{\tau_0^-}} 1_{\{w_{\tau_0^-} \in\bb(L,
\lambda
)\}}, \tau_0^- < \tau_L^+\bigr].
\]

Let us write $a_j:= \sum_{u\in\I_j} (1+e^{\varrho_-
\Delta
V(u)})$, $j\ge1$, in this proof. Then
%
\begin{equation} \label{1w0} 1_{\{ w_{\tau_0^-} \in\bb(L,\lambda)
\}}
\le\sum_{j=1}^{\tau_0^-} \lambda^{- r} a_j ^r e^{- r \varrho_-
(L- S_{j-1})},
\end{equation}
which yields that
\begin{eqnarray*} \e_x\bigl[Z_b[0,L] \bigr] &\le& \lambda^{- r}
e^{\varrho_- x} \Q_x^{(\varrho_-)}\Biggl [ e^{ - \varrho_-
S_{\tau
_0^-}} \sum_{j=1}^{\tau_0^-} a_j ^r e^{- r \varrho_- (L-
S_{j-1})}, \tau_0^- < \tau_L^+\Biggr] \\
&\le& c \lambda^{- r} e^{\varrho_- x} e^{ (\varrho_+- \varrho_-)
(x-L) },
\end{eqnarray*}
by applying (\ref{est9sub}) to $\gamma=\varrho_+-\varrho_-,
p=1$ and $b= r \varrho_- > \gamma$ [the integrability hypothesis is
satisfied thanks to (\ref{hypsub}) and the choice of $r$: $ \Q
^{(\varrho
_-)}[ (1+1_{\{S_1<0\}} e^{- \varrho_- S_1}) a_1^r] \le
\e
[ \sum_{ \vert u \vert=1} (1+ e^{ \varrho_- V(u)} )]^{ r+1}
< \infty$, and $\Q^{(\varrho_-)}[ e^{ r \varrho_- S_1} ]=
e^{ \psi( \varrho_- ( 1+ r ))} < \infty$]. This proves (i).

(ii) Remark that $\#\cL_{b, L}[0]= \sum_{ u \in\H(L)\cap\bb(L,
\lambda
)} \#\cL^{(u)}[0]$, where $\cL^{(u)}[0]$ denotes the set of leaves
which are descendants of $u$. By the branching property, conditioned on
$ \H(L)\cap\bb(L, \lambda)$, $(\#\cL^{(u)}[0])_{u\in\H(L)\cap
\bb(L,
\lambda)} $ are independent and are distributed as $\#\cL[0]$ under
$\mathbf{P}
_{V(u)}$. It follows from Lemma \ref{lmomentsubL} (with $k=1$) that
\begin{eqnarray*} \e_x \bigl( \#\cL_{b, L}[0]\bigr)  &\le& c \e_x
\biggl[ \sum_{u \in\H(L)\cap\bb(L, \lambda)} e^{ \varrho_- V(u)}
\biggr] \\
&=& c e^{\varrho_- x} \Q_x^{(\varrho_-)}\bigl [ w_{\tau_L^+}
\in\bb
(L, \lambda), \tau_L^+ < \tau_0^-\bigr],
\end{eqnarray*}
by the change of measure from $\mathbf{P}_x$ to $ \Q
_x^{(\varrho
_-)}$. By (\ref{1w0}) (with $\tau_L^+$ instead of~$\tau_0^-$), the
above probability under $ \Q_x^{(\varrho_-)}$ is less than
\begin{eqnarray*} && \lambda^{- r} \Q_x^{(\varrho_-)} \Biggl[
\sum
_{j=1}^{\tau_L^+} a_j ^r e^{- r \varrho_- (L- S_{j-1})}, \tau_L^+
< \tau_0^-\Biggr] \\
&&\qquad\le \lambda^{- r} \sum_{j\ge1} \Q_x^{(\varrho_-)} \bigl[
e^{- r
\varrho_- (L- S_{j-1})}, j \le\min\bigl(\tau_L^+, \tau_0^-\bigr) \bigr]
\Q
_x^{(\varrho_-)}\bigl [a_j ^r\bigr],
\end{eqnarray*}
since for each $j$, $a_j$ is independent of $(S_{j-1}, j
\le\min(\tau_L^+, \tau_0^-))$; moreover $\Q_x^{(\varrho_-)} [a_j
^r]= \Q^{(\varrho_-)} [a_j ^r]=c' < \infty$ as in (i). Then we have
\[
\e_x\bigl[Z_b[0,L] \bigr] \le c c' e^{\varrho_- x} \lambda^{- r}
\sum
_{j\ge1} \Q_x^{(\varrho_-)} \bigl[ e^{- r \varrho_- (L-
S_{j-1})},
j \le\min\bigl(\tau_L^+, \tau_0^-\bigr) \bigr],
\]
which by an application of (\ref{est3sub}) (with $r \varrho
_- > \gamma:= \varrho_+-\varrho_-$) gives (ii).
\end{pf}

Let $M_{\infty}^{(\varrho_-)}$ be the almost sure limit of
$M_{n}^{(\varrho_-)}:= \sum_{\vert u \vert=n} e^{ \varrho_- V(u)}$. By
\cite{big77,lyons97},
$M_{\infty}^{(\varrho_-)}$ is almost surely positive on the event
$\{\cT=\infty\}$. From \cite{Liu00}, we know that there exists a
constant $c_{\varrho_-}$ such that
%
\begin{equation}\label{liuvarrho-}
\mathbf{P}\bigl(M_{\infty}^{(\varrho_-)}>t\bigr) \sim c_{\varrho_-}
t^{-\varrho_+/\varrho_-}, \qquad t \to\infty.
\end{equation}

We mention that the constant $c_{\varrho_-} $ is given
in \cite{JOC10}, Theorem 4.10:
\[
c_{\varrho_-}= {1\over
\varrho_+ \psi'(\varrho_+) } \e\biggl[ \biggl(
\sum_{|u|=1} e^{\varrho_- V(u)} M_{\infty}^{(\varrho_-, u)}\biggr)^
{\varrho_+/\varrho_-} - \sum_{|u|=1} e^{\varrho_+ V(u)}
\bigl(M_{\infty}^{(\varrho_-, u)}\bigr)^ {\varrho_+/\varrho_-} \biggr],
\]
where under $\mathbf{P}$ and conditioned on $\{ V(u), |u|=1\}$,
$(M_{\infty}^{(\varrho_-, u)})_{|u|=1}$ are i.i.d. copies of
$M_{\infty}^{(\varrho_-)}$.

\begin{lemma}[(Subcritical case)] \label{Lconjecture2c} As $t \to
\infty
$, the law of $\#\cL[0]$ under
$\mathbf{P}_t$, the number of descendants absorbed at $0$ of a particle
starting from $t$, normalized by $e^{\varrho_- t}$ converges in
distribution to $c_{\mathrm{sub}}^* M^{(\varrho_-)}_\infty$ where
\[
c_{\mathrm{sub}}^*
= { \Q^{(\varrho_-)}[
e^{-\varrho_- S_{\tau_0^-}} ] -1 \over\varrho_- \Q
^{(\varrho
_-)}[-S_{\tau_0^-} ] }.
\]
\end{lemma}

\begin{pf}
The proof is
similar to that of Lemma \ref{Lconjecture2}; we only point
out the main difference and omit the details. Recall that $ \cL[a]
:= \{ u\in\cT\dvtx|u|=\tau_a^-(u) \} $. By linear translation, it
is enough to prove that $ e^{- \varrho_- t} \# \cL[-t]$ converges
in law to $c_{\mathrm{sub}}^* M^{(\varrho_-)}_\infty$. Let
$M^{(\varrho_-)}_{\cL[-t]}:=\sum_{u\in\cL[-t]}e^{\varrho_- V(u)}$,
which converges~almost~surely to $M^{(\varrho_-)}_\infty$. On the
other hand, we have $M^{(\varrho_-)}_{\cL[-t]} =\break e^{- \varrho_- t}
\sum_{u\in\cL[-t]}e^{\varrho_- (V(u)+t)}$. Similarly to the proof of
Lemma \ref{Lconjecture2}, we apply Theorem 6.3 in Nerman \cite
{nerman81} (with $\alpha=\varrho_-$ there) and obtain that on $\{\cT
=\infty\}$, almost surely
\[
{ \sum_{u\in\cL[-t]}e^{\varrho_-
(V(u)+t)} \over\#\cL[-t]} \to \varrho_- {\Q^{(\varrho_-)}
[- S_{\tau_0^-}
] \over\Q^{(\varrho_-)}[ e^{-\varrho_-S_{\tau_0^-}}
] -1},\qquad
t \to\infty,
\]
which
easily yields the lemma.
\end{pf}

\begin{lemma} \label{Ltailsub} For any $\lambda>0$, let $\widehat
\mu
_{\lambda,\infty}:=\sum_{i=1}^{\hat\zeta_{\lambda}} \delta_{\{
x_i\}} $ be
the point process defined in Proposition \ref{pyaglom} associated
with $\B(\theta):= ( {1\over\lambda} \int\theta(dx) (1+
e^{\varrho_-
x}))^{1/\varrho_-}$ for $\theta\in\Omega_f$. Let $(M_{\infty
}^{(\varrho_-,i)},
i\ge1)$ be a sequence of i.i.d. random variables of common law that of
$(M_{\infty}^{(\varrho_-)}, \mathbf{P})$, independent of $\widehat
\mu
_{\lambda
,\infty}$. As $t
\to\infty$, we have
\[
\Q\Biggl(\sum_{i=1}^{\hat\zeta_\lambda} e^{ \varrho_- x_i}
M_\infty^{(\varrho_-,i)}
> t\Biggr) \sim c_{\varrho_-} \Q\biggl[ \int\widehat\mu
_{\lambda,
\infty} (dx) e^{ \varrho_+ x} \biggr] t^{- \varrho_+/\varrho_-}.
\]
\end{lemma}

We mention that as $\lambda\to\infty$, $\Q[ \int\widehat\mu
_{\lambda, \infty} (dx) e^{ \varrho_+ x} ] \to{1\over\Q
[\Re^{-1}
]}$ by (\ref{eqref-muinfty2})\break and~(\ref{eqref-R}).

\begin{pf}
Let $\Lambda_{L, \lambda}:= \sum_{u \in
\H(L) \cap\gg(L,\lambda) } e^{\varrho_- (V(u)-L)}$.\vspace*{-1.5pt} By Proposition
\ref{pyaglom}, under $\mathbf{P}_x(\cdot\vert H(L)>0)$, $\Lambda_{L,
\lambda
}$ converges in law to $ \int\widehat\mu_{\lambda, \infty} (dx) e^{
\varrho_- x}= \sum_{i=1}^{\hat\zeta_\lambda} e^{ \varrho_- x_i}$ (some
tightness is required here but we omit the details since the arguments
are similar to the critical case).
By Lemma \ref{lmomentsubZ}(ii), the family $( \Lambda_{L,
\lambda},\break
\mathbf{P}_x(\cdot\vert H(L)>0))$ is bounded in $L^{q}$ with $q={\rho
_+\over
\rho_-} +\delta_2$, hence
%
\begin{equation}\label{ahah3} \Q\Biggl[ \sum_{i=1}^{\hat\zeta
_\lambda}
e^{ \varrho_- x_i} \Biggr]^{q} < \infty.
\end{equation}

This together with (\ref{liuvarrho-}) allows us to apply
Lemma \ref{Lconvolution} to $p=
{\varrho_+\over\varrho_-}$ and yields the desired asymptotic result.
\end{pf}

We now prove Theorem \ref{tleaves} in the subcritical case.

\begin{pf*}{Proof of Theorem \ref{tleaves}(ii)}

\textit{Lower bound of Theorem \ref{tleaves}}(ii):
The proof of the lower bound goes in the same way as that of Theorem
\ref{tleaves}(i) by using Proposition \ref{pyaglom} and Lemma~\ref{Lconjecture2c}.
Let $A>0$. Consider $n
\to\infty$, let $L_A:={1\over\varrho_-} \log n -A$ and $\lambda:=
e^{\varrho_- A}$. We keep
the same notation $H_g(L_A), (\#\cL^{(i)}[0], 1\le i \le
H_g(L_A))$: Recall \eqref{defbad} and $H_{g}(L_A):=\#\H_{g}(L_A)$ with
%
\begin{equation} \label{defhgl2} \H_g(L_{ A}):= \H(L_{ A}) \cap\gg
\bigl(L_n, e^{ \varrho_- A}\bigr).
\end{equation}

We define as well $B^{(i)}:=\#\cL^{(i)}[0] e^{-\varrho_-
V(u^{(i)})} $ for $u^{(i)} \in\H(L_A)$, and $E(L_A)$ the event that $B^{(i)}
> (1-\varepsilon)M_{\infty}^{(\varrho_-,i)}$, $\forall i$ with
small $\varepsilon>0$. Repeating the proof of the lower bound of
Theorem \ref{tleaves}(i), and using Proposition \ref{pyaglom} and
Lemma~\ref{Lconjecture2c}, we get that for any
$A>0$,
%
\begin{eqnarray} \label{67}
&& \liminf_{n\to\infty} n^{\varrho_+/\varrho_-} \nonumber
\mathbf{P}_x\Biggl(\sum_{i=1}^{H_g(L_A)}\#\cL^{(i)}[0] > n \Biggr) \\
&&\qquad\ge {\Q[\Re^{-1}]\over
C_R }R(x) e^{\varrho_+ x} e^{\varrho_+ A}
\Q\Biggl(\sum_{i=1}^{\hat\zeta_\lambda} e^{ \varrho_- x_i}
M^{(\varrho_-,
i)}_\infty> {1\over c_{\mathrm{sub}}^*} e^{ \varrho_- A}\Biggr) \\
&&\qquad=: {\Q[\Re^{-1}]\over
C_R }R(x)e^{\varrho_+ x}C_s(A),\nonumber
\end{eqnarray}
where $\widehat\mu_{\lambda,\infty}:=\sum_{i=1}^{\hat
\zeta
_\lambda} \delta_{\{x_i\}} $ is
the point process as in Lemma \ref{Ltailsub} (with $\lambda:=
e^{\varrho_- A}$ there) and $c_{\mathrm{sub}}^*$ is defined in Lemma \ref
{Lconjecture2c}. The same also holds for the upper bound, hence for
any $A>0$,
%
\begin{equation}\label{subHG}\qquad  \lim_{n\to\infty} n^{\varrho_+/
\varrho_-}
\mathbf{P}_x\Biggl(\sum_{i=1}^{H_g(L_A)}\#\cL^{(i)}[0] > n \Biggr) = {\Q
[\Re
^{-1}]\over
C_R }R(x)e^{\varrho_+ x}C_s(A).
\end{equation}

Since $\mathbf{P}_x ( \#\cL[0] >n )\ge
\mathbf{P}_x(\sum_{i=1}^{H(L_A)}\#\cL^{(i)}[0] >n )$, we get
that for
any \mbox{$A>0$},
%
\begin{equation}\label{ahah1}
\liminf_{n\to\infty} n^{\varrho_+/\varrho_-} \mathbf{P}_x \bigl( \#
\cL[0] >n
\bigr) \ge{\Q[\Re^{-1}]\over
C_R }R(x)e^{\varrho_+ x}C_s(A).
\end{equation}

\textit{Upper bound of Theorem \ref{tleaves}}(ii): By Lemma \ref
{lsub-bad} and Lemma \ref{lmomentsubZ}(i) with $L:=L_A= {1\over
\varrho_-} \log n - A$, $\lambda:= e^{ \varrho_- A}$ and $k^*:=
\lfloor
{ \varrho_+\over\varrho_-} \rfloor+1$, we obtain the following
estimate: For any $\varepsilon>0$,
\begin{eqnarray*}
\mathbf{P}_x \bigl( Z_g[0, L_A] \ge\varepsilon n \bigr) &\le&( \varepsilon
n)^{-k^*} c e^{ A(\varrho_- k^*- \varrho_+- \delta_2 \varrho_-)}
e^{\varrho_+ x + (\varrho_- k^*- \varrho_+) L_A}\\
& =& c_{\varepsilon,
x}
n^{ - \varrho_+/\varrho_-} e^{ - \delta_2 \varrho_- A}
\end{eqnarray*}
and
\begin{eqnarray*}
\mathbf{P}_x \bigl( Z_b[0, L_A] \ge\varepsilon n \bigr) &\le&{1\over
\varepsilon
n} c e^{ - A(\varrho_+-\varrho_-+ \delta^* \varrho_-/2)}
e^{\varrho
_+ x + (\varrho_- - \varrho_+) L_A} \\
&=& c_{\varepsilon, x} n^{ -
\varrho
_+/\varrho_-} e^{ - \delta^* \varrho_- A/2},
\end{eqnarray*}
with the same estimate for $ \mathbf{P}_x ( \cL_{b,
L_A}[0] \ge
\varepsilon n ) $. Since $Z[0, L_A]= Z_g[0, L_A]+Z_b[0, L_A]$, we
obtain that for any $\varepsilon>0$,
\[
\limsup_{ A \to\infty} \limsup_{n\to\infty} n^{ \varrho
_+/\varrho_-}
\mathbf{P}_x \bigl( Z[0, L_A] + \cL_{b, L_A}[0] \ge3 \varepsilon n
\bigr) =0.
\]

From here and using the fact that $ \#\cL[0] = Z[0, L_A] +
\cL
_{b, L_A}[0] +\break  \sum_{i=1}^{H_g(L_A)}\#\cL^{(i)}[0] $, we deduce from
(\ref{subHG}) that for any $A>0$,
\[
\limsup_{n\to\infty} n^{\varrho_+/\varrho_-} \mathbf{P}_x \bigl( \#
\cL[0] >n
\bigr) \le{\Q[\Re^{-1}]\over
C_R }R(x)e^{\varrho_+ x}C_s(A) + o_A(1),
\]
with $o_A(1) \to0$ as $A \to\infty$ (in fact exponentially
fast). This together with the lower bound (\ref{ahah1}) yields that $
\lim_{n\to\infty} n^{\varrho_+/\varrho_-} \mathbf{P}_x ( \#\cL
[0] >n )$
exists and is finite. Then, a fortiori, $\lim_{A\to\infty}C_s(A) $
also exists and is some finite constant. This proves Theorem \ref
{tleaves}(ii).
\end{pf*}

We end this section by giving the proof of Lemma \ref{Lconst}.

\begin{pf*}{Proof of Lemma \ref{Lconst}}
By (\ref{defcr}), $C_R=
1/\Q(\tau_0^-=\infty)$. Recall (\ref{67}). It suffices to show that
%
\begin{equation}\label{ahah2} \lim_{A\to\infty}C_s(A) = {c_{\varrho_-}
\over\Q[\Re^{-1}]}
\bigl(c_{\mathrm{sub}}^*\bigr)^{\varrho_+/\varrho_-}.
\end{equation}

The lower bound follows from the monotonicity: the random point measure
$\widehat\mu_{\lambda, \infty}$ is stochastically increasing in
$\lambda$; Then for any $A > A_0$, $\lambda= e^{ \varrho_-
A}>\lambda
_0:= e^{\varrho_- A_0}, $ $\widehat\mu_{\lambda, \infty}$
stochastically dominates $\widehat\mu_{\lambda_0, \infty}:= \sum
_{i=1}^{\hat\zeta_{\lambda_0}} \delta_{\{x_i\}}$, hence
\[
C_s(A) \ge e^{\varrho_+ A}
\Q\Biggl(\sum_{i=1}^{\hat\zeta_{\lambda_0}} e^{ \varrho_- x_i}
M^{(\varrho_-,
i)}_\infty> {1\over c_{\mathrm{sub}}^*} e^{ \varrho_- A}\Biggr).
\]
Applying Lemma \ref{Ltailsub} to $\lambda_0$ yields that
for any $\lambda_0:=e^{ \varrho_- A_0}$,
\[
\liminf_{A\to\infty} C_s(A) \ge c_{\varrho_-} \Q\biggl[ \int
\widehat\mu_{\lambda_0, \infty} (dx) e^{\varrho_+ x}\biggr]
\bigl(c_{\mathrm{sub}}^*\bigr)^{\varrho_+/\varrho_-}.
\]
Letting $\lambda_0 \to\infty$, the above expectation term
converges to $1/\Q[\Re^{-1}]$ and proves the lower bound.

To derive the upper bound, by Lemma \ref{lmomentsubZ}(iii) and
Theorem \ref{thmyaglom-critiq}(ii), we get that under $\mathbf
{P}(\cdot
| \H
(L)>0)$, $ \sum_{u \in\H(L)} e^{\varrho_- (V(u) -L)} $ is bounded in
$L^{k^*}$ and converges in law to $\sum_{i=1}^{\hat\zeta_\infty} e^{
\varrho_- x_i} $, where $\widehat\mu_\infty= \sum_{i=1}^{\hat
\zeta
_\infty} \delta_{\{x_i\}}$. Therefore
\[
\Q\Biggl[ \sum_{i=1}^{\hat\zeta_\infty} e^{ \varrho_- x_i}
\Biggr]^{k^* } < \infty,
\]
which in view of Lemma \ref{Lconvolution} and (\ref
{liuvarrho-}) yields, as $A \to\infty$,
\[
e^{\varrho_+ A}
\Q\Biggl(\sum_{i=1}^{\hat\zeta_\infty} e^{ \varrho_- x_i}
M^{(\varrho_-,
i)}_\infty> {1\over c_{\mathrm{sub}}^*} e^{ \varrho_- A}\Biggr) \to
{c_{\varrho_-} \over\Q[\Re^{-1}]} \bigl(c_{\mathrm{sub}}^*\bigr)^{\varrho
_+/\varrho
_-}.
\]
Since $\widehat\mu_\infty$ stochastically dominates
$\widehat\mu_{A, \infty}$, this gives the desired upper bound for
$C_s(A)$ and completes the proof of the lemma.
\end{pf*}






\section{Proofs of the technical lemmas}
\label{proofsRW}

\subsection{\texorpdfstring{Proof of Lemma \protect\ref{lemlaw-convergence-reverse}}{Proof of Lemma 4}}
Obviously we may assume that $\|F\|_\infty\le1$ throughout the proof
of (i) and (ii).

\textit{Proof of part} (i). Since $\ppp( \tau_t^+ > K)
\to
1$ as $t \to\infty$, it is enough to show that
%
\begin{eqnarray}\label{i1} &&\lim_{t \to\infty} \eee\bigl[ 1_{\{\tau_t^+
> K\}} F\bigl( T_t^+,(S_{\tau_t^+}-S_{\tau_t^+-j})_{1\le j \le K}\bigr) \bigr]
\nonumber
\\[-8pt]
\\[-8pt]
\nonumber
&&\qquad=
\eee\bigl[ F\bigl(U \hat S_{\hat\sigma}, (\hat S_j)_{1 \le j \le K}\bigr) \bigr].
\end{eqnarray}

Recall
that $(\sigma_n, H_n)_{n\ge1}$ are the strict ascending ladder
epochs and ladder heights of $S$. Since for
some (unique) $n \ge1$, $\tau_{t}^+=\sigma_n$ and $T_t^+= H_n -t$, we
can write
\begin{eqnarray*}
B_t &:=& \eee\bigl[ 1_{\{\tau_t^+ > K\}} F\bigl( T_t^+,(S_{\tau
_t^+}-S_{\tau
_t^+-j})_{1\le j \le K}\bigr) \bigr]
\\
&=& \sum_{n \ge1} \eee\bigl[ 1_{\{H_{n-1} \le t < H_n \}}
1_{\{K< \sigma_n \}} F\bigl(H_n -t, (S_{\sigma_n}-S_{\sigma_n-j})_{1\le j
\le K}\bigr) \bigr].
\end{eqnarray*}
Let us choose some integer $m>K$. Notice that $\sigma_n-\sigma_{n-m}>K$
and $\sigma_n>K$ for $n \ge m$. Since the previous sum for $n < m$ is
smaller than $\ppp(H_m >t)$ which tends to $0$ when $t$ tends to
infinity, we get
\begin{eqnarray*}
B_t&=&\sum_{n \ge m} \eee\bigl[ 1_{\{H_{n-1} \le t < H_n
\}} F\bigl(H_n -t, (S_{\sigma_n}-S_{\sigma_n-j})_{1\le j \le K}\bigr) \bigr] + o_t(1)
\\
&=:& B'_t + o_t(1),
\end{eqnarray*}
with $|o_t(1)|\le\mathbf{P}(H_m >t) \to0$ as $t \to\infty$.
Applying the
strong Markov property at the stopping time
$\sigma_{n-m}$, we obtain that
\begin{eqnarray*}
B'_t&=&\sum_{n \ge m} \eee\bigl[ 1_{\{H_{n-m} \le t \}}
\eee_{H_{n-m}}\bigl[ 1_{\{H_{m-1} \le t < H_m \}}
F\bigl(H_m -t,(S_{\sigma_m}-S_{\sigma_m-j})_{1\le j \le K}\bigr) \bigr] \bigr]
\\
&=&\sum_{n \ge m} \eee\bigl[ 1_{\{H_{n-m} \le t \}} g(t- H_{n-m})
\bigr],
\end{eqnarray*}
with
\[
g(x):= \eee\bigl[ 1_{\{H_{m-1} \le x < H_m \}}
F\bigl(H_m -x, (S_{\sigma_m}-S_{\sigma_m-j})_{1\le j \le K}\bigr) \bigr]\qquad
\forall x\ge0.
\]

Therefore
%
\begin{equation}\label{Bt} B'_t=\int_{0}^{t} g(t-x) \,{d}u (x),
\end{equation}
with $u(x)= \sum_{n
\ge0} \ppp(H_n \le x). $ Let us check that $g$ is directly Riemann
integrable on $\r_+$.
Recall that a function $g$ is directly Riemann integrable (see Feller
\cite{feller}, page 362) if $g$ is continuous almost everywhere and satisfies
%
\begin{equation}
\label{eqdRi1}
\sum_{n=0}^\infty\sup_{n \le x \le n+1} \bigl|g(x)\bigr| <\infty.
\end{equation}
Observe first that $\|F\|_\infty\le1$ implies $\| g \|_\infty\le1$.
Now recall that $H_1$ is integrable. Therefore,
\[
\sum_{n\ge0} \sup_{n \le x \le n+1} \bigl|g(x)\bigr| \le\sum_{n\ge0} \ppp(H_m
\ge n) =1+ \eee[ H_m ] = 1+m \eee[H_1]<\infty,
\]
%
yielding (\ref{eqdRi1}). Now we prove that $g$ is a.e. continuous.
For $z \in\r_+^K$, denote by $D(z) \subset\r_+^*$ the set on
which $F(\cdot, z)$ is discontinuous. By assumption, $D(z)$ is at most
countable for any real $z$, hence
$D((S_{\sigma_m}-S_{\sigma_m-j})_{1\le j \le K})$ is a random set
(maybe empty) at most countable; the same is true for the random
set
\[
\Upsilon:= \bigcup_{n=1}^\infty\bigl\{ H_n - z\dvtx z \in D\bigl((S_{\sigma
_m}-S_{\sigma_m-j})_{1\le j \le
K}\bigr) \cup\{0\}\bigr\}.
\]

In other words, we may represent $\Upsilon$ by a sequence
of random variables taking values in $\r$. It follows that
\[
{\mathscr D}:= \bigl\{y\dvtx\ppp( y \in\Upsilon)
>0 \bigr\}\qquad \mbox{is at most countable.}
\]

We claim that for any $x \in\r_+^*\setminus{\mathscr D}$, $g$ is
continuous at $x$. In fact, for any sequence $(x_n)_n$ such that $x_n
\to x $ as $n \to\infty$, let $\xi_n:= 1_{\{H_{m-1} \le x_n < H_m \}}
F(H_m -x_n, (S_{\sigma_m}-S_{\sigma_m-j})_{1\le j \le K})$ and
$\xi:=1_{\{H_{m-1} \le x < H_m \}}F(H_m -x,
(S_{\sigma_m}-\break S_{\sigma_m-j})_{1\le j \le K})$, we shall show that
as $n \to\infty$,
%
\begin{equation}\label{xi311} \xi_n \to\xi\qquad
\mathrm{a.s.},
\end{equation}
which in view of the dominated
convergence theorem, implies that $g(x_n) \to g(x)$ and the desired
continuity of $g$ at $x$. To prove (\ref{xi311}), first we remark
that
%
\begin{eqnarray}\label{10as}\quad  \limsup_{n\to\infty} \bigl\vert1_{\{H_{m-1}
\le x_n < H_m \}} - 1_{\{H_{m-1} \le x < H_m
\}} \bigr\vert&\le&1_{\{H_{m-1} = x \}} + 1_{\{H_m = x
\}}
\nonumber
\\[-8pt]
\\[-8pt]
\nonumber
& =&0\qquad \mathrm{a.s.},
\end{eqnarray}
since $x \notin{\mathscr D} $ [hence a fortiori
$\ppp(H_n=x)=0$ for all $n \ge1$]. Second,
\[
\ppp\bigl( H_m -x \in
D\bigl((S_{\sigma_m}-S_{\sigma_m-j})_{1\le j \le K}\bigr) \bigr) \le\ppp( x \in
\Upsilon) =0,
\]
since $x \notin{\mathscr D}$. In words, almost surely,
$H_m-x \notin D((S_{\sigma_m}-S_{\sigma_m-j})_{1\le j \le K}) $,
which implies that $F(\cdot, (S_{\sigma_m}-S_{\sigma_m-j})_{1\le j
\le
K})$ is continuous at $H_m -x$; hence $F(H_m -x_n, (S_{\sigma
_m}-S_{\sigma_m-j})_{1\le j \le K}) \to F(H_m -x, (S_{\sigma
_m}-S_{\sigma_m-j})_{1\le j \le K})$ a.s. when $n \to\infty$. This and
(\ref{10as}) yield (\ref{xi311}) and the continuity of $g$ on $\r
_+^*\setminus{\mathscr D} $. Then $g$ is directly Riemann integrable.


Going back to (\ref{Bt}), we apply the renewal
theorem (see Feller \cite{feller}, page~363) and obtain that
\[
\lim_{t\to\infty} B'_t =\frac{1}{\eee[ H_1 ] }
\int_{0}^{\infty} g(x) \,{d}x,
\]
which implies
\begin{eqnarray*}
\lim_{t\to\infty} B_t &=& \frac{1}{\eee[ H_1 ] } \eee
\biggl[
\int_0^{ H_m- H_{m-1}}
F\bigl(H_m- H_{m-1}-x, (S_{\sigma_m}-S_{\sigma_m-j})_{1\le j \le
K}
\bigr) \,d x
\biggr] \\
&=& \frac{1}{\eee[ H_1 ] } \eee\bigl[ (H_m- H_{m-1})
F\bigl( U (H_m- H_{m-1}), (S_{\sigma_m}-S_{\sigma_m-j})_{1\le j \le
K}\bigr)
\bigr],
\end{eqnarray*}
by using the independent uniform variable $U$.

Finally since the random segments
$\{(S_{\sigma_k+j}-S_{\sigma_k})_{0 \le j \le
\sigma_{k+1}-\sigma_k };\break 0 \le k <m\}$ are i.i.d.,
Tanaka's construction [see (\ref{eqtanakaconstruction})] implies that
under $\ppp$ the segment of the
random walk $(S_n)_{n\ge0}$ up to time $\sigma_m$ viewed from
$(\sigma_m,S_{\sigma_m})$ in reversed time and reflected in the
$x$-axis, that is, $(S_{\sigma_m}-S_{\sigma_m-j})_{0\le j \le K}$, has
the same law as $(\zeta_j)_{0 \le j \le K}$. Moreover
since with this ``partial'' construction $H_m-H_{m-1}$ corresponds
to the value of the reversed and reflected process at time $\widetilde
\sigma=
\sup\{n \ge1 \dvtx\zeta_n=\min_{1 \le i \le n} \zeta_i\}$, we
obtain that
\begin{eqnarray*}
&& \frac{1}{\eee[ H_1 ] } \eee\bigl[ (H_m- H_{m-1})
F\bigl( U (H_m- H_{m-1}), (S_{\sigma_m}-S_{\sigma_m-j})_{1\le j \le
K}\bigr)
\bigr] \\
&&\qquad= \frac{1}{\eee[ H_1 ] } \eee\bigl[ \zeta
_{\widetilde\sigma
}
F\bigl( U \zeta_{\widetilde\sigma}, (\zeta_j)_{1\le j \le K}\bigr)
\bigr]= \eee\bigl[ F\bigl(U \hat S_{\hat\sigma}, (\hat S_j)_{1 \le j \le K}\bigr)
\bigr],
\end{eqnarray*}
by using (\ref{eqdef-xi}). This proves (\ref{i1}) and part
{(i)} of the lemma.

\textit{Proof of part} {(ii)}. Write for notational
convenience $\widetilde S^{(t)}_j:= S_{\tau_t^+} - S_{\tau_t^+-j}$ when
$1\le j \le\tau_t^+$.
Note that part (i) of the lemma implies
%
\begin{equation}
\label{eqci-conseq-part-i} \qquad\lim_{L \to\infty}\eee\bigl[ 1_{\{K<
\tau
_{L}^+\}} F \bigl(T_L^+,\bigl(\widetilde S^{(L)}_j\bigr)_{1\le j \le K}\bigr)
\bigr] =\eee\bigl[F \bigl(U \hat S_{\hat\sigma},(\hat S_j)_{1\le j \le
K}\bigr) \bigr]=: C_{ F}.
\end{equation}
Using the absolute continuity between $\ppp^{+}$ and $\ppp$ up to the
stopping time $\tau_t^+$ [the martingale $(R(S_j)1_{(j < \tau^+_t)},
j\le\tau_t^+)$ is uniformly integrable thanks to Lemma \ref
{lemRW}(ii) and (iv)], we can write
\[
\eee^{+} \bigl[ 1_{\{K< \tau_t^+\}} F \bigl(T_t^+,\bigl(\widetilde
S^{(t)}_j\bigr)_{1 \le j \le K}\bigr) \bigr] =\eee\bigl[R(S_{\tau_t^+}) 1_{\{K<
\tau_t^+<\tau_0^-\}} F \bigl(T_t^+,\bigl(\widetilde S^{(t)}_j\bigr)_{1 \le j
\le K}\bigr) \bigr].
\]

We treat first the case $\e[S_1]=0$. Combining parts (iii) and (iv)
of Lemma \ref{lemRW}, we deduce from the above equality that as $t
\to
\infty$,
\begin{eqnarray*}
\eee^{+} \bigl[ 1_{\{K< \tau_t^+\}} F \bigl(T_t^+,\bigl(\widetilde
S^{(t)}_j\bigr)_{1\le j \le K}\bigr) \bigr] &\sim& C_R t \eee\bigl[ 1_{\{K<
\tau
_t^+<\tau_0^-\}} F \bigl(T_t^+,\bigl(\widetilde S^{(t)}_j\bigr)_{1\le j
\le K}\bigr) \bigr]\\
&=:& A_t.
\end{eqnarray*}
Let us now introduce $\ell_t:=t-2t^{\gamma}$ with
$(1+\delta/2)^{-1}<\gamma<1$ and observe that
$\tau_{\ell_t}^+<\tau_0^-$ on the event $\{\tau_t^+<\tau_0^-\}$.
Recalling that part (ii) of Lemma \ref{lemRW} says that $(T^+_t, t\ge0)$
is bounded in $L^p$ for all $1<p < 1+\delta$, we get $\ppp(
T^+_{\ell_t} > t^\gamma) \le c t^{-\gamma p}=o(t^{-1})$ by
choosing $p$ such that $\gamma p>1$. Therefore we obtain
\begin{eqnarray*}\label{eqcritic1}
A_t &=& C_R t \eee\bigl[ 1_{\{K< \tau_t^+<\tau_0^-\}} 1_{\{
S_{\tau
^+_{\ell_t}} \le t- t^\gamma\}} F\bigl(T_t^+,\bigl(\widetilde
S^{(t)}_j\bigr)_{1\le j \le K}\bigr) \bigr] + o_t(1)
\\
&=& A'_t + A''_t + o_t(1),
\end{eqnarray*}
where $o_t(1) \to0$ as $t \to\infty$ and
\begin{eqnarray*}
A'_t &:=& C_R t \eee\bigl[ 1_{\{ \tau_t^+<\tau_0^-\}} 1_{\{
S_{\tau
^+_{\ell_t}} \le t- t^\gamma\}} 1_{\{
\tau^+_t - \tau^+_{\ell_t}> K \}} F \bigl(T_t^+,\bigl(\widetilde
S^{(t)}_j\bigr)_{1\le j \le K}\bigr) \bigr],
\\
A''_t &:=& C_R t \eee\bigl[ 1_{\{K< \tau_t^+<\tau_0^-\}}
1_{\{S_{\tau^+_{\ell_t}} \le t- t^\gamma\}} 1_{\{
\tau^+_t - \tau^+_{\ell_t} \le K \}} F \bigl(T_t^+,\bigl(\widetilde
S^{(t)}_j\bigr)_{1\le j \le K}\bigr) \bigr].
\end{eqnarray*}
Applying the strong Markov property at the stopping time
$\tau_{\ell_t}^+$ yields
\[
A'_t= C_R t \eee\bigl[ 1_{\{\tau_{\ell_t}^+<\tau_0^-\}}
1_{\{S_{\tau^+_{\ell_t}} \le t- t^\gamma\}}
f(S_{\tau^+_{\ell_t}}) \bigr],
\]
where
%
\begin{equation}
\label{eqcritic5} f(x):= \eee_{x}\bigl[ 1_{\{K<
\tau_{t}^+<\tau_0^-\}} F \bigl(T_t^+,\bigl(\widetilde S^{(t)}_j\bigr)_{1\le j \le K}\bigr)
\bigr].
\end{equation}
Then, writing
\[
\eee_{x}\bigl[ 1_{\{K< \tau_{t}^+\}} F \bigl(T_t^+,\bigl(\widetilde
S^{(t)}_j\bigr)_{1\le j \le K}\bigr)\bigr]=\eee \bigl[ 1_{\{K<
\tau_{L}^+\}} F\bigl(T_L^+,\bigl(\widetilde S^{(L)}_j\bigr)_{1\le j \le K}\bigr)\bigr],
\]
with $L=t-x$, equation (\ref{eqci-conseq-part-i}) yields
%
\begin{equation}
\label{eqresuadmit}\qquad  \max_{x \in[\ell_t; t- t^\gamma]}
\bigl\vert\eee_{x}\bigl[ 1_{\{K< \tau_{t}^+\}}
F\bigl(T_t^+,\bigl(\widetilde S^{(t)}_j\bigr)_{0\le j \le K}\bigr)\bigr] - C_{ F}
\bigr\vert\longrightarrow0,\qquad t \to\infty,
\end{equation}
from which we deduce
\[
\max_{x \in[\ell_t; t- t^\gamma]} \bigl\vert f(x) -
C_{ F} \bigr\vert\longrightarrow0, \qquad t \to\infty,
\]
since uniformly in $x \ge\ell_t$, $\ppp_x(\tau_0^-<\tau_{t}^+)=
\ppp(
\tau^-_{-x} < \tau^+_{t-x})\le\ppp( \tau^-_{- \ell_t} < \tau
^+_{t^\gamma})= o_t(1)$. Furthermore, observing that $\ppp(\tau
_{\ell
_t}^+<\tau_0^-)
\sim{1\over C_R t}$ [see part (v) of Lemma~\ref{lemRW} and recall that
$\ell_t=t-2t^{\gamma}$ with $\gamma<1$] and $\ppp(t-
S_{\tau^+_{\ell_t}} \le t^\gamma)= \ppp(T^+_{\ell_t}> t^\gamma
)=o(t^{-1})$ imply
$\ppp(\tau_{\ell_t}^+<\tau_0^-; S_{\tau^+_{\ell_t}} \le t-
t^\gamma) \sim1/C_R t$, and when $t$ tends to infinity, we obtain
%
\begin{equation}
\label{eqcritic2} A'_t \longrightarrow C_{ F},\qquad t \to
\infty.
\end{equation}
Similarly, the strong Markov property applied at the stopping time
$\tau_{\ell_t}^+$ implies
\[
A''_t \le C_R t \eee\bigl[ 1_{\{\tau_{\ell_t}^+<\tau_0^-\}}
1_{\{S_{\tau^+_{\ell_t}} \le t- t^\gamma\}}
\ppp_{S_{\tau^+_{\ell_t}}} \bigl(\tau^+_{t} \le K\bigr) \bigr].
\]
Moreover, observe that
%
\begin{equation}
\label{eqsouscrititc1} \sup_{x \le t- t^\gamma} \ppp_{x} \bigl(\tau^+_{t}
\le K\bigr) \le\ppp_{t- t^\gamma} \bigl(\tau^+_{t} \le K\bigr)= \ppp\bigl(\tau
^+_{t^\gamma}
\le K\bigr)=o_t(1),
\end{equation}
which implies
%
\begin{equation}
A''_t \le C_R t \ppp\bigl(\tau_{\ell_t}^+<\tau_0^-,
S_{\tau^+_{\ell_t}} \le t- t^\gamma\bigr) \ppp\bigl(\tau^+_{t^\gamma}
\le
K\bigr)=o_t(1), \label{eqcritic3}
\end{equation}
by recalling that $\ppp(\tau_{\ell_t}^+<\tau_0^-;
S_{\tau^+_{\ell_t}} \le t- t^\gamma) \sim{1\over C_R t}$. Combining
(\ref{eqcritic2}), (\ref{eqcritic3}) and recalling
(\ref{eqcritic1}), we obtain $A_t \to C_{F}$, when $t \to
\infty$, which concludes the proof of part (ii) in the case
$\eee[S_1]=0$.

The case $\eee[S_1]>0$ is similar but easier. Indeed, combining parts
(iii) and (iv) of Lemma \ref{lemRW} implies
\begin{eqnarray*}
\eee^{+}\bigl [ 1_{\{K< \tau_t^+\}} F\bigl(T_t^+,\bigl(\widetilde
S^{(t)}_j\bigr)_{1\le j \le K}\bigr) \bigr]& \sim& C_R \eee\bigl[ 1_{\{K<
\tau
_t^+<\tau_0^-\}} F\bigl(T_t^+,\bigl(\widetilde S^{(t)}_j\bigr)_{1\le j \le K}\bigr) \bigr]\hspace*{-15pt}
\\
&=:& \widetilde A_t.\hspace*{-15pt}
\end{eqnarray*}
Recalling that $\ell_t=t-2t^{\gamma}$ and that part (ii) of Lemma
\ref{lemRW} implies $\ppp( T^+_{\ell_t} > t^\gamma)=o_t(1)$, and
we get
%
\begin{eqnarray}
\widetilde A_t&=& C_R \eee\bigl[ 1_{\{K< \tau_t^+<\tau_0^-\}}
1_{\{
S_{\tau^+_{\ell_t}} \le t- t^\gamma\}} F\bigl(T_t^+,\bigl(\widetilde
S^{(t)}_j\bigr)_{1\le j \le K}\bigr) \bigr]+o_t(1)
\nonumber\hspace*{-35pt}
\\[-8pt]
\\[-8pt]
\nonumber
&=& C_R \eee\bigl[ 1_{\{ \tau_t^+<\tau_0^-\}} 1_{\{S_{\tau
^+_{\ell
_t}} \le t- t^\gamma\}} 1_{\{ \tau^+_t - \tau^+_{\ell_t} > K \}}
F\bigl(T_t^+,\bigl(\widetilde S^{(t)}_j\bigr)_{1\le j \le K}\bigr) \bigr]+o_t(1),\hspace*{-35pt}
\end{eqnarray}
the last equality being a consequence of (\ref{eqsouscrititc1}),
which still holds in the case $\eee[S_1]>0$. Then, the strong Markov
property yields
%
\begin{equation}
\widetilde A_t= C_R \eee\bigl[ 1_{\{\tau_{\ell_t}^+<\tau_0^-\}}
1_{\{S_{\tau^+_{\ell_t}} \le t- t^\gamma\}}
f(S_{\tau^+_{\ell_t}}) \bigr]+o_t(1),
\end{equation}
where we recall that the function $f$ is defined by
(\ref{eqcritic5}). Now the strategy is exactly the same as for the
previous case. Indeed, since $\ppp_x(\tau_0^-<\tau_{t}^+)= o_t(1)$
(uniformly in $x \ge\ell_t$) is still true,
(\ref{eqci-conseq-part-i}) implies $\max_{x \in[\ell_t; t-
t^\gamma]} \vert f(x) - C_{F} \vert\to0$, when $t$ tends to
$\infty$. Combining this with part (v) of Lemma \ref{lemRW} [which
implies $\ppp(\tau_{\ell_t}^+<\tau_0^-; S_{\tau^+_{\ell
_t}} \le
t- t^\gamma) \to1/C_R$] yields $\widetilde A_t \to C_{ F}$, when $t
\to\infty$. This completes the proof of part~(ii) of the lemma and
completes the proof of Lemma
\ref{lemlaw-convergence-reverse}. 

\begin{pf*}{Proof of Lemma \ref{lpositivewalk}} We may assume that
$p$ equals some integer, say, $m\ge1$. Indeed, for any $m-1< p\le m$,
by the concavity,
\[
\eee_x\Biggl[ \sum_{k=0}^{\tau_t^+ -1} a_{k+1} e^{\kappa(S_{k}-t) }
\Biggr] ^p \le\eee_x\Biggl[ \sum_{k=0}^{\tau_t^+ -1}
(a_{k+1})^{p/m}
e^{\kappa p (S_{k}-t)/m } \Biggr] ^m.
\]
Applying (\ref{zb2new}) to $ ((a_{k+1})^{p/m}, S_k- S_{k-1})$
with integer $m$ yields the general case~$p$.

Now, we consider $p=m$ is some integer and prove (\ref{zb2new}). First,
\begin{eqnarray*} \eee\Biggl[ \sum_{k=0}^{\tau_t^+ -1} e^{\kappa
(S_{k}-t) } \Biggr] \le \sum_{k=0}^\infty\eee\bigl[ 1_{\{
\overline
S_k \le t \}} e^{\kappa(\overline S_k -t) } \bigr] = \int_0^t e^{
-\kappa(t-y)} \,{d}u(y),
\end{eqnarray*}
where $\overline S_k:= \max\{S_j\dvtx0\le j \le k \} $ and
\[
u(y):= \sum_{n=0}^\infty\ppp( \overline S_n \le y),\qquad y
\ge0.
\]

Remark that $ u$ is finite and satisfies the following
renewal equation (see Heyde~\cite{heyde}, Theorem~1):
\[
u(y)=1_{\{0\le y\}} + F \ast u (y),\qquad y\ge0,
\]
with $F(s):=\ppp( S_1 \le s), s\in{\mathbb R}. $ According to the
renewal theorem (see Heyde \cite{heyde}, Theorem 2 or Feller \cite
{feller} page 362 (1.17) and page 381), $ \int_0^t e^{ - \kappa(t-y)}
\,{d}u(y)= O(1)$ as $t \to\infty$ (the limit exists in the
nonarithmetic case). By linear transformation, we obtain that for any
$\kappa>0$, $ \eee_x[ \sum_{k=0}^{\tau_t^+ -1} e^{\kappa(S_{k}-t)
} ] $ is uniformly bounded for all $x \le t$.

We now prove the lemma by induction on $m$. By independence,\break $ \eee
_x[ \sum_{k=0}^{\tau_t^+ -1} a_{k+1} e^{\kappa(S_{k}-t) }
] = \sum_{k\ge0} \eee_x[ e^{\kappa(S_{k}-t) }, k < \tau_t^+ -1
] \eee[ a_1 ]$ is bounded by some constant (the law of $a_{k+1}$
does not depend on $x$). This proves the lemma in the case $m=1$.

Let $m\ge2$ and assume that the lemma holds for $1, \ldots, m-1$. Write
$\chi_i:= \sum_{k=i}^{\tau_t^+ -1} a_{k+1} e^{\kappa(S_{k}-t) } $
for $0\le i < \tau_t^+$ and $\chi_{\tau_t^+}:=0$. Note that
\begin{eqnarray*}
(\chi_0)^m &=& \sum_{i=0}^{\tau_{t}^{+} -1} \bigl[ (\chi_i)^m -
(\chi
_{i+1})^m\bigr]\\
&=& \sum_{j=0}^{m-1} \pmatrix{m\cr j} \sum
_{i=0}^{\tau
_t^+ -1} a_{i+1}^{m-j} e^{ (m-j) \kappa(S_i -t)} (\chi_{i+1})^{
j}.
\end{eqnarray*}

Applying Markov's property at $i+1$, we get
\begin{eqnarray*} \eee_x\bigl[\chi_0^m \bigr] & = &\sum_{j=0}^{m-1}
\pmatrix{m\cr j} \eee_x \Biggl[ \sum_{i=0}^{\tau_t^+ -1}
a_{i+1}^{m-j}
e^{ (m-j) \kappa(S_i -t)} \eee_{S_i} \bigl[ (\chi_{i+1})^{ j}\bigr]
\Biggr] \\
&\le& c \sum_{j=0}^{m-1} \eee_x \Biggl[ \sum_{i=0}^{\tau_t^+ -1}
a_{i+1}^{m-j} e^{ (m-j) \kappa(S_i -t)} \Biggr],
\end{eqnarray*}
since by the induction hypothesis $ \eee_{S_i} [ (\chi
_{i+1})^{j}] $ is bounded by some constant. The last expectation is
again uniformly bounded (the case $m=1$ of the lemma), which proves
that the lemma holds for $m$, as desired.
\end{pf*}

\subsection{\texorpdfstring{Proof of Lemma \protect\ref{lemestcrit}}{Proof of Lemma 6}}

For $a \in\r$, denote as before by $T_a^+:=S_{\tau_a^+}-a>0$ (resp.,
$T^-_a:= a- S_{\tau_a^-}>0$)
the overshoot (resp., undershoot) at level $a$. Clearly the overshoot
$T_a^+$ is also the overshoot at the level $a$ for the strict
ascending ladder heights
$(H_n)$. By assumption (\ref{Hyp-expmoment}), $\max(S_1, 0)$ has finite
$\eta$-exponential moment. This in view of Doney \cite{D80} implies that
$ \eee[ e^{ \delta H_1}] < \infty$ for any $0< \delta< \eta$.
Applying Chang (\cite{C94}, Proposition 4.2) yields \eqref{ta+}.
Similarly for the undershoot $T_a^->0$: since $\max(-S_1, 0)$ has a
finite $(1+\eta)$-exponential moment, again \eqref{ta-} follows from
Chang (\cite{C94}, Proposition 4.2).

By (\ref{ta+}) and (\ref{ta-}),
$ \max_{0\le k \le\tau_{0}^- \wedge\tau_{L }^+} | S_k|\le L + T_{L
}^+ + T_{0}^- $ is integrable under $\ppp_a$. By applying the optional
stopping theorem, we get\vspace*{-1pt}
\[
a= \eee_a [ S_{ \tau_{0}^- \wedge
\tau_{L }^+} ] = \eee_a \bigl[ ( S_{ \tau_{0}^-} - S_{ \tau_{L
}^+}) 1_{\{ \tau_{0}^-< \tau_{L }^+\}}\bigr] + \eee_a[ S_{ \tau_{L
}^+}].
\]
Observe that $\eee_a[ S_{ \tau_{L }^+}]= L + \eee_a[
T^+_{L}] \le L + c $ by (\ref{ta+}). Since $S_{
\tau_{0}^-} - S_{ \tau_{L }^+} < -L$, this implies \eqref{pal1}.
Exactly doing the same and using (\ref{ta-}), we get \eqref{pal2}.

Let us mention that by considering the martingale
$(S_j^2 - \operatorname{Var}(S_1) j)_{j\ge1}$, which is uniformly integrable on
$[0, \tau_0^- \wedge\tau_L^+]$, we can find some constant $c>0$
such that for all $L>1$ and $0\le a \le L$,\vspace*{-1pt}
%
\begin{equation}\label{c4ea} \eee_a\bigl[ \tau_0^- \wedge\tau_L^+
\bigr] \le c L^2.
\end{equation}


{(i) Proof of (\ref{pal})}. If $L-a \ge{L\over3}$, we deduce from
(\ref{ta-}) that\break $\eee_a [ e^{- S_{\tau_0^-}} 1_{\{\tau_0^- <
\tau_L^+\}}] \le\eee_a [ e^{- S_{\tau_0^-}} ] \le
c$ which is less than $c' {L-a +1\over L}$ if
$c'\ge3 c$.

Let $0<L-a < {L\over3}$. Note that under $\ppp_a$, $\tau^-_0< \tau
^+_L$ implies that $ \tau^-_{L/2} \le\tau^-_0< \tau^+_L$. Then by the
strong Markov property at $\tau^-_{L/2}$,\vspace*{-1pt}
\begin{eqnarray*}
\eee_a \bigl[ e^{- S_{\tau_0^-}} 1_{\{\tau_0^- < \tau_L^+\}}
\bigr] &= &
\eee_a\bigl [ e^{- S_{\tau_0^-}} 1_{\{\tau_{L/2}^- \le\tau_{0}^- <
\tau_L^+\}} \bigr]
\\
&= & \eee_a \bigl[ 1_{\{\tau_{L/2}^- < \tau_L^+\}} \eee_{S_{\tau
_{L/2}^-}} \bigl[ e^{- S_{\tau_0^-}} 1_{\{\tau_{0}^- < \tau_L^+\}}
\bigr] \bigr]
\\
&\le& \eee_a \bigl[ 1_{\{\tau_{L/2}^- < \tau_L^+\}} \bigl(c+
e^{-S_{\tau_{L/2}^-}} 1_{\{S_{\tau_{L/2}^-}<0\}} \bigr) \bigr],
\end{eqnarray*}
where we use the fact that for all $z:=
S_{\tau_{L/2}^- }\ge0$, $\eee_z [ e^{- S_{\tau_0^-}}
1_{\{\tau_{0}^- < \tau_L^+\}} ] \le\eee_z [ e^{- S_{\tau_0^-}}]
\le c$ by (\ref{ta-}). Since
$S_{\tau_{L/2}^-}<0$ means that $T^-_{L/2} \ge L/2$, we deduce
from (\ref{ta-}) that
\[
\eee_a \bigl[ e^{-S_{\tau_{L/2}^-}} 1_{\{
S_{\tau_{L/2}^-}<0\}} \bigr] = \eee_a \bigl[ e^{ {L/2} +
T^-_{L/2}} 1_{\{T^{-}_{L/2} \ge L/2\}}\bigr] \le c
e^{- \delta L/2}.
\]

This together with (\ref{pal1}) gives that
\begin{eqnarray*}
\eee_a \bigl[
e^{- S_{\tau_0^-}} 1_{\{\tau_0^- < \tau_L^+\}}\bigr] &\le& c
\ppp_a\bigl( \tau_{L/2}^- < \tau_L^+\bigr)+ c e^{-
\delta L/2} \\
&=& c
\ppp_{a-L/2}\bigl( \tau_{0}^- < \tau_{L/2}^+\bigr)+ c e^{-
\delta L/2} \\&\le& c { L-a + c'\over(L/2)} + c e^{- \delta L/2}
\\
&\le& c''{ L-a + 1 \over L}.
\end{eqnarray*}

{(ii) Proof of (\ref{foot3})}. Let us show that $\eee[\sum_{j
=0}^{\tau_0^- -1} e^{- \delta S_j}]<\infty$,
\[
\eee\Biggl[ \sum_{j
=0}^{\tau_0^- -1} e^{- \delta S_j}\Biggr]
=
\sum_{j\ge0} \eee\bigl[e^{- \delta S_j}, j<\tau_0^-\bigr]
\le
\sum_{j\ge0} c (1+j)^{-3/2} <\infty,
\]
where we used Theorem 4 (and Theorem 6 if $S_1$ is lattice)
of \cite{vatutin-wachtel} for the bound of $\eee[e^{- \delta
S_j}, j<\tau_0^-]$. Let $(H_n^-,\sigma_n^-)_{n\ge0}$ be the
strict ascending ladder heights and epochs of $-S$ (with $\sigma
_0^-:=0$). For $a> 0$, we notice that
\begin{eqnarray*}
\eee_a\Biggl [ \sum_{j
=0}^{\tau_0^- -1} e^{- \delta S_j}\Biggr]
&=& \eee\Biggl[ \sum_{j
=0}^{\tau_{-a}^- -1} e^{- \delta(a+S_j)}\Biggr] \\
&=&\sum_{n= 0}^\infty\eee\biggl[\sum_{\sigma_n^-\le j <\sigma_{n+1}^-}
e^{-\delta(a+S_j)} 1_{\{H_n^-\le a\}} \biggr]\\
&=&
\sum_{n= 0}^{\infty} \eee\bigl[ e^{-\delta(a-H_n^-)} 1_{\{
H_n^-\le
a\}} \bigr] \eee\Biggl[ \sum_{j
=0}^{\tau_0^- -1} e^{- \delta S_j}\Biggr],
\end{eqnarray*}
by applying the strong Markov property at $\sigma_n^-$. We
showed that\break $\eee[ \sum_{j
=0}^{\tau_0^- -1} e^{- \delta S_j}] <\infty$. On the other hand,
Lemma \ref{lpositivewalk} applied to the random walk $(H_n^-)_{n\ge0}$
says that
\[
\sup_{a>0} \sum_{n= 0}^{\infty} \eee\bigl[ e^{-\delta(a-H_n^-)}
1_{\{
H_n^-\le a\}} \bigr] < \infty.
\]
Hence $\sup_{a\ge0} \eee_a [ \sum_{j
=0}^{\tau_0^- -1} e^{- \delta S_j}]< \infty$. Similary, by
considering the random walk $L-S_\cdot$, we get that $\eee_a [
\sum_{j =0}^{\tau_L^+-1} e^{- \delta(L-S_j)}] $ is uniformly
bounded by some constant. This proves (\ref{foot3}).

{(iii)} Proof of (\ref{diftau0}). Considering the value of the time
$\tau_0^-$, then using Markov's property, we have
\begin{eqnarray*}
\eee_a\bigl[e^{S_{\tau_0^- - 1} - S_{\tau_0^-}} \bigr]
&=&
\sum_{k\ge1} \eee_a\bigl[e^{S_{k-1} - S_k}1_{\{\tau_0^-=k\}}\bigr]\\
&=&
\sum_{k\ge1} \eee_a\bigl[h(-S_{k-1})1_{\{\tau_0^-\ge k\}}\bigr]
\end{eqnarray*}
where for any $y\in\r$, $h(y):= \eee[e^{-S_1}1_{\{S_1 \le
y\}
}]\le e^{\delta y} \eee[e^{-(1+\delta) S_1}] =c e^{\delta y}$ for
$\delta>0$ small enough. Hence,
\[
\eee_a\bigl[e^{S_{\tau_0^- - 1} - S_{\tau_0^-}} \bigr]
\le
c \eee_a\Biggl[\sum_{k=0}^{\tau_0^- - 1} e^{-\delta S_{k}}\Biggr]
\]
and (\ref{diftau0}) follows from (\ref{foot3}).

{(iv) Proof of (\ref{L-a+1}) and (\ref{a+1})}: Clearly
(\ref{a+1}) follows from (\ref{L-a+1}) by considering the random
walk $(L-S_j)_{j\ge0}$. It suffices to prove (\ref{L-a+1}). If $L-a
\ge L/3$, there is nothing to prove since $\eee_a [ \sum_{0\le j <
\tau_0^- \wedge\tau_L^+} e^{- \delta S_j} ] \le\eee_a [
\sum_{0\le j < \tau_0^- } e^{- \delta S_j} ] $ is less than
some constant by (\ref{foot3}).

Considering $L-a < L/3$, we have
\begin{eqnarray*}
&& \eee_a \biggl[
\sum_{0\le j < \tau_0^- \wedge\tau_L^+} e^{- \delta S_j} \biggr]
\\
&&\qquad=
\eee_a \biggl[1_{\{ \tau_{L/2}^- \ge\tau_0^- \wedge\tau_L^+\}}
\sum_{0\le j < \tau_0^- \wedge\tau_L^+} e^{- \delta S_j} \biggr]\\
&&\qquad\quad{} +
\eee_a \biggl[1_{\{ \tau_{L/2}^- < \tau_0^- \wedge\tau_L^+\}} \sum
_{0\le
j < \tau_0^- \wedge\tau_L^+} e^{- \delta S_j} \biggr]
\\
&&\qquad\le \eee_a\bigl[ e^{-\delta L/2} \bigl(\tau_0^- \wedge\tau
_L^+\bigr)\bigr] +
\eee_a \biggl[1_{\{ \tau_{L/2}^- < \tau_0^- \wedge\tau_L^+\}}
\sum_{\tau_{L/2}^- \le j < \tau_0^- \wedge\tau_L^+} e^{- \delta S_j}
\biggr]
\\
&&\qquad\le c L^2 e^{-\delta L/2} + \eee_a\biggl[1_{\{ \tau_{L/2}^- < \tau_0^-
\wedge\tau_L^+\}}
\eee_{ S_{\tau_{L/2}^-} } \biggl[ \sum_{0 \le j < \tau_0^- \wedge
\tau_L^+} e^{- \delta S_j}\biggr] \biggr],
\end{eqnarray*}
by using (\ref{c4ea}) and the strong Markov property at
$\tau_{L/2}^-$. Let $x:=S_{\tau_{L/2}^-} < L/2$. If $x<0$, then
under $\ppp_x$, $\tau_0^-=0$ and $\e_{x } [\sum_{0 \le j <
\tau_0^-
\wedge\tau_L^+} e^{- \delta S_j}] =0$, whereas if $0\le x <
L/2$, $\eee_{x } [\sum_{0 \le j < \tau_0^- \wedge\tau_L^+} e^{-
\delta S_j}] \le c $ by (\ref{foot3}). Then we get
\begin{eqnarray*} \eee_a \biggl[ \sum_{0\le j < \tau_0^- \wedge\tau_L^+}
e^{- \delta
S_j} \biggr] &\le& c L^2 e^{-\delta L/2} + c \ppp_a\bigl(
\tau_{L/2}^- < \tau_0^- \wedge\tau_L^+\bigr) \\
&\le& c L^2 e^{-\delta L/2} + c \ppp_a\bigl(
\tau_{L/2}^- < \tau_L^+\bigr) \\
&\le& c L^2 e^{-\delta L/2} + c { L - a
+c'\over L/2},
\end{eqnarray*}
by using (\ref{pal1}). This proves (\ref{L-a+1}).

{(v) Proof of (\ref{eal})}: By monotonicity, it is sufficient to
prove $ (\ref{eal})$ for \mbox{$0<\delta<\eta$}. Then notice that
\begin{eqnarray*}
&&\eee_a \biggl[ e^{- S_{\tau_0^-}} 1_{\{ \tau_0^- < \tau_L^+\}}
\sum_{
0\le j < \tau_0^-} e^{ - \delta(L-S_j)}\biggr]\\
&&\qquad= \sum_{n=1}^\infty\eee_a \biggl[ 1_{\{ n \le\tau_L^+ \wedge\tau
_0^-, S_n <0\}} e^{-S_n } \sum_{ 0\le j < n} e^{ - \delta
(L-S_j) }\biggr].
\end{eqnarray*}

Applying Markov's property of $S$ at $n-1$ and using the
fact that for all $x \ge0$, $
\eee_x[e^{-S_1 } 1_{\{S_1 <0 \}}] = \eee[e^{ -x- S_1} 1_{\{ S_1<-x\}}]
\le
c(\delta) e^{-(1+\delta) x }$ by (\ref{Hyp-expmoment}) (recall that
$0<\delta<\eta$), we get that
%
\begin{eqnarray}\label{iv00}
&&
\eee_a \biggl[ e^{- S_{\tau_0^-}} 1_{\{ \tau_0^- < \tau_L^+\}}
\sum_{
0\le j < \tau_0^-} e^{ - \delta(L-S_j)}\biggr] \nonumber
\\
&&\qquad \le c \sum_{n=1}^\infty\eee_a \biggl[ 1_{\{ n \le\tau_L^+
\wedge\tau_0^- \}} e^{- (1+\delta) S_{n-1}} \sum_{ 0\le j < n} e^{
-\delta(L-S_j) }\biggr]
\\
&&\qquad= c \sum_{j=0}^\infty\eee_a \biggl[ 1_{\{ j<\tau_L^+
\wedge\tau_0^-\}} e^{- \delta(L-S_j)} \eee_{S_j}\biggl[
\sum_{0\le m < \tau_L^+ \wedge\tau_0^- } e^{- (1+\delta)
S_m}\biggr]\biggr]
,\nonumber
\end{eqnarray}
where the last equality follows from Markov's
property at $j$. Applying~(\ref{L-a+1}) and~(\ref{a+1}), we get that
\begin{eqnarray*}
&& \eee_a\biggl [ e^{- S_{\tau_0^-}} 1_{\{ \tau_0^- < \tau_L^+\}}
\sum_{ 0\le j < \tau_0^-} e^{ - \delta(L-S_j)}\biggr]
\\
&&\qquad\le c \sum_{j=0}^\infty\eee_a \biggl[ 1_{\{ j<\tau_L^+ \wedge
\tau
_0^-\}} e^{- \delta(L-S_j)} c { L- S_j +1 \over L}\biggr]
\\
&&\qquad\le {c' \over L} \eee_a\biggl [ \sum_{0\le j<\tau_L^+ \wedge\tau_0^-
} e^{- {(\delta/2) } (L-S_j)} \biggr]
\le c { a+1\over L^2},
\end{eqnarray*}
proving (\ref{eal}).

We mention that (\ref{iv00}) also holds with $\delta=0$, which
implies that
%
\begin{eqnarray}
\label{iv01} \eee_a\bigl [ e^{- S_{\tau_0^-}} 1_{\{ \tau_0^- <
\tau_L^+\}} \tau_0^- \bigr] \le c \eee_a\bigl[\tau_0^- \wedge\tau_L^+\bigr]
\le c' L^2
\nonumber
\\[-8pt]
\\[-8pt]
 \eqntext{\forall L\ge1, 0\le a\le L.}
\end{eqnarray}

\subsection{\texorpdfstring{Proof of Lemmas \protect\ref{lemestsub} and \protect\ref{Laux1}}
{Proof of Lemmas 7 and 8}}

Keeping the notation $T_{a}^-$ for the undershoot at level $a$, we have
as before for any $0<r <\eta_1$,
%
\begin{equation}\label{ta-sub}
\ppp_b\bigl(T_a^- >x\bigr) \le c(r )e^{-r x} \qquad\forall a\le b,
\forall x>0.
\end{equation}

\begin{pf*}{Proof of Lemma \ref{lemestsub}}
{(i) Proof of (\ref{est1sub})}. It is a straightforward consequence
of (\ref{ta-sub}).

{(ii) Proof of (\ref{est3sub})}. Let us introduce the tilted
measure $\tilde\ppp_a$ defined by ${d \tilde\ppp_a \over d \ppp_a}
|_{\sigma(S_0,\ldots,S_n)}:= e^{\gamma(S_n-S_0)}$. Under
$\tilde
\ppp_a$, the random walk drifts to $+\infty$. We write
\begin{eqnarray*}
&&\eee_a\biggl[ \sum_{0\le\ell< \tau_L^+}(1+L-S_{\ell})^{\alpha
}e^{r
S_{\ell}} \biggr]
\\
&&\qquad=
\sum_{\ell\ge0} \eee_a\bigl[(1+L-S_\ell)^\alpha e^{r S_\ell} 1_{\{
\ell
<\tau_L^+\}} \bigr]\\
&&\qquad=
e^{\gamma a}\sum_{\ell\ge0} \tilde\eee_a\bigl[(1+L-S_\ell)^\alpha
e^{(r
-\gamma) S_\ell} 1_{\{\ell<\tau_L^+\}}\bigr]\\
&&\qquad=
e^{\gamma a}e^{(r -\gamma) L} \tilde\eee_a\biggl[ \sum_{0\le\ell<
\tau_L^+} (1+L-S_\ell)^\alpha e^{(r -\gamma) (S_{\ell}-L)}\biggr]\\
&&\qquad\le
c e^{\gamma a}e^{(r -\gamma) L} \tilde\eee_a \biggl[ \sum_{0\le
\ell<
\tau_L^+} e^{(r -\gamma) (S_{\ell}-L)/2}\biggr].
\end{eqnarray*}
Therefore, we only have to show that
\[
\sup_{a\ge0} \tilde\eee_a\biggl[ \sum_{0\le\ell< \tau_L^+}
e^{(r
-\gamma) (S_{\ell}-L)/2}\biggr] \le c,
\]
which is done by the same argument as in the proof of (\ref{foot3}).

{(iii) Proof of (\ref{est7sub})}. We have
\begin{eqnarray*}
\eee_a\Biggl[ \sum_{\ell=0}^{\min(\tau_0^-,\tau_L^+)} (1+L- S_\ell
)^\alpha e^{\gamma S_{\ell}} \Biggr]
&=&
e^{\gamma a} \tilde\eee_a\Biggl[\sum_{\ell=0}^{\min(\tau_0^-,\tau_L^+)}
(1+L- S_\ell)^\alpha\Biggr] \\
&=&
e^{\gamma a} \tilde\eee\Biggl[\sum_{\ell=0}^{\min(\tau_{-a}^-,\tau
_{L-a}^+)} (1+L-a- S_\ell)^\alpha\Biggr].
\end{eqnarray*}

Remark that $(1+L-a- S_\ell)^\alpha\le c (1+L-a )^\alpha+ c \vert
S_\ell\vert^\alpha1_{\{ S_\ell<0\}}$ and that $\tilde\eee
[\sum
_{\ell\ge0}\vert S_\ell\vert^\alpha1_{\{ S_\ell<0\}}] <
\infty$
(indeed observe that for any $\gamma' \in(0,\gamma)$ there exists
$c(\alpha, \gamma')$ such that $\sum_{\ell\ge0}\vert S_\ell\vert
^\alpha1_{\{ S_\ell<0\}} \le c(\alpha, \gamma') \sum_{\ell\ge0} e^{-
\gamma' S_{\ell}}$, whose expectation under $\tilde\ppp$ is finite;
see Kesten \cite{K73}). Therefore, we get
\begin{eqnarray*}
\tilde\eee\Biggl[\sum_{\ell=0}^{\min(\tau_{-a}^-,\tau_{L-a}^+)}
(1+L-a- S_\ell)^\alpha\Biggr] &\le& c' (1+L-a)^\alpha \tilde\eee
\bigl[ \tau_{L-a}^+ \bigr] + c'\\
& \le& c (1+L-a)^{ \alpha+1},
\end{eqnarray*}
which completes the proof of the lemma.
\end{pf*}

\begin{pf*}{Proof of Lemma \ref{Laux1}} First, we remark that it
is enough to prove the lemma for integer $p$. In fact let us assume
that (i) holds for any integer $p$ satisfying the hypothesis in (i).
Now for $0 \le p < {\gamma\over b}$, we choose an arbitrary integer
$k$ larger than~$p$. Then \eqref{est8sub} holds for any $(\tilde
a_\cdot
, \tilde b)$ [in lieu of $(a_\cdot, b)$] satisfying the hypothesis in
(i): $0 \le k < \gamma/ \tilde b $ and $ \eee[ (1+ 1_{\{S_1< 0\}}
e^{- \eta S_1}) \tilde a_1^k ] < \infty$. Observe that $\tilde a_\ell:=
(a_\ell)^{p/k}$ for any $\ell\ge1$ and $\tilde b:= {pb\over k}$
fulfill the above hypothesis, hence
\[
\eee_x \Biggl[ e^{ - \eta S_{ \tau^-_0}} \Biggl( \sum_{\ell=1}^{
\tau
^-_0} e^{ {(p b/ k)} S_{\ell-1}} (a_\ell)^{p/k}\Biggr)^k \Biggr]
\le
c_k e^{ \tilde b k x}= c_k e^{ b p x}\qquad \forall x \ge0.
\]
Since $k \ge p$, we have by concavity that
\[
\eee_x \Biggl[ e^{ - \eta S_{ \tau^-_0}} \Biggl( \sum_{\ell=1}^{
\tau
^-_0} e^{ b S_{\ell-1}} a_\ell\Biggr)^p \Biggr] \le\eee_x\Biggl [
e^{ -
\eta S_{ \tau^-_0}} \Biggl( \sum_{\ell=1}^{ \tau^-_0} e^{ {(p b/ k)}
S_{\ell-1}} (a_\ell)^{p/k}\Biggr)^k \Biggr] \le c_k e^{ b p x}.
\]
Hence it is enough to show (i) with an integer $p$. The same
is true for (ii).

Now we assume $p$ is an integer, and we shall use Markov's property to
expand the power. Let either $\chi:=\tau_0^-$ or $\chi:= \min(\tau_0^-
, \tau^+_L)$ and consider a measurable function $f\dvtx\r\to\r_+$.
Define
\[
A_{\chi, f}(x, k ):= \eee_x \Biggl[ e^{ - \eta S_{ \tau^-_0}} \Biggl(
\sum
_{\ell=1}^{ \chi}f( S_{\ell-1}) a_\ell\Biggr)^k \Biggr],\qquad
k\ge0, x \in\r,
\]
and we mention that $A_{\chi, f}(x,0)= e^{ -\eta x}$ if
$x<0$, $A_{\chi, f}(x,k)=0$ if $x<0$ and $k\ge1$. Let $k\ge1$ and
$Y_i:= \sum_{\ell=i}^{ \chi}f( S_{\ell-1}) a_\ell$ for $1\le i
\le\tau_0^-$, $Y_{\chi+1}:=0$. Then
\[
Y_1^k= \sum_{i=1}^{\chi} \bigl(Y_i^k - Y_{i+1}^k\bigr)= \sum_{r=1}^k C_k^r
\sum
_{i=1}^{\chi} \bigl(f(S_{i-1})\bigr)^r (a_i)^r (Y_{i+1})^{k-r}.
\]

Applying Markov's property at $i$ gives that
%
\begin{eqnarray} \label{axkxk}
 A_{\chi, f}(x,k)&= &\sum_{r=1}^k C_k^r \sum
_{i=1}^{\infty} \eee_x \bigl[ 1_{\{i\le\chi\}} \bigl(f(S_{i-1})\bigr)^r
(a_i)^r A_{\chi, f}(S_i, k-r) \bigr]
\nonumber
\\[-8pt]
\\[-8pt]
\nonumber
&=& B_{\chi}(x, k) + C_{\chi}(x,k),
\end{eqnarray}
with
\begin{eqnarray*}
B_{\chi}(x,k) &:=& \sum_{r=1}^{k} C_k^r \sum_{i=1}^{\infty} \eee_x
\bigl[ 1_{\{ i \le\chi, S_i \ge0\}} \bigl(f(S_{i-1})\bigr)^r (a_i)^r A_{\chi,
f}(S_i, k-r) \bigr],
\\
C_{\chi}(x,k) &:=& \sum_{i=1}^{\infty} \eee_x \bigl[ 1_{\{ i \le
\chi
, S_i <0\}} \bigl(f(S_{i-1})\bigr)^k (a_i)^k e^{- \eta S_i} \bigr].
\end{eqnarray*}

In the rest of the proof of the lemma, we shall use twice the notation
$A_\chi(x,k)$, $B_\chi(x,k) $, $C_\chi(x,k)$ but without the subscript
$\chi$ and take $\chi=\tau_0^-$, $f(y)= e^{ b y}$ in the proof of (i)
and $\chi=\min(\tau_0^-, \tau_L^+)$, $f= (L-y+1)^\alpha e^{ b y}$ in
the proof of (ii).

{Proof of (i)}. Let in this proof $A(x,k)= \eee_x [ e^{ -
\eta S_{ \tau^-_0}} ( \sum_{\ell=1}^{ \tau^-_0} e^{ b S_{\ell-1}}
a_\ell)^k ] $. We prove~(\ref{est8sub}) by induction on $k$.

The case $k=0$ follows from (\ref{est1sub}). Let $1\le k < \gamma/b$
and assume that we know that $A(x,j)\le c_j e^{ j b x}$ for all
$0\le j \le k-1$ and $x\ge0$. We have to show that $A(x, k) \le c_k
e^{k b x}$.

Using the induction hypothesis, $A(S_{\ell}, k-r) \le c_{k-r} e^{
(k-r) b S_{\ell}} $ if $S_\ell\ge0$. From~(\ref{axkxk}), we have
\begin{eqnarray*} B(x,k) &\le&c \sum_{r=1}^{k} \sum_{ \ell\ge1 }
\eee
_x \bigl[ e^{ k b S_{\ell-1}} (a_{\ell})^r e^{ (k-r) b \Delta
S_\ell}, \ell\le\tau_0^- \bigr] \\
&\le& c \sum_{r=1}^{k} \sum_{ \ell\ge1 } \eee_x\bigl [ e^{ k b
S_{\ell-1}} (a_{\ell})^r e^{ (k-r) b \Delta S_\ell} \bigr],
\end{eqnarray*}
with $ \Delta S_\ell:=S_\ell- S_{\ell-1}$ for $\ell\ge1$.
By the independence of $(a_\ell, \Delta S_\ell)$, we get that
\begin{eqnarray*} B(x,k) & \le& c \sum_{r=1}^{k} \eee_x\bigl[ (a_{1})^r
e^{ (k-r) b \Delta S_1} \bigr] \sum_{ \ell\ge1 } \eee_x
\bigl[ e^{ k b S_{\ell-1}} \bigr] \\
&=& c e^{ k b x} \sum_{r=1}^{k} \eee\bigl[ (a_{1})^r e^{ (k-r)
b S_1} \bigr] \sum_{ \ell\ge1 } \bigl( \eee\bigl[ e^{ k b
S_1} \bigr] \bigr)^{\ell-1}.
\end{eqnarray*}
Observe that
%
\begin{equation}
\label{eqrefereeasked}\qquad
\sum_{r=1}^{k} \eee\bigl[ (a_{1})^r e^{ (k-r) b S_1} \bigr]
\le\eee\bigl[ \bigl( a_1+ e^{ b S_1}\bigr)^k \bigr] \le2^k \bigl( \eee
\bigl[ a_1^k\bigr] + \eee\bigl[ e^{ k b S_1} \bigr] \bigr) <
\infty,
\end{equation}
and $\eee[ e^{ k b S_1} ] < 1$ since $k <
\gamma/b$. Hence $ B(x,k) \le c_k e^{ k b x}$.

It remains to deal with $C(x,k)$. Observe from (\ref{axkxk}) that
\begin{eqnarray*} C(x,k) &= & \sum_{i=1}^\infty\eee_x\bigl [ e^{b k
S_{i-1} } (a_i)^k 1_{\{\tau_0^- > i-1\}} 1_{\{S_i<0\}} e^{- \eta
S_i}\bigr]
\\ &= & \sum_{i=1}^\infty\eee_x \bigl[ e^{b k S_{i-1} } 1_{\{ \tau
_0^- > i-1\}} \eee_{S_{i-1}} \bigl[ 1_{\{S_1<0\}} (a_1)^k e^{- \eta S_1}\bigr]
\bigr],
\end{eqnarray*}
by Markov's property at $i-1$. Since $y:= S_{i-1} >0$,
\[
\eee_{y}\bigl[ 1_{\{S_1<0\}} (a_1)^k e^{- \eta S_1}\bigr] = e^{- \eta y} \eee\bigl[
1_{\{S_1<-y\}} (a_1)^k e^{- \eta S_1}\bigr] \le\eee\bigl[ 1_{\{S_1<0\}} (a_1)^k
e^{- \eta S_1}\bigr].
\]

It follows that
\[
C(x,k) \le c \sum_{i=1}^\infty\eee_x \bigl[ e^{b k S_{i-1} } \bigr]
\le c' e^{ b k x},
\]
since $bk < \gamma$. This yields that $A(x,k) = B(x,k) + C(x,k) \le c
e^{ b k x}$ proving~(\ref{est8sub}).

\textit{Proof of} (ii). Write in this proof
\[
A(x, j):=
\eee_x\Biggl [ e^{ - \eta S_{ \tau^-_0}} \Biggl( \sum_{\ell=1}^{
\min(\tau
^-_0, \tau_L^+)} (1+L-S_{\ell-1})^\alpha e^{ b S_{\ell-1}} a_\ell
\Biggr)^j \Biggr], \qquad x\in\r, j\ge0.
\]

We mention that $A(x,0)= e^{ - \eta x}$ if $x<0$ and for
$j\ge1$, $A(x,j)=0$ if $x<0$ or $x>L$.

From (\ref{axkxk}), $A(x,k) = B(x,k) + C(x,k)$ with
%
\begin{eqnarray}\label
{bxkxk}  B(x,k) &=& \sum_{r=1}^{k} C_k^r \sum_{j \ge1} \eee_x
\bigl[ (1+L-S_{j-1})^{\alpha r} e^{ r b S_{j-1} } (a_j)^r
\nonumber
\\[-8pt]
\\[-8pt]
\nonumber
&&\hspace*{63pt}{}\times A(S_j, k
-r) 1_{\{ j < \min(\tau_0^-, \tau_L^+ )\}}\bigr], \\
 C(x,k) &=& \sum_{i=1}^\infty\eee_x\bigl[ (L- S_{i-1} +1)^{\alpha k}
a_i^k e^{ b k S_{i-1}} e^{- \eta S_i} 1_{\{ i = \tau_0^- < \tau
_L^+\}
}\bigr]. \label{cxkxk}
\end{eqnarray}

We now prove (\ref{est9sub}) by induction on $p$, where $p$ equals some
integer $m\ge1$.

First, let $m < \gamma/b$, and assume (\ref{est9sub}) holds for all
$A(x,j)$ with $0\le j \le m-1$.
By (\ref{bxkxk}),
\begin{eqnarray*}
&& B(x, m) \le c \sum_{r=1}^{m} \sum_{j \ge1} \eee_x\bigl [
(1+L-S_{j-1})^{\alpha r} e^{ r b S_{j-1} } (a_j)^r (1+L - S_j)^{
\alpha(m-r)}
\\
&&\hspace*{232pt}{}\times e^{ b (m-r) S_j}, j < \tau_L^+ \bigr].
\end{eqnarray*}
Write as before $\Delta S_j=S_{j}- S_{j-1}$. Notice that for
any $j < \tau_L^+$, $(1+L - S_j)^{ \alpha(m-r)} e^{ b (m-r) \Delta
S_j } \le c + c (1+L - S_{j-1})^{ \alpha(m-r)} e^{ b (m-r) \Delta
S_j } $. By the independence of $(a_j, \Delta S_j)$, it is easy to see
that the above expectation under $\eee_x$ is less than
\[
c \eee\bigl[ a_1^r \bigl(1+ e^{ b (m-r) S_1}\bigr)\bigr] \eee_x \bigl[
\bigl(1+L-S_{j-1}\bigr)^{\alpha m} e^{ mb S_{j-1} }, j < \tau_L^+ \bigr].
\]

Since $ \eee[ a_1^r (1+ e^{ b (m-r) S_1})] <\infty$ by
(\ref
{eqrefereeasked}), this implies that
%
\begin{eqnarray}\label{c1+l-x}
 B(x, m) & \le& c' \sum_{j\ge1} \eee_x \bigl[
(1+L-S_{j-1})^{\alpha m} e^{ m b S_{j-1} }, j < \tau_L^+ \bigr]
\nonumber\\
&=& c' e^{ mb x} \sum_{j\ge1} \eee\bigl[ (1+L-x -
S_{j-1})^{\alpha
m} e^{ m b S_{j-1} }, j < \tau_{L-x}^+ \bigr]\\
&\le& c (1+ L-x)^{\alpha m} e^{ mb x}, \nonumber
\end{eqnarray}
where the last estimate follows from the facts that for $j <
\tau_{L-x}^+ $, $ (1+L-x - S_{j-1})^{\alpha m}\le c (1+L-x)^{\alpha m}
+ c \vert S_{j-1}\vert^{\alpha m}$ and that $\sum_{j\ge1} \eee[
\vert S_{j-1}\vert^{\alpha m} e^{ m b S_{j-1} } ] < \infty$
(since $mb < \gamma$).

By Markov's property at $i-1$,
\begin{eqnarray*}
&&C(x,m)= \sum_{i=1}^\infty\eee_x\bigl[ (L- S_{i-1} +1)^{\alpha m}
e^{ b
m S_{i-1}} \eee_{S_{i-1}}\bigl[ 1_{\{S_1<0\}} a_1^m e^{- \eta S_1}\bigr],\\
&&\hspace*{223pt} i -1
< \tau_0^- < \tau_L^+\bigr].
\end{eqnarray*}

As in the proof of (i), $ \eee_{S_{i-1}}[ 1_{\{S_1<0\}}
a_1^m e^{- \eta S_1}] $ is less than some constant, hence
%
\begin{eqnarray}\label{cxm1} C(x,m)&\le& c \sum_{i=1}^\infty\eee_x\bigl[ (L-
S_{i-1} +1)^{\alpha m} e^{ b m S_{i-1}}, i -1 < \tau_0^- < \tau
_L^+\bigr]
\nonumber
\\[-8pt]
\\[-8pt]
\nonumber
&\le& c' (1+ L-x)^{\alpha m} e^{ mb x}, \nonumber
\end{eqnarray}
by (\ref{c1+l-x}). Therefore, $A(x,m)= B(x,m)+ C(x,m) \le c
(1+ L-x)^{\alpha m} e^{ mb x}$ proving the case $m$.

Consider now the case when $\gamma/b=m$ is an integer. Since $m -r <
\gamma/b$ for any $1\le r \le m$, $B(y, m-r) \le
c_{m-r, \alpha} (1+L-y)^{ \alpha(m-r)} e^{ (m-r) y}$ for $0\le y \le
L$. By~(\ref{bxkxk}),
\begin{eqnarray*} && B(x, m) \le c \sum_{r=1}^{m} \sum_{j \ge1}
\eee_x\bigl [ (1+L-S_{j-1})^{\alpha r} e^{ r b S_{j-1} } (a_j)^r
(1+L-S_j)^{ \alpha(m-r)} \\
&&\hspace*{208pt}{}\times e^{ (m-r) S_j } 1_{\{ j< \min( \tau_0^-,
\tau
_L^+ )\}}\bigr].
\end{eqnarray*}

Repeating the same argument as before, we get that
\[
B(x, m) \le c' \eee_x \Biggl[ \sum_{j =1}^{ \min(\tau_0^-, \tau_L^+)}
(1+L-S_{j-1})^{\alpha m} e^{ m b S_{j-1} } \Biggr] \le c (1+L-x)^{
1+ \alpha m} e^{ mb x},
\]
by (\ref{est7sub}). According to (\ref{cxm1}), we get the
same estimate for $C(x, m)$, which proves the case $m=\gamma/b$.

It remains to deal with the case $ m > \gamma/b$.
Let $m_1:= \lfloor\gamma/b \rfloor+1 $ be the least integer larger
than $\gamma/b$ and assume that $\eee[a_1^{m_1}] <\infty$, $\eee[ e^{
b (m_1-1) S_1}] < \infty$. We check that (\ref{est9sub}) is satisfied
for $m=m_1$: applying (\ref{bxkxk}) and using the already proved
results for $A(x, m_1-r)$ (since $m_1-r \le\gamma/b$), we get that
$B(x, m_1) $ is bounded by
\begin{eqnarray*}
&&c \sum_{r=1}^{m_1 } \sum_{j \ge1} \eee_x \bigl[
(1+L-S_{j-1})^{\alpha
r} e^{ r b S_{j-1} } (a_j)^r (1+L - S_j)^{ 1+\alpha(m_1-r)}\\
&&\hspace*{198pt}{}\times e^{ b
(m_1-r) S_j} 1_{\{ j < \tau_L^+\}} \bigr],
\end{eqnarray*}
(the extra 1 in the power comes from the possibility that
$m_1-1=\gamma/b$). As before, we get that
\[
B(x, m_1) \le c' \sum_{j \ge1} \eee_x \bigl[ (1+L-S_{j-1})^{1+
\alpha m_1} e^{ m_1 b S_{j-1} }, j < \tau_L^+ \bigr] \le c
e^{ \gamma(x-L) + m_1 b L},
\]
by applying (\ref{est7sub}). The same estimate holds for $C(x, m_1)$ by
using (\ref{c1+l-x}). This proves that (\ref{est9sub}) holds for
$m=m_1$. The other $m> m_1$ can be treated by induction on $m$, and by
using the same arguments as before, we omit the details.
\end{pf*}

\subsection{\texorpdfstring{Proofs of Lemmas \protect\ref{Lhut0}, \protect\ref{lemnonkilled},
\protect\ref{lemvarphi-infinity} and \protect\ref{lemcontinuity}}
{Proofs of Lemmas 9, 10, 11 and 12}}\label{pr-abc}

We give in this subsection the proofs of these lemmas used in the proof
of Theorem \ref{thmyaglom-critiq}.

\begin{pf*}{Proof of Lemma \ref{Lhut0}} Write in this proof
%
\begin{eqnarray}\label{upsi} A_{\scriptsize (\ref{upsi})}&:=& \Biggl\{\sum_{k=1}^{\tau_t^+-K}
\sum
_{u \in\I_k} H^{u}(t)>0\Biggr\},
\nonumber
\\[-8pt]
\\[-8pt]
\nonumber
 B_{\scriptsize(\ref{upsi})}&:=& \bigl\{
\beta
_t ( w_{\tau_t^+}) \le\tau_t^+-K\bigr\}.
\end{eqnarray}

Let us first observe that Markov's inequality together with part (i)
of Corollary~\ref{cmany2} imply
%
\begin{equation}
\label{eqconexpectUpsilon}
\Q_x^+ (A_{\scriptsize(\ref{upsi})} | \G_{\infty} )
\le
\sum_{k=1}^{\tau_t^+-K} \sum_{u \in\I_k} \pi\bigl(V(u),t\bigr),
\end{equation}
with
\[
\pi(x,t):=\e_x \bigl[ H(t) \bigr] 1_{\{ x \le t\} } + 1_{\{ x >
t\} }.
\]
Furthermore, part (ii) of Corollary \ref{cmany2} yields for any $x\le
t$
\[
\e_x \bigl[ H(t) \bigr] = R(x) e^{\varrho x} \Q^{+}_x \biggl[ { e^{-
\varrho
S_{\tau_t^+}}\over R(S_{\tau_t^+})} 1_{\{ \tau_t^+< \tau_0^-\}
}\biggr]
\le\frac{R(x)}{R(t)} e^{\varrho x} e^{-\varrho t} \le e^{\varrho(x-t)},
\]
from which we deduce that $\pi(x,t) \le e^{\varrho(x-t)} 1_{\{ x \le
t\} } + 1_{\{ x > t\} } \le e^{\varrho(x-t)}$.
Therefore, we obtain
\[
\Q_x^+ (A_{\scriptsize(\ref{upsi})} | \G_{\infty} )
\le\sum
_{k=0}^{\tau_t^+-K-1} e^{\varrho(S_{k}-t)} \sum_{u \in\I_{k+1}}
e^{\varrho\Delta V(u)}.
\]

On the other hand, by the definition of $\beta_t(w_{\tau_t^+})$ [see
(\ref{blu})],
\[
1_{ B_{\fontsize{6.6}{8.6}\selectfont(\ref{upsi})}} \le\sum_{k=0}^{\tau_t^+-K-1} e^{\varrho
(S_{k}-t)} \bigl( \B(w_{k+1})\bigr)^\varrho.
\]
It follows that
%
\begin{equation} \label{715}
\Q_x^+ (A_{\scriptsize(\ref{upsi})} \cup B_{\scriptsize(\ref{upsi})} | \G
_{\infty} ) \le\sum_{k=0}^{\tau_t^+-K-1} e^{\varrho(S_{k}-t)}
b_{k+1}:= \Upsilon(t),
\end{equation}
with $b_{k+1}:=\sum_{u \in\I_{k+1}} e^{\varrho\Delta V(u)}
+ (\B(w_{k+1}))^\varrho$. Recall that under $\Q_x^+$, $(S_k,
b_k)_{k\ge0}$ is a Markov chain; see Proposition \ref{pchange-proba}.
Fix a $\lambda>0$. Then we claim that the following double limits equal
zero:
%
\begin{equation} \label{qxtK}
\limsup_{K\to\infty}\limsup_{t\to\infty} \Q_x^+\bigl( \exists
k< \tau
_t^+-K\dvtx t-S_k < \lambda, \tau^+_t >K\bigr) =0.
\end{equation}

In fact, let $t$ be large, and observe that
\[
\Q_x^+\bigl( \exists k< \tau_t^+-K\dvtx t-S_k < \lambda, \tau^+_t
>K\bigr) \le\Q_x^+\bigl( \tau^+_{t- \lambda} + K < \tau
_t^+\bigr),
\]
which by Markov's property at $\tau^+_{t- \lambda} $, is less than $
\sup_{ t-\lambda< y < t } \Q_y^+ ( K < \tau_t^+)$. By the
absolute continuity between $\Q^+_y$ and $\Q_y$,
\begin{eqnarray*}
\Q_y^+ \bigl( K < \tau_t^+\bigr) &=& \Q_y\biggl [ 1_{\{K < \tau_t^+
\wedge\tau_0^-\}} {R(S_K ) \over R(y)} \biggr] \le{ R(t) \over R(y)}
\Q_y\bigl( \tau^+_t > K\bigr)\\
& =& { R(t) \over R(y)} \Q\bigl( \tau^+_{t-y} > K\bigr).
\end{eqnarray*}

It follows that
\begin{eqnarray*}
\limsup_{t\to\infty} \Q_x^+\bigl( \exists k< \tau_t^+-K\dvtx t-S_k <
\lambda, \tau^+_t >K\bigr) &\le&\Q\bigl( \tau^+_{\lambda} > K\bigr)
\limsup
_{t\to\infty} { R(t) \over R(t- \lambda)} \\
&=& \Q\bigl( \tau^+_{\lambda}
> K\bigr),
\end{eqnarray*}
which goes to $0$ as $K \to\infty$. This proves (\ref{qxtK}).

Let
\[
E_1(t, K): = \bigl\{ \forall k < \tau_t^+- K\dvtx t- S_k \ge\lambda,
\tau
_t^+ >K\bigr\}.
\]
We have $\Q_x^+( \tau^+_t >K) \to1$ as $t \to\infty$, which in view
of (\ref{qxtK}) yields that
for any small $\varepsilon>0$, there exists some $K_0=K_0(\varepsilon,
\lambda)>0$ such that for all $K \ge K_0$, there exists some $t_0(K,
\varepsilon, \lambda)$ satisfying
%
\begin{equation} \label{eq827} \Q_x^+\bigl( E_1(t, K)^c\bigr) \le\varepsilon\qquad
\forall t \ge t_0.
\end{equation}

We claim that there exists some small $\delta>0$ such that
%
\begin{eqnarray}
 \sup_{z\ge0} \Q_z^+\bigl[b_1 ^\delta\bigr]&<& \infty, \label{zb1} \\
 \limsup_{t \to\infty} \Q_x^+\Biggl[ \sum_{k=0}^{\tau_t^+ -1} e^{
\kappa (S_{k}-t) } \Biggr] &<& \infty\label{zb2}
\end{eqnarray}
for any $\kappa>0$.

Assuming for the moment (\ref{zb1}) and (\ref{zb2}), we prove the lemma
as follows: define
\[
E_2(t, K):= \bigcap_{k=0}^{ \tau^+_t-K-1} \bigl\{ b_{k+1} \le e^{
{(\varrho/2)} (t -S_k)} \bigr\} \cap\bigl\{ \tau_t^+ >K\bigr\}.
\]
By (\ref{715}) and on $E_2(t, K) \cap E_1(t, K)$ which is $\G_\infty
$-measurable,
\[
\Q_x^+ (A_{\scriptsize(\ref{upsi})} \cup B_{\scriptsize(\ref{upsi})} | \G
_{\infty} ) \le\Upsilon(t) \le\sum_{k=0}^{\tau_t^+-K-1}
e^{ {
(\varrho/2)} (S_{k}-t) },
\]
which is less than $ e^{- \varrho \lambda/4} \sum_{k=0}^{\tau
_t^+-K-1} e^{ {( \varrho/4)} (S_{k}-t) }$ since on $E_1(t, K)$, $S_k
-t \le- \lambda$ for $k < \tau_t^+-K$. This with (\ref{eq827}) imply
that for all $t \ge t_0$,
%
\begin{eqnarray} \label{eq830}\qquad&& \Q_x^+ (A_{\scriptsize(\ref{upsi})} \cup B_{\scriptsize(\ref{upsi})}
)
\nonumber
\\[-8pt]
\\[-8pt]
\nonumber
&&\qquad\le \varepsilon+ \Q_x^+ \bigl( E_2(t,K)^c \cap E_1(t, K)\bigr) + e^{-
\varrho \lambda/4} \Q_x^+\Biggl[ \sum_{k=0}^{\tau_t^+ -1}
e^{{(\varrho/4)} (S_{k}-t) } \Biggr].
\end{eqnarray}

On the other hand, fix the constant $ \delta>0$ in (\ref
{zb1}), and we have
\begin{eqnarray*} \Q_x^+ \bigl( E_2(t,K)^c \cap E_1(t, K)\bigr) &\le&
\Q
_x^+ \biggl[1_{E_1(t, K)} \sum_{k< \tau_t^+- K} (b_{k+1})^\delta
e^{ - { (\delta\varrho/2)} ( t-S_k)} \biggr] \\
&\le& e^{ -
\delta\varrho\lambda/4} \Q_x^+ \biggl[1_{E_1(t, K)} \sum_{k<
\tau_t^+- K} (b_{k+1})^\delta e^{ - { (\delta\varrho/4) } (
t-S_k)} \biggr] \\
&\le& e^{ - \delta\varrho\lambda/4} \Q_x^+ \Biggl[ \sum_{k<
\tau
_t^+} (b_{k+1})^\delta e^{ - { (\delta\varrho/4) } ( t-S_k)}
\Biggr].
\end{eqnarray*}

Applying Markov's property at $k$ gives that
\begin{eqnarray*}
\Q_x^+\Biggl[ \sum_{k=0}^{\tau_t^+ -1} e^{{ (\delta\varrho/4) }
(S_{k}-t) } (b_{k+1})^\delta\Biggr]&= & \sum_{k=0}^\infty\Q
_x^+\bigl[
1_{\{ k <\tau_t^+ \}} e^{{ (\delta\varrho/4) } (S_{k}-t) } \Q
^+_{S_k} \bigl(b^\delta_1\bigr)\bigr] \\
&\le& \sup_{z\ge0} \Q_z^+\bigl[b_1^\delta\bigr] \Q_x^+\Biggl[ \sum
_{k=0}^{\tau
_t^+ -1} e^{{ (\delta\varrho/4) } (S_{k}-t) } \Biggr].
\end{eqnarray*}

By (\ref{zb1}) and (\ref{zb2}), we get some constant $c$ independent of
$\lambda$ and $ t$ (the constant $c$ may depend on $x$, $\delta$) such
that $\Q_x^+[ \sum_{k=0}^{\tau_t^+ -1} e^{{( \delta\varrho
/4)
} (S_{k}-t) } (b_{k+1})^\delta] \le c$ and $ \Q_x^+[
\sum
_{k=0}^{\tau_t^+ -1} e^{{(\varrho/4)} (S_{k}-t) } ] \le c$.
Going back to (\ref{eq830}), we obtain that for all $K\ge K_0$,
\[
\limsup_{t\to\infty} \Q_x^+ (A_{\scriptsize(\ref{upsi})} \cup B_{\scriptsize(\ref
{upsi})} ) \le\varepsilon+ c e^{- \delta\varrho\lambda/4}+
c e^{-\varrho \lambda/4}.
\]

Letting $\lambda\to\infty$ and $\varepsilon\to0$, we get that
\[
\limsup_{K\to\infty} \limsup_{t\to\infty} \Q_x^+ (A_{\scriptsize(\ref{upsi})}
\cup B_{\scriptsize(\ref{upsi})} )=0.
\]

It remains to show (\ref{zb1}) and (\ref{zb2}). By (\ref{Q+}),
\begin{eqnarray*} \Q_z^+\bigl[b_1^\delta\bigr] &=& \e_z \biggl[{ e^{- \varrho z}
\over R(z)} \sum_{|u|=1} 1_{\{V(u) \ge0\}} R\bigl(V(u)\bigr) e^{\varrho V(u)}
\biggl( \sum_{v \neq u} e^{\varrho(V(v)-z)} + \B(u)^\varrho
\biggr)^\delta
\biggr]
\\
&=& \e\biggl[{ 1 \over R(z)} \sum_{|u|=1} 1_{\{V(u) \ge-z\}}
R\bigl(V(u) +z\bigr) e^{\varrho V(u)} \biggl(\sum_{v \neq u} e^{\varrho
V(v)} +
\B(u)^\varrho\biggr)^\delta\biggr]
\\
&\le& \cases{\displaystyle
c \e\biggl[ \sum_{|u|=1} \bigl(1+ \bigl|V(u)\bigr|\bigr) e^{\varrho V(u)}\biggl ( \biggl(\sum
_{|v|=1} e^{\varrho V(v)}\biggr)^\delta+ \B(u)^{\delta\varrho}\biggr)
\biggr],\vspace*{2pt}\cr
 \qquad\mbox{(critical case),} \vspace*{2pt}\cr
\vspace*{2pt}\cr
\displaystyle c \e\biggl[ \biggl(\sum_{|u|=1} e^{\varrho V(u)}\biggr)^{1+\delta} +
\biggl(\sum_{|u|=1} e^{\varrho V(u)}\biggr) \B(u)^{\delta\varrho} \biggr],\vspace*{2pt}\cr
 \qquad
\mbox{(subcritical case),}}
\end{eqnarray*}
since $R(z) \sim C_R z$ in the critical case and $R(z) \sim
C_R$ in the subcritical case as $z \to\infty$. If $\delta>0$ is
sufficiently small, the later expectations are finite by (\ref{hypb})
together with (\ref{hypcrit}) and (\ref{hypsub}), respectively, which
yields (\ref{zb1}).

To show (\ref{zb2}), we deduce from the absolute continuity between
$\Q
_x^+$ and $\Q_x$ that
%
\begin{eqnarray}
\Q_x^+\Biggl[ \sum_{k=0}^{\tau_t^+ -1} e^{ \kappa (S_{k}-t) }
\Biggr]
= \sum_{k=0}^\infty\Q_x\biggl[ 1_{\{ k <\tau_t^+ \wedge\tau^-_0
\}}
e^{\kappa(S_{k}-t) } {R(S_k)\over R(x)} \biggr]. \label{zb3}
\end{eqnarray}

Let us distinguish the critical and subcritical cases: In
the critical case, $\Q[S_1]=0$ and $R(z) \sim C_R z$ as $z \to\infty$.
There exists some constant $c$ such that for all $t\ge1$, the RHS of
(\ref{zb3}) is less than
\[
c t \sum_{k=0}^\infty\Q_x\bigl[ 1_{\{ k <\tau_t^+ \wedge\tau^-_0
\}}
e^{\kappa(S_{k}-t) } \bigr] = c t \Q_x\Biggl[ \sum_{k=0}^{\tau_t^+
\wedge\tau^-_0 -1} e^{\kappa(S_{k}-t) } \Biggr].
\]

Applying (\ref{a+1}) with $ L=t$ and $\delta= \kappa$ [this $\delta$
has nothing to do with\break that in~(\ref{zb1})] gives that $\Q_x[
\sum
_{k=0}^{\tau_t^+ \wedge\tau^-_0 -1} e^{\kappa(S_{k}-t) }]
\le
c_5 {x+1\over t}$. Hence\break $ \Q_x^+[ \sum_{k=0}^{\tau_t^+ -1}
e^{\kappa(S_{k}-t) } ] \le c (x+1) $ for all $t \ge1$. This
proves (\ref{zb2}) in the critical case.

In the subcritical case, we note that $\Q[S_1]>0$ and $R(\cdot)$ is
bounded. By (\ref{zb3}), we get that for some constant $c>0$,
\[
\Q_x^+\Biggl[ \sum_{k=0}^{\tau_t^+ -1} e^{\kappa(S_{k}-t) } \Biggr]
\le c \sum_{k=0}^\infty\Q_x\bigl[ 1_{\{ k <\tau_t^+ \}} e^{\kappa
(S_{k}-t) } \bigr],
\]
which, according to Lemma \ref{lpositivewalk} is uniformly bounded by
some constant. This completes the proof of (\ref{zb2}) and hence that
of Lemma \ref{Lhut0}.
\end{pf*}

\begin{pf*}{Proof of Lemma \ref{lemnonkilled}}
Observe that
\[
\bigl\{\tau_t^+ > K\bigr\}\cap\Gamma^c(t, K) \subset\bigcup_{k \in(\tau_t^+-
K, \tau_t^+]} \bigcup_{ u \in\I_k} \bigl\{ \exists v \in\cT
^{(u)} \dvtx|u| \le\tau_0^-(v) < \tau_t^+(v)=|v|\bigr\}.
\]

Recall that $\G_{\cc_t} = \sigma\{ (\Delta V(u), u\in
\I
_k), V(w_k), w_k, \I_k, 1\le k \le\tau_t^+ \}$. For any event $F
\in\G_{\cc_t} $, we deduce from Corollary \ref{cmany2} that
\[
\Q_x^+\bigl( \bigl\{\tau_t^+ > K\bigr\}\cap\Gamma^c(t, K)\bigr) \le\Q_x^+\bigl(
F^c\bigr) + \Q_x^+ \biggl[ 1_F \sum_{k \in(\tau_t^+- K, \tau_t^+]}
\sum
_{u \in\I_k} f\bigl(V(u)\bigr)\biggr],
\]
with $f(y):= \mathbf{P}_{y} ( \exists v \dvtx\tau_0^-(v) < \tau_t^+(v)=
|v|)= \mathbf{P}( \exists v \dvtx\tau_{-y}^-(v) < \tau_{t-y}^+(v)=
|v|)$. [We
mention that $f(y)=0$ if $y>t$.] For any $y \le t$, by the branching
property at $\tau_{-y}^-(v)$, $f(y) \le\sup_{z \le- y} \mathbf
{P}_z (
\exists
u\dvtx\tau^+_{t-y}(u) < \infty) = \mathbf{P}( \exists u\dvtx\tau
^+_{t}(u) <
\infty
):= \eta(t)$ which converges to $0$ since the (nonkilled) branching
random walk $V$ goes to $-\infty$. Therefore,
\[
\Q_x^+\bigl( \bigl\{\tau_t^+ > K\bigr\}\cap\Gamma^c(t, K)\bigr) \le\Q_x^+\bigl(
F^c\bigr) + \eta(t) \Q_x^+ \Biggl[ 1_F \sum_{k \in(\tau_t^+- K,
\tau
_t^+]} \# \I_k \Biggr].
\]

Consider an arbitrary $\varepsilon>0$. By Lemma \ref
{lemlaw-convergence-reverse}(ii), $(S_{\tau_t^+} - S_{\tau_t^+ -i},
1\le i \le K)$ converges in law, and hence there exists some $\lambda=
\lambda(\varepsilon, K) >0$ such that for all large~$t$ (in particular,
$t > 4 \lambda$),
\[
\Q_x^+ (F_1):= \Q_x^+ \biggl( \bigl\{\tau_t^+ >K\bigr\} \cap\bigcap_{k \in
(\tau
_t^+- K, \tau_t^+]} \bigl\{ S_k > t- \lambda, |S_k - S_{k-1} |\le
\lambda\bigr\}\biggr) > 1 - \varepsilon,
\]
with obvious definition of the event $F_1$. Let $C>0$, and define
\[
F_2:= F_1 \cap\bigl\{ \forall k \in\bigl(\tau_t^+- K, \tau_t^+\bigr]\dvtx\#\I_k
\le C \bigr\}.
\]

Hence for all sufficiently large $t$, $\Q_x^+(\tau_t^+\le K) \le
\varepsilon$ and
%
\begin{eqnarray}\label{qxtk2} \Q_x^+\bigl( \Gamma^c(t, K)\bigr)
&\le&
2 \varepsilon+ \Q_x^+\bigl( F_1 \cap F_2^c\bigr) + \eta(t) \Q_x^+ \biggl[ 1_{
F_2} \sum_{k \in(\tau_t^+- K, \tau_t^+]} \# \I_k \biggr]
\nonumber
\\[-8pt]
\\[-8pt]
\nonumber
&\le& 2 \varepsilon+ \Q_x^+\bigl( F_1 \cap F_2^c\bigr) + C K \eta(t),
\end{eqnarray}
with $\eta(t) \to0$ as $t\to\infty$. By (\ref{hypcrit}) and
(\ref{hypsub}), we can find a sufficiently small $\delta>0$ such that
$\Q[(\# \I_1)^\delta]= \e[ (\nu-1)^\delta\sum_{ |u|=1}
e^{\varrho V(u)}]:= c < \infty$. Observe that
\begin{eqnarray*}
\Q_x^+\bigl( F_1 \cap F_2^c\bigr) &\le& C^{- \delta} \Q_x^+\biggl[ 1_{F_1}
\sum
_{k \in(\tau_t^+- K, \tau_t^+]} (\# \I_k )^\delta\biggr] \\
&\le&
C^{-\delta} \sum_{k \ge1} \Q_x^+ \bigl[ 1_{\{ |S_k - S_{k-1}
|\le
\lambda, S_{k-1} > t-\lambda, \tau^+_t \ge k\}} (\# \I_k
)^\delta\bigr]
\\
&=& C^{-\delta} \sum_{k \ge1} \Q_x \biggl[ { R(S_k) \over R(x)}
1_{\{
|S_k - S_{k-1} |\le\lambda, S_{k-1} > t-\lambda, k \le\tau^+_t
\wedge
\tau_0^-\}} (\# \I_k )^\delta\biggr]
\\
& \le& C^{-\delta} \sum_{k \ge1}{ R(t+\lambda) \over R(x)} \Q_x
\bigl[ 1_{\{ S_{k-1} > t-\lambda, k \le\tau^+_t \wedge\tau_0^-\}}
(\# \I_k )^\delta\bigr],
\end{eqnarray*}
since $R$ is nondecreasing and $S_k \le t+\lambda$. By
Corollary \ref{cmany1}(i), under $\Q_x$, $\# \I_k $ is independent of
$\{ S_{k-1} > t-\lambda, k \le\tau^+_t \wedge\tau_0^-\}$ and has the
same law as $\#\I_1$; moreover $ \Q_x[ (\# \I_1)^\delta] = \Q[ (\#
\I
_1)^\delta]=:c < \infty$. Using the fact that $R(t+\lambda) \le2
R(t-\lambda) $ for all large $t$, we have
\begin{eqnarray*}
\Q_x^+\bigl( F_1 \cap F_2^c\bigr) & \le& c C^{ -\delta} \sum_{k \ge1} {
R(t+\lambda) \over R(x)} \Q_x [ 1_{\{ S_{k-1} > t-\lambda, k
\le
\tau^+_t \wedge\tau_0^-\}} ]
\\
&\le& 2 c C^{ -\delta} \sum_{k \ge1} \Q_x \biggl[ { R(S_{k-1} )
\over R(x)} 1_{\{ S_{k-1} > t-\lambda, k \le\tau^+_t \wedge\tau
_0^-\}
} \biggr]
\\
&= & 2 c C^{ -\delta} \Q_x ^+ \Biggl[ \sum_{k =1}^{ \tau^+_t }
1_{\{ S_{k-1} > t-\lambda\}} \Biggr].
\end{eqnarray*}

Observe that $\Q_x ^+ [ \sum_{k =1}^{ \tau^+_t } 1_{\{ S_{k-1} >
t-\lambda\}} ] \le\Q_x ^+ [ \sum_{k =1}^{ \tau^+_t } e^{
\varrho( S_{k-1} - (t-\lambda))} ] $ which\break by~(\ref{zb2}) is
smaller than some constant $ c= c(\lambda, x)< \infty$. Going back to
(\ref{qxtk2}), we get that
\[
\Q_x^+\bigl( \Gamma^c(t, K)\bigr) \le2 \varepsilon+ 2 c
C^{-\delta
}+ C K \eta(t).
\]
Letting $t \to\infty$, $C \to\infty$ and then $\varepsilon\to0$
($\delta$ being fixed), we prove Lemma \ref{lemnonkilled}.
\end{pf*}

\begin{pf*}{Proof of Lemma \ref{lemvarphi-infinity}}
First, note that there is nothing to prove in the subcritical case
[since ${\mathscr R}(t) \equiv1$ by (\ref{def-Rt})]. It remains to
consider the critical case, thus $\varrho=\varrho_*$ and ${\mathscr R}(t)=t$
for all $t\ge0$. For notational convenience, write
\begin{eqnarray*}
A&:=& \exp\Biggl\{ - f(t_0) 1_{\D_{1,K}} - \sum_{i=1}^K 1_{\D_{i,K}}
\sum
_{j=1}^{m^{(i)}} \bigl\langle f, \overline\mu
_{s_i-t_0-x_j^{(i)}}^{(i,j)} \bigr\rangle\Biggr\}, \\
B&:=& e^{ \varrho_* t_0} +\sum_{i=1}^K \sum_{j=1}^{m^{(i)}} \int e^{
\varrho_* z} \overline\mu_{s_i-t_0-x_j^{(i)}}^{(i,j)} (\mathrm
{d}z), \\
D&:=& t_0 e^{\varrho_* t_0} +\sum_{i=1}^K \sum_{j=1}^{m^{(i)}} \int z
e^{ \varrho_* z} \overline\mu_{s_i-t_0-x_j^{(i)}}^{(i,j)}
({d}z).
\end{eqnarray*}
Then
\begin{eqnarray*}
\varphi_{t,K}\bigl(t_0,s_1,\ldots,s_K,\theta^{(1)},\ldots,\theta
^{(K)}\bigr) &=& \e\biggl[ { A \over B+ {(1/ t)} D} \biggr], \\
\varphi_{\infty,K}\bigl(t_0,s_1,\ldots,s_K,\theta^{(1)},\ldots
,\theta
^{(K)}\bigr) &= & \e\biggl[ { A \over B } \biggr].
\end{eqnarray*}
Since $f \ge0$, $A \le1$, and we get that
\begin{eqnarray*}
&&\bigl\vert\varphi_{t,K}\bigl(t_0,s_1,\ldots,s_K,\theta^{(1)},\ldots
,\theta
^{(K)}\bigr) - \varphi_{\infty,K}\bigl(t_0,s_1,\ldots,s_K,\theta
^{(1)},\ldots,\theta^{(K)}\bigr) \bigr\vert\\
&&\qquad
\le{1\over t} \e\biggl[ { D \over B^2}\biggr].
\end{eqnarray*}
We are going to prove that
\[
{ D \over B^2} \le{ 1\over\varrho_*}\qquad \mathrm{a.s.}
\]
Indeed, notice first that the nonkilled branching random walk $V$ goes
to $-\infty$, $\mu_{s_i-t_0-x_j^{(i)}}^{(i,j)} ({d}z) $ is an a.s.
finite measure on $\r_+$, and $t_0 e^{\varrho_* t_0} \le{1\over
\varrho
_*} e^{ 2 \varrho_* t_0} $ for any $t_0> 0$. Second, let $\zeta_{i,
j}:= \sup\{a>0\dvtx\int_{[a, \infty)} \mu_{s_i-t_0-x_j^{(i)}}^{(i,j)}
({d}
z)>0\}$. Note that $ \zeta_{i, j} \le{1\over\varrho_*} e^{\varrho_*
\zeta_{i, j}} \le{1\over\varrho_*} \int e^{ \varrho_* z}
\overline\mu_{s_i-t_0-x_j^{(i)}}^{(i,j)} ({d}z)$. It follows that
$\int z e^{ \varrho_* z} \overline\mu
_{s_i-t_0-x_j^{(i)}}^{(i,j)} ({d}z) \le\zeta_{i, j} \int e^{
\varrho
_* z} \overline\mu_{s_i-t_0-x_j^{(i)}}^{(i,j)} ({d}z) \le
{1\over
\varrho_*} ( \int e^{ \varrho_* z} \overline\mu
_{s_i-t_0-x_j^{(i)}}^{(i,j)} ({d}z))^2$. Hence
\[
D \le{1\over\varrho_*} e^{ 2 \varrho_* t_0} + {1\over\varrho_*}
\sum
_{i=1}^K \sum_{j=1}^{m^{(i)}} \biggl( \int e^{ \varrho_* z}
\overline
\mu_{s_i-t_0-x_j^{(i)}}^{(i,j)} ({d}z)\biggr)^2 \le{B^2\over
\varrho_*}
,
\]
yielding that $\vert\widetilde\varphi_{t,K}
(T_t^+,S^{(t)}_{1},\ldots,S^{(t)}_{K}) - \widetilde\varphi_{\infty
,K}(T_t^+,S^{(t)}_{1},\ldots,S^{(t)}_{K}) \vert\le{1 \over t
\varrho_*}$ and proving Lemma \ref{lemvarphi-infinity}.
\end{pf*}

\begin{pf*}{Proof of Lemma \ref{lemcontinuity}} We prove the
following stronger statement: For any $K \ge1 $,
%
\begin{eqnarray} \label{phikzt}
 && \lim_{t\to\infty} \Q_x^+\bigl[\widetilde\varphi_{\infty
,K}\bigl(T_t^+,S^{(t)}_{1},\ldots,S^{(t)}_{K}\bigr) 1_{\{\tau^+_t >K\}} \bigr]
\nonumber\hspace*{-35pt}
\\[-8pt]
\\[-8pt]
\nonumber
&&\qquad= \Q\biggl[ \frac{\exp\{ - f(U \hat S_{\hat\sigma}) 1_{\D
_{1,K}} -
\sum_{i=1}^K 1_{\D_{i,K}} \sum_{j=1}^{\widetilde\nu_i} \langle f,
\overline\mu_{\hat
S_{i} - U \hat S_{\hat\sigma}-\widetilde X_j^{(i)}}^{(i,j)} \rangle
\}}{ e^{\varrho U \hat S_{\hat\sigma}} +\sum_{i=1}^K
\sum_{j=1}^{\widetilde\nu_i} \int e^{\varrho z} \overline\mu
_{\hat
S_{i} - U \hat S_{\hat\sigma}-\widetilde X_j^{(i)}}^{(i,j)} ({d}z)}
\biggr],\hspace*{-35pt}
\end{eqnarray}
which implies Lemma \ref{lemcontinuity} by letting $K \to
\infty$. Define
\begin{eqnarray*}
\pounds_i(\mathbf{s}, {\bolds{\theta}})&:=& \min_{i\le j \le K} \bigl(
s_j -
\log\B(\theta_j)\bigr),\qquad 1\le i \le K, \\
A(t_0, \mathbf{s}, {\bolds{\theta}})&:=& \exp\Biggl\{ - f(t_0) 1_{\{
\pounds
_1(\mathbf{s}, {\bolds{\theta}}) \ge t_0 \}} - \sum_{i=1}^K 1_{\{
\pounds
_i(\mathbf{s}, {\bolds{\theta}}) \ge t_0\} } \sum_{j=1}^{m^{(i)}}
\bigl\langle f,
\overline\mu_{s_i-t_0-x_j^{(i)}}^{(i,j)} \bigr\rangle\Biggr\}, \\
B(t_0, \mathbf{s}, {\bolds{\theta}})&:=& \ee^{\varrho t_0} +\sum_{i=1}^K
\sum_{j=1}^{m^{(i)}} \int e^{ \varrho z} \overline
\mu_{s_i-t_0-x_j^{(i)}}^{(i,j)} ({d}z)
\end{eqnarray*}
for $\mathbf{s}:= (s_1, \ldots, s_K)$, ${\bolds{\theta}}:= (\theta
_1, \ldots,
\theta
_K)$, with $\theta_i= \sum_{j=1}^{m^{(i)} } \delta_{\{x_j^{(i)}\}}$,
$1\le i \le K$. Denote by $\Theta(\mathbf{s})$ a random variable taking
values in $\Omega_f^{\otimes K}$ with law\break $\prod_{i=1}^K \Xi
_{s_i-s_{i-1}} (d \theta^{(i)}) $. Then (recalling $s_0:=0$)
\begin{eqnarray*}\widetilde
\varphi_{\infty,K}(t_0,\mathbf{s} ) & = &\int \e\biggl[ {
A(t_0, \mathbf{s}, {\bolds{\theta}}) \over B(t_0, \mathbf{s},
{\bolds{\theta}})}
\biggr] \prod_{i=1}^K \Xi_{s_i-s_{i-1}} \bigl(d \theta^{(i)}\bigr) \\
&=& \e\biggl[ { A(t_0, \mathbf{s}, \Theta(\mathbf{s})) \over B(t_0,
\mathbf{s},
\Theta(\mathbf{s}))} \biggr],\qquad (t_0, \mathbf{s}) \in\r_+^* \times\r
_+^K.
\end{eqnarray*}

Plainly the function $\widetilde\varphi_{\infty,K}$ is bounded by
$1$. Therefore Lemma \ref{lemcontinuity} will be a consequence of
Lemma \ref{lemlaw-convergence-reverse}
if we have checked that for any fixed $\mathbf{s} \in\r_+^K$, the
function $t_0 \to\widetilde\varphi_{\infty,K}(t_0,\mathbf{s}
) $ is continuous excepted on a set that is at most countable.


To this end, we study at first the continuity of $y \to\langle f,
\overline\mu_{y}^{(i,j)} \rangle$ which are i.i.d. copies of
$\langle
f, \overline\mu_{y}\rangle$. Recall that $\langle f, \overline\mu
_{y}\rangle=\sum_{u \in\cc_y } f(V(u)-y)$ for any fixed $y>0$. Let us
consider $\widetilde\tau^+_t(u):=\inf\{k\dvtx V(u_k) \ge t \}$ and define
the associated optional line $\widetilde\cc_t$ just like (\ref{defcct}).
By the definition of the stopping line $\widetilde\cc_y$ and the
continuity of $f$, we immediately obtain
%
\begin{eqnarray}
\label{eqlem-killnonkill}
\limsup_{k \to\infty} \bigl\vert\langle f, \overline\mu_{y_k}
\rangle
- \langle f, \overline\mu_{y} \rangle\bigr\vert&\le& f(0) \sum_{u} 1_{\{
\widetilde\tau_y^+= \vert u \vert, V(u)=y\}}
\nonumber
\\[-8pt]
\\[-8pt]
\nonumber
&=& f(0) \sum_{u\in
\widetilde\cc_y} 1_{\{ V(u)=y\}}
\end{eqnarray}
for any sequence $(y_k)_k$, such that $y_k \to y$ when $k \to\infty$.
On the other hand,
Corollary \ref{cmany1}(ii) also holds for this family of optional
lines by replacing $n$ by $\widetilde\tau^+_t$. Then we take the
expectation (under $\mathbf{P}$) in (\ref{eqlem-killnonkill}) and
obtain that
%
\begin{equation}
\label{eqpre-lem-con1}
\e\Bigl[ \limsup_{k \to\infty} \bigl\vert\langle f, \overline
\mu_{y_k} \rangle- \langle f, \overline\mu_{y} \rangle\bigr\vert
\Bigr] \le f(0) e^{-\varrho y} \Q(S_{\widetilde
\tau_y^{+}}=y),
\end{equation}
where $\widetilde\tau_y^{+}:=\inf\{n \ge0 \dvtx S_n \ge y\}$.
Denoting as before by $(H_n)_{n\ge1}$ the
(strict) ascending ladder heights of $S$, we remark that
\[
\Lambda_1:= \bigl\{ y\dvtx\Q(S_{\widetilde\tau_y^{+}}=y) >0 \bigr\}
\subset
\bigcup_{n=1}^\infty\bigl\{ y\dvtx\Q( H_n= y) >0\bigr\} \qquad\mbox{is countable.}
\]
%
%
Then by (\ref{eqpre-lem-con1}), $y \to\langle f,
\overline
\mu_{y} \rangle$ is continuous (in $L^1$ hence a fortiori in
probability) on $y \notin\Lambda_1$. The same holds for $y \to
\langle f, \overline\mu^{(i, j)}_{y} \rangle$ with any $i, j \ge1$.
Now we write explicitly $\Theta(\mathbf{s})$ by a random vector
$\Theta
(\mathbf{s}) = (\theta_1, \ldots, \theta_K)$ with $\theta_i:= \sum
_{j=1}^{M^{(i)} } \delta_{\{X_j^{(i)}\}}$ and the associated random
variables $\pounds_i(\mathbf{s}, {\bolds{\theta}})$, $1\le i \le K$. [The
random variables $M^{(i)} $ take values in $\n$, $X_j^{(i)}$ in $\r$,
and $\pounds_i(\mathbf{s}, {\bolds{\theta}})$ in $\r\cup\{\infty
\}$.]
Observe that all the following three events are countable:
\begin{eqnarray*} \Lambda_2&:=&\bigcup_{i=1}^K \bigl\{ x\dvtx{\mathbf
P}\bigl( x=
X^{(i)}_j, \mbox{ for some } 1\le j \le M^{(i)}\bigr) >0 \bigr\}, \\
\Lambda_3&:=&\bigcup_{i=1}^K \bigl\{ x\dvtx\mathbf{P}\bigl( x= \pounds
_i(\mathbf{s},
{\bolds{\theta}})\bigr) >0 \bigr\}, \\
\Lambda_4&:=& \Lambda_3 \cup \bigcup_{i=1}^K \{s_i - x -
y\dvtx x
\in\Lambda_2, y \in\Lambda_1\}.
\end{eqnarray*}

We claim that $ \varphi_{\infty,K}(t_0,\mathbf{s} )$ is
continuous on $t_0 \notin\Lambda_4$. To check this, we fix $t_0
\notin\Lambda_3$ and take a sequence $t_n \to t_0$ as $n \to\infty$. Let
\[
E:=\bigcup_{i=1}^K \bigcup_{j=1}^{ M_j^{(i)}} \bigl\{ X_j^{(i)} \in s_i
-t_0 - \Lambda_1\bigr\} \cup\bigl\{ \pounds_i(\mathbf{s}, {\bolds{\theta
}})= t_0\bigr\}.
\]

Since $t_0 \notin\Lambda_4$, we deduce from the definition
of $\Lambda_2$ that $\mathbf{P}(E)=0$. Observe that on $E^c$, $ s_i
-t_0 -
X_j^{(i)} \notin\Lambda_1 $ and $t_0 \neq \pounds_i(\mathbf{s},
{\bolds{\theta}})$, hence $
A(t_n, \mathbf{s}, {\bolds{\theta}}) 1_{E^c} \to A(t_0, \mathbf{s},
{\bolds{\theta}
})1_{E^c} $ in probability. In other words, $
A(t_n, \mathbf{s}, {\bolds{\theta}}) \to A(t_0, \mathbf{s}, {\bolds
{\theta}}) $ in
probability, and the same holds for $ B(t_n, \mathbf{s}, {\bolds
{\theta}}) $.
By the dominated convergence theorem, when $n \to\infty$,
\begin{eqnarray*}
\varphi_{\infty,K}(t_n,\mathbf{s} ) &=& \e\biggl[ { A(t_n,
{\mathbf
s}, \Theta(\mathbf{s})) \over B(t_n, \mathbf{s}, \Theta(\mathbf
{s}))} \biggr]
\\
&\to&
\e\biggl[ { A(t_0, \mathbf{s}, \Theta(\mathbf{s})) \over B(t_0, \mathbf{s},
\Theta
(\mathbf{s}))} \biggr] =\varphi_{\infty,K}(t_0,\mathbf{s} ),
\end{eqnarray*}
proving the desired continuity at any $t_0\notin\Lambda
_3$. Then we can apply Lemma \ref{lemlaw-convergence-reverse} and get
Lemma \ref{lemcontinuity}.
\end{pf*}

\subsection{\texorpdfstring{Proof of Lemma \protect\ref{Lconvolution}}{Proof of Lemma 16}}
Throughout the proof, $\delta>0$ is taken to be sufficiently small.

{Proof of (i)}. Let us write $f(x):= - \log\eee e^{ -x \Gamma_1}
$ for $x \ge0$; By a Tauberian theorem,
\[
f(x) \sim a { x \over\log(1/x)},\qquad x \to0.
\]
Let $A_x:= \{ \max_{ 1\le i \le\xi} Y_i \le x^{ - 1 + {\delta/2}}\}$ ($\max_\varnothing=0$). Then for $x>0$,
\[
\ppp\bigl( A_x^c \bigr) \le\eee\sum_{i=1}^\xi x^{ (1+\delta) (1-
{\delta
/2})} Y_i^{ 1+\delta} = c x^{ (1+\delta) (1- {\delta/2})}
=o\bigl( x^{ 1+ \delta/3}\bigr),\qquad x\to0,
\]
since $\delta>0$ is small. By independence of $(\Gamma_i)$,
we have
%
\begin{eqnarray}\label{xyi3} \eee\bigl[ e^{ - x \sum_{i=1}^\xi Y_i \Gamma_i }
\bigr] &=& \eee\bigl[ e^{- \sum_{i=1}^\xi f(x Y_i ) } \bigr]
\nonumber
\\[-8pt]
\\[-8pt]
\nonumber
&=& \eee
\exp
\Biggl[ - \sum_{i=1}^\xi f(x Y_i ) 1_{ A_x} \Biggr] + o\bigl( x^{ 1+
\delta
/3}\bigr).
\end{eqnarray}
Define
\[
\Upsilon_x:= { \log(1/x)\over x} \sum_{i=1}^\xi f(x Y_i ) 1_{
A_x},\qquad 0< x <1.
\]

Plainly as $x \to0$, $\Upsilon_x \to a \sum_{i=1}^\xi Y_i$
almost surely. Notice that on $A_x$, $x Y_i \le x^{\delta/2}$, which
together with the asymptotic properties of $f$ implies that for all $0<
x< x_0$ with $x_0$ sufficiently small, $f(x Y_i) \le2 a { x Y_i \over
\log(1/(x Y_i))} \le{ 4 a \over\delta} { x Y_i \over\log(1/x)}$,
for all $1\le i \le\xi$. Hence
\[
{ \log(1/x)\over x} \bigl( 1- e^{- { x/\log(1/x)}
\Upsilon
_x} \bigr) \le\Upsilon_x \le{ 4 a \over\delta} \sum_{i=1}^\xi Y_i.
\]

By the dominated convergence theorem,
\[
{ \log(1/x)\over x} \Biggl( 1- \eee\exp\Biggl[ - \sum_{i=1}^\xi f(x
Y_i ) 1_{ A_x} \Biggr] \Biggr) \to a \eee\sum_{i=1}^\xi Y_i.
\]
This and (\ref{xyi3}) yield that as $x \to0$, $ { \log
(1/x)\over x} ( 1- \eee[ e^{ - x \sum_{i=1}^\xi Y_i
\Gamma
_i } ] ) \to a \eee\sum_{i=1}^\xi Y_i$ which
implies (i) by a Tauberian theorem.

\textit{Proof of} (ii).
Define $W:= \sum_{i=1}^\xi Y_i $ and let $\lambda>1$ and
$0<\varepsilon<a/2$. By conditioning on $(Y_i)_{1\le i\le\xi}$ and
using the tail of $\Gamma_i$, we have that for large $t$,
\begin{eqnarray*} \ppp\Biggl( \sum_{i=1}^\xi Y_i \Gamma_i > t
\Biggr)
&\ge& \ppp\Bigl(
\max_{1 \le i \le\xi} (Y_i \Gamma_i) > t, W \le\lambda\Bigr)
\\
&\ge& \eee\Biggl[ 1_{\{W \le\lambda\}}\Biggl ( 1- \prod_{i=1}^\xi\biggl( 1-
{ (a-\varepsilon) Y_i^p \over t^p}\biggr) \Biggr)\Biggr] \\& \ge&
(a-2\varepsilon
) \eee\Biggl[1_{\{ W \le\lambda\}} \sum_{i=1}^\xi Y_i^p \Biggr] t^{-p},
\end{eqnarray*}
which implies that
\[
\liminf_{t \to\infty} t^p \ppp\Biggl( \sum_{i=1}^\xi Y_i
\Gamma
_i > t \Biggr) \ge(a-2\varepsilon) \eee\Biggl[1_{\{ W \le\lambda\}}
\sum
_{i=1}^\xi Y_i^p \Biggr].
\]
Letting $\varepsilon\to0$ and then $\lambda\to\infty$ yields the
lower bound.

To prove the upper bound, we remark that by considering $ {c + Y_i\over
c}$ instead of $Y_i$ (with $c>0$), we can assume without loss of
generality that almost surely $Y_i\ge1$ (if $i\le\xi$).

By Markov's inequality ($\delta$ being small),
%
\begin{equation}\label{last11}\ppp\bigl( W > t^{1-\delta/2}\bigr) \le t^{ -
(p+\delta) (1-\delta/2)} \eee\bigl[W^{p+\delta}\bigr] = o\bigl(t^{-p}\bigr).
\end{equation}

Let $\varepsilon>0$ be small, and define
%
\begin{eqnarray}\label{last12}
 A_{\scriptsize(\ref{last12})}&:=& \Bigl\{ \max_{1\le i \le\xi} (Y_i
\Gamma_i) \le\varepsilon t\Bigr\}, \qquad B_{\scriptsize(\ref{last12})}:=\Biggl\{
\sum
_{i=1}^\xi Y_i \Gamma_i \ge t\Biggr\},
\nonumber
\\[-8pt]
\\[-8pt]
\nonumber
C_{\scriptsize(\ref{last12})}&:=&
\bigl\{ W
\le t^{ 1-\delta/2}\bigr\}.
\end{eqnarray}

By conditioning on ${\mathbb Y}:= \sigma\{ Y_i, 1\le i \le
\xi
, \xi\}$, we get that
\[
\ppp( A_{\scriptsize(\ref{last12})}
\cap B_{\scriptsize(\ref{last12})} \cap C_{\scriptsize(\ref{last12})} ) \le t^{ -p -
\delta} \eee\Biggl[ 1_{C_{\fontsize{6.6}{8.6}\selectfont(\ref{last12})}} \eee\Biggl[ \Biggl( \sum
_{i=1}^\xi
Y_i \Gamma_i \Biggr)^{p+\delta}1_{A_{\fontsize{6.6}{8.6}\selectfont(\ref{last12})}} \bigg| {\mathbb
Y}\Biggr]
\Biggr].
\]

By convexity, $( \sum_{i=1}^\xi y_i \Gamma_i )^{p+\delta}
\le
(\sum_{i=1}^\xi y_i)^{p+\delta-1} \sum_{i=1}^\xi y_i \Gamma
_i^{p+\delta
}$ for any $y_i\ge0$. Observe that by using the tail of $\Gamma_i$,
\[
\eee\bigl[ \Gamma_i^{p+\delta} 1_{\{ \Gamma_i \le{ \varepsilon t
/
y_i}\}}
\bigr] \le\int_0^{ \varepsilon t /y_i} (p+\delta) x^{ p+\delta-1}
\ppp(\Gamma_i >x) \,d x \le{ 2 (p+\delta)\over\delta} ( \varepsilon t
/y_i)^\delta,
\]
for all large $t$ and $y_i \le t^{1-\delta/2}$. It follows
that for any $0< \varepsilon<1$,
%
\begin{equation} \label{last22}\quad  \ppp( A_{\scriptsize(\ref{last12})}
\cap B_{\scriptsize(\ref{last12})} \cap C_{\scriptsize(\ref{last12})} ) \le c_{p,
\delta
} t^{-p} \varepsilon^\delta \eee\Biggl[ W^{ p+\delta-1} \sum
_{i=1}^\xi Y_i^{1-\delta}\Biggr].
\end{equation}

Since $Y_i\ge1$, the above expectation is less than $\eee
[W^{p+\delta}]$ which is finite.

Let $1< q <p$ and $p-q < 1/2$. Using Markov's inequality and
conditioning on ${\mathbb Y}$, we obtain
\begin{eqnarray*}
&& \ppp\bigl( \bigl\{\exists i\le\xi\dvtx \varepsilon t < \Gamma_i Y_i <
(1- \varepsilon) t \bigr\} \cap B_{\scriptsize(\ref{last12})} \cap C_{\scriptsize(\ref
{last12})} \bigr)
\\
&&\qquad\le \ppp\biggl( \biggl\{\exists i\le\xi\dvtx \Gamma_i Y_i >\varepsilon t
, \sum_{j\neq i} Y_j \Gamma_j > \varepsilon t \biggr\} \cap C_{\scriptsize(\ref
{last12})} \biggr)
\\ &&\qquad\le (\varepsilon t) ^{-1-q} \eee\Biggl[ \sum_{i=1}^\xi Y_i
\Gamma_i \biggl(\sum_{j\neq i} Y_j \Gamma_j\biggr)^q 1_{ C_{\fontsize{6.6}{8.6}\selectfont(\ref{last12})}}
\Biggr]
\\ &&\qquad\le (\varepsilon t) ^{-1-q} \eee\Biggl[ \sum_{i=1}^\xi Y_i
\biggl(\sum_{k \neq i} Y_k\biggr) ^{q-1} \biggl(\sum_{j\neq i} Y_j \Gamma
_j^{q}\Gamma
_i\biggr) 1_{ C_{\fontsize{6.6}{8.6}\selectfont(\ref{last12})}} \Biggr]
\\ &&\qquad \le (\varepsilon t) ^{-1-q} \eee[\Gamma_1] \eee\bigl[\Gamma
_1^{q}\bigr] \eee\bigl[ W^{ 1+q} 1_{ C_{\fontsize{6.6}{8.6}\selectfont(\ref{last12})}} \bigr],
\end{eqnarray*}
since $(\sum_{j\neq i} Y_j \Gamma_j)^q \le(\sum_{k \neq i} Y_k)
^{q-1} (\sum_{j\neq i} Y_j \Gamma_j^{q})$ for all $i$ by the convexity
inequality and since the $ \Gamma_j$'s are i.i.d. and independent of
${\mathbb Y}$.
Furthermore, observe that $ \eee[ W^{ 1+q} 1_{ C_{\fontsize{6.6}{8.6}\selectfont(\ref{last12})}}
] \le\eee[ W^{ p+\delta}] t^{ (1+q -p -
\delta)
(1-\delta/2)}$. Therefore, we obtain
\[
\ppp\bigl( \bigl\{\exists i\le\xi\dvtx \varepsilon t < \Gamma_i Y_i < (1-
\varepsilon) t\bigr\} \cap B_{\scriptsize(\ref{last12})} \cap C_{\scriptsize(\ref{last12})}
\bigr)
\le c_{\varepsilon, q} t^{ -p - (1+q-p)\delta/2}.
\]

This combined with (\ref{last11}) and (\ref{last22}) yields
that, for all large $t$,
\begin{eqnarray*}
\ppp(B_{\scriptsize(\ref{last12})} )
& \le& \ppp\Bigl( \max_{1 \le i \le\xi} ( Y_i \Gamma_i ) >
(1-\varepsilon) t, C_{\scriptsize(\ref{last12})}
\Bigr) + c'_{p, \delta} t^{-p} \varepsilon^\delta+ o\bigl(t^{-p}\bigr) \\
&\le& \eee\Biggl[ \sum_{i=1}^\xi{ (a+\varepsilon) Y_i ^p \over
(1-\varepsilon)^p t ^p }1_{\{W \le t^{1-\delta/2}\}} \Biggr] + c'_{p,
\delta} t^{-p} \varepsilon^\delta+
o\bigl(t^{-p}\bigr).
\end{eqnarray*}

It follows that
\[
\limsup_{t \to\infty} t^p \ppp\Biggl(
\sum_{i=1}^\xi Y_i \Gamma_i > t \Biggr) \le\eee\Biggl[ \sum
_{i=1}^\xi{ (a+\varepsilon) Y_i ^p \over(1-\varepsilon)^p } \Biggr] +
c'_{p, \delta} \varepsilon^\delta,
\]
where $\delta>0$ is fixed. Letting $\varepsilon\to0$ yields the upper
bound and completes the proof of the lemma.

\section{Notation}\label{SNotation}
\mbox{\\}

\textit{Tree}

$\cT$: genealogical tree;

$\varnothing$: root;

$|u|$: generation of the vertex $u$;

$\nu(u)$: number of children of $u$;

$\F_{n}$: sigma-field of the branching random walk up to generation $n$.


\textit{Branching random walk}

$(V(u), u\in\cT)$: branching random walk;

$\L$: point process on $\mathbb{R}$ governing the positions of the
offspring of an individual;

$\cL[a]$: set of vertices absorbed below level $a$;

$Z[0,L]$: number of vertices in $\cL[0]$ that did not touch level $L$;

$\cc_{t}$: set of vertices absorbed above level $t$;\eject

$\bar\mu_{t}$: point process on $\mathbb{R}$ composed of the
overshoots of the vertices in $\cc_{t}$;

$\tau_{t}^-(u)$, respectively, $\tau_{t}^+(u)$: for a vertex $u$,
hitting time of $(-\infty,t)$, respectively, $(t,\infty)$, on its
ancestral line ($\infty$ if no such time).\vspace*{6pt}


\textit{Killed branching random walk}

$\zz$: set of nonkilled vertices;

$Z$: cardinal of $\zz$;

$\H(t)$: set of nonkilled vertices absorbed above level $t$;

$H(t)$: cardinal of $\H(t)$;

$\mu_{t}$: point process on $\mathbb{R}$ composed of the overshoots of
the vertices in $\H(t)$;

$\widehat\mu_{\infty}$: limit in distribution of $\mu_{t}$
conditioned upon being nonempty.\vspace*{6pt}


\textit{Good and bad vertices}

$\B(u)$: function controlling the jumps of the offspring of $u$;

$\beta_{L}(u)$: gives the first time there is an \textit{atypical} jump.
$\beta_{L}(u)=\infty$ means that vertex $u$ is a \textit{good} vertex;

$\H_{\B}(t)$: set of vertices in $\H(t)$ which are \textit{good};

$\bar\mu_{\B,t}$: the point process $\bar\mu_{t}$ restricted to
\textit{good} vertices;

$\mu_{\B,t}$: the point process $\mu_{t}$ restricted to \textit
{good} vertices;

$Z_{g}[0,L]$: number of \textit{good} vertices in $Z[0,L]$;

$Z_{b}[0,L]$: number of \textit{bad} vertices in $Z[0,L]$.\vspace*{6pt}


\textit{One-dimensional random walk}

$S_{n}$: one-dimensional random walk;

$R(x)$: renewal function of $S_{n}$; see (\ref{Rx});

$\tau_{t}^+$: hitting time of $(t,+\infty)$;

$\tau_{t}^-$: hitting time of $(-\infty,t)$;

$T_{t}^+$: overshoot at level $t$;

$T_{t}^-$: undershoot at level $t$.\vspace*{6pt}

\textit{Spine decomposition}

$w_{n}$: spine at generation $n$;

$\I_{n}$: brothers of $w_{n}$;

$S_{n}$: position of $w_{n}$;

$\G_{n}$: sigma-field generated by $w_{k}, V(w_{k}),\I_{k}$ for
$k\le n$;

$\Q_{x}$: defined by ${{d}\Q_{x} \over d\mathbf
{P}_{x}}_{\mid
\F_{n}}:= \ee
^{-\rho x}\sum_{|u|=n} \ee^{\rho V(u)}$. Under $\Q_{x}$, the spine
is a
centered random walk;

$\Q_{x}^+$: defined by ${{d}\Q_{x}^+ \over d{\mathbf
P}_{x}}_{\mid\F_{n}}:=
{1\over R(x)}\ee^{-\rho x}\sum_{|u|=n} R(V(u))\ee^{\rho
V(u)}{\mathbf
1}_{\{
\tau_{0}^-(u)>|u|\}}$. Under $\Q_{x}^+$, the spine is a centered random
walk conditioned to stay positive;

$\Q_{x}^{(\varrho_{-})}$: defined by ${{d}\Q_{x}^{(\varrho
_{-})} \over
d\mathbf{P}_{x}}_{\mid\F_{n}}:= \ee^{-\varrho_{-} x}\sum_{|u|=n}
\ee
^{\varrho
_{-} V(u)}$. Under $\Q_{x}^{(\varrho_{-})}$, the spine is a random walk
with negative drift.\vspace*{6pt}

\textit{Martingales}

$\partial W_{n}:= - \sum_{|u|=n} V(u)\ee^{\varrho_{*} V(u)} $;

$M_{n}^*:= \sum_{|u|=n} R(V(u))\ee^{\varrho V(u)} $;

$M_{n}^{(\varrho_{-})}:= \sum_{|u|=n} \ee^{\varrho_{-} V(u)}$.

\section*{Acknowledgments}
We thank Pascal Maillard for
helpful discussions and an anonymous referee for her/his careful
reading on the first version of this paper.

%


\printaddresses

\end{document}